\documentclass[10pt]{amsart}
\usepackage{amsfonts}
\usepackage{amsmath}
\usepackage{amsthm}
\usepackage{amssymb}
\usepackage{latexsym}
\usepackage{verbatim}
\newenvironment{aenume}{%
  \begin{enumerate}%
  }{\end{enumerate}}
\usepackage[mathscr]{eucal}

\makeatletter
\@ifclasslater{amsart}{1999/11/24}{}{
\renewcommand*\subjclass[2][1991]{%
  \def\@subjclass{#2}%
  \@ifundefined{subjclassname@#1}{%
    \ClassWarning{\@classname}{Unknown edition (#1) of Mathematics
      Subject Classification; using '1991'.}%
  }{%
    \@xp\let\@xp\subjclassname\csname subjclassname@#1\endcsname
  }%
}
\renewcommand{\subjclassname}{%
  \textup{1991} Mathematics Subject Classification}
\@xp\let\csname subjclassname@1991\endcsname \subjclassname
\@namedef{subjclassname@2000}{%
  \textup{2000} Mathematics Subject Classification}
}
\makeatother

\newtheorem{cor}[equation]{Corollary}
\newtheorem{Corollary}[equation]{Corollary}
\newtheorem{lem}[equation]{Lemma}

\newtheorem{prop}[equation]{Proposition}
\newtheorem{Proposition}[equation]{Proposition}

\newtheorem{Lemma}[equation]{Lemma}
\newtheorem*{Claim}{Claim}
\newtheorem{thm}[equation]{Theorem} 
\newtheorem{Theorem}[equation]{Theorem}
\newenvironment{pf}{\proof}{\endproof}

\theoremstyle{definition}
\newtheorem{defn}[equation]{Definition}
\newtheorem{Definition}[equation]{Definition}
\newtheorem{rem}[equation]{Remark}
\newtheorem{Remark}[equation]{Remark}
\newtheorem{Remarks}[equation]{Remarks}

\theoremstyle{remark}

\numberwithin{equation}{section}
\def\endeq{\end{eqnarray}}
\def\eneq{\end{eqnarray}}
\def\endeqn{\end{eqnarray*}}
\def\eneqn{\end{eqnarray*}}
\def\eq{\begin{eqnarray}}
\def\eqn{\begin{eqnarray*}}
\def\beq{\begin{eqnarray}}
\def\beqn{\begin{eqnarray*}}
\def\endeq{\end{eqnarray}}
\def\endeqn{\end{eqnarray*}}

\def\on{\operatorname}

\def\cl{{\rm{cl}}}
\def\tD{\tilde\Delta}
\def\re{{\rm re}}

\def\tal{\tilde \alpha}
\def\tom{\tilde \omega}
\def\tomp{\tilde \omega^\prime}

\newcommand{\ea}{\end{array}}
\newcommand{\tQ}{{\widetilde Q}}
\newcommand{\tP}{{\widetilde P}}

\newcommand{\thmref}[1]{Theorem~\ref{#1}}
\newcommand{\secref}[1]{Sect.~\ref{#1}}
\newcommand{\subsecref}[1]{\S\ref{#1}}
\newcommand{\lemref}[1]{Lemma~\ref{#1}}
\newcommand{\propref}[1]{Proposition~\ref{#1}}
\newcommand{\corref}[1]{Corollary~\ref{#1}}
\newcommand{\remref}[1]{Remark~\ref{#1}}

\newcommand{\nc}{\newcommand}
\newcommand{\rnc}{\renewcommand}
\nc{\cal}{\mathcal}
\nc{\goth}{\mathfrak}
\rnc{\bold}{\mathbf}
\renewcommand{\frak}{\mathfrak}
\renewcommand{\Bbb}{\mathbb}
\nc{\xn}{X_N^{(r)}}
\nc{\Usb}{U'_q(\fb)}

\nc{\tL}{\tilde L}

\nc{\set}{{\,;\,}}
\nc{\Z}{{\bold Z}}

\nc\Uez{U_q^+(\ag)_\Z}
\nc\Uz{\U(\ag)_\Z}
\nc\atnt{A_{2n}^{(2)}}
\nc\att{A_{2}^{(2)}}

\nc{\seps}{\varepsilon^*}
\nc{\sphi}{\varphi^*}
\nc\wt{{\on{wt\,}}}

\nc{\Cal}{\cal}
\nc{\Xp}[1]{X^+(#1)}
\nc{\Xm}[1]{X^-(#1)}

\nc{\ch}{\mbox{ch}}

\nc{\J}{{\cal J}}
\nc{\C}{{\bold C}}
\nc{\Q}{{\bold Q}}

\nc{\N}{{\bold N}}
\nc{\BZ}{{\Z}}
\nc{\BQ}{{\Bbb Q}}

\nc\lan{\langle}
\nc\ran{\rangle}
\nc\bsl{\backslash}
\nc\mto{\mapsto}
\nc\lra{\leftrightarrow}
\nc\hra{\hookrightarrow}
\nc\sm{\smallmatrix}
\nc\esm{\endsmallmatrix}
\nc\sub{\subset}
\nc\ti{\tilde}
\nc\nl{\newline}
\nc\fra{\frac}
\nc\und{\underline}
\nc\ov{\overline}
\nc\ot{\otimes}
\nc\bbq{\bar{\bq}_l}
\nc\bcc{\thickfracwithdelims[]\thickness0}
\nc\ad{\text{\rm ad}}
\nc\Ad{\text{\rm Ad}}
\nc\Hom{\text{\rm Hom}}
\nc\End{\text{\rm End}}
\nc\Ind{\text{\rm Ind}}
\nc\Res{\text{\rm Res}}
\nc\Ker{\text{\rm Ker}}
\rnc\Im{\text{Im}}
\nc\sgn{\text{\rm sgn}}
\nc\tr{\text{\rm tr}}
\nc\Tr{\text{\rm Tr}}
\nc\supp{\text{\rm supp}}
\nc\card{\text{\rm card}}
\nc\bst{{}^\bigstar\!}
\nc\he{\heartsuit}
\nc\clu{\clubsuit}
\nc\spa{\spadesuit}
\nc\di{\diamond}

\nc\al{\alpha}
\nc\bet{\beta}
\nc\ga{\gamma}
\nc\de{\delta}
\nc\ep{\varepsilon}
\nc\eps{\varepsilon}
\nc\io{\iota}
\nc\om{\omega}
\nc\si{\sigma}
\rnc\th{\theta}
\nc\ka{\kappa}
\nc\lam{\lambda}
\nc\la{\lambda}
\nc\ze{\zeta}

\nc\vp{\varpi}
\nc\vt{\vartheta}
\nc\vr{\varrho}

\nc\Ga{\Gamma}
\nc\De{\Delta}
\nc\Om{\Omega}
\nc\Si{\Sigma}
\nc\Th{\Theta}

\nc\boa{\bold a}
\nc\bob{\bold b}
\nc\boc{\bold c}
\nc\boC{\bold C}
\nc\bod{\bold d}
\nc\boe{\bold e}
\nc\bof{\bold f}
\nc\bog{\bold g}
\nc\boh{\bold h}
\nc\boi{\bold i}
\nc\boj{\bold j}
\nc\bok{\bold k}
\nc\bol{\bold l}
\nc\bom{\bold m}
\nc\bon{\bold n}
\nc\boo{\bold o}
\nc\bop{\bold p}
\nc\boq{\bold q}
\nc\bor{\bold r}
\nc\bos{\bold s}
\nc\bou{\bold u}
\nc\bov{\bold v}
\nc\bow{\bold w}
\nc\boz{\bold z}

\nc\ba{\bold A}
\nc\bb{\bold B}
\nc\bc{\bold C}
\nc\bd{\bold D}
\nc\be{\beta}
\nc\bg{\bold G}
\nc\bh{\bold H}
\nc\bi{\bold I}
\nc\bj{\bold J}
\nc\bk{\bold K}
\nc\bl{\bold L}
\nc\bm{\bold M}
\nc\bn{\bold N}
\nc\bo{\bold O}
\nc\bp{\bold P}
\nc\bq{\bold Q}
\nc\br{\bold R}
\nc\bs{\bold S}
\nc\bt{\bold T}
\nc\bu{\bold U}
\nc\bv{\bold V}
\nc\bw{\bold W}
\nc\bz{\bold Z}
\nc\bx{\bold X}

\nc\ca{\Cal A}
\nc\cb{\Cal B}
\nc\cc{\Cal C}
\nc\cd{\Cal D}
\nc\ce{\Cal E}
\nc\cf{\Cal F}
\nc\cg{\Cal G}
\rnc\ch{\Cal H}
\nc\ci{\Cal I}
\nc\cj{\Cal J}
\nc\ck{\Cal K}
\nc\cm{\Cal M}
\nc\cn{\Cal N}
\nc\co{\Cal O}
\nc\cp{\Cal P}
\nc\cq{\Cal Q}

\nc\cs{\Cal S}
\nc\ct{\Cal T}
\nc\cu{\Cal U}
\nc\cv{\Cal V}
\nc\cz{\Cal Z}
\nc\cx{\Cal X}
\nc\cy{\Cal Y}

\nc\fa{\frak a}
\nc\fb{\frak b}
\nc\fc{\frak c}
\nc\fd{\frak d}
\nc\fe{\frak e}
\nc\ff{\frak f}
\nc\g{\frak g}
\nc\fg{\frak g}
\nc\fh{\frak h}
\nc\h{\mathfrak h}
\nc\fj{\frak j}
\nc\fk{\frak k}
\nc\fl{\frak l}
\nc\fm{\frak m}
\nc\fn{\frak n}
\nc\fo{\frak o}
\nc\fp{\frak p}
\nc\fq{\frak q}
\nc\fr{\frak r}
\nc\fs{\frak s}
\nc\ft{\frak t}
\nc\fu{\frak u}
\nc\fv{\frak v}
\nc\fz{\frak z}
\nc\fx{\frak x}
\nc\fy{\frak y}

\nc\fA{\frak A}
\nc\fB{\frak B}
\nc\fC{\frak C}
\nc\fD{\frak D}
\nc\fE{\frak E}
\nc\fF{\frak F}
\nc\fG{\frak G}
\nc\fH{\frak H}
\nc\fJ{\frak J}
\nc\fK{\frak K}
\nc\fL{\frak L}
\nc\fM{\frak M}
\nc\fN{\frak N}
\nc\fO{\frak O}
\nc\fP{\frak P}
\nc\fQ{\frak Q}
\nc\fR{\frak R}
\nc\fS{\frak S}
\nc\fT{\frak T}
\nc\fU{\frak U}
\nc\fV{\frak V}
\nc\fZ{\frak Z}
\nc\fX{\frak X}
\nc\fY{\frak Y}

\nc\te{\tilde e}
\nc\tf{\tilde f}
\nc\qinti[1]{[#1]_i}
\nc\q[1]{[#1]_q}
\nc\qcoeff[3]{\left[ \begin{smallmatrix} {#1}& \\ {#2}& \end{smallmatrix}
\negthickspace \right]_{#3}}

\nc\isom{\cong} 
\nc{\isoto}{\xrightarrow{\sim}}
\nc{\Sym}{{\goth S}}
\nc{\Lzero}{\Lambda_0}

\nc{\ve}{\varepsilon}
\nc{\Li}{\Lambda_i}
\nc{\vtimesz}{v_\ve \otimes z^m}

\nc{\asltwo}{\widehat{\goth{sl}_2}}
\nc\ag{\widehat{\goth{g}}}  
\nc\ed[2][]{E_{#1, #2 \delta}}

\newcommand{\Um}{\widetilde{\mathbf U}} 
\nc{\tU}{\Um}   

\newcommand{\B}{\mathcal B} 
\newcommand{\La}{{\mathscr L}} 

\newcommand{\aW}{\hat W} 
\newcommand{\eW}{\widetilde W} 

\nc{\U}{{\mathbf U}}
\nc\Up{\U^+}

\nc\Ua{{{}_\ca\mathbf U}}

\newcommand{\Ui}{\Ua}
\newcommand{\Uiv}{{\mathbf U}^{(i)}}

\newcommand{\A}{\mathbf A}
\newcommand{\Ainfty}{\mathbf A_\infty}
\newcommand{\Azero}{\mathbf A_0}
\newcommand{\aR}{\mathscr R}
\nc\car{\mathscr R}
\nc{\aC}{\mathscr C}

\newcommand{\Bm}{\tilde\B} 
\newcommand{\cP}{P_{\cl,+}^{0}}  
\newcommand{\Lm}{\tilde\La} 

\newcommand{\Uml}{\Um[\lambda]}
\newcommand{\Umgl}{\Um[{}^>\lambda]}
\newcommand{\Umgel}{\Um[{}^\ge\lambda]}
\newcommand{\Lml}{\Lm[\lambda]}
\newcommand{\Lmgl}{\Lm[{}^>\lambda]}
\newcommand{\Lmgel}{\Lm[{}^\ge\lambda]}

\newcommand{\Bml}{\Bm[\lambda]}
\newcommand{\Gm}{G(\Bm)}
\newcommand{\Gml}{G(\Bml)}

\newcommand{\hb}{\hat \beta}

\nc\bocind{\N^{\aR_0}(\lambda)}
\nc\bocindp{\N^{\aR_0}(\lambda)'}

\nc\Wa{W_{\rm aff}}   
\nc\utl{\tilde u_\lambda}   
\nc\Vtl{\tilde V(\lambda)}   
\nc\Wtl{\breve V(\lambda)}   
\nc\bmp{\boc_{-_p}}
\nc\bpp{\boc_{+_p}}
\nc\bzp{\boc_{0_p}}

\newcommand{\GL}{{GL}}
\newcommand{\SL}{{SL}}
\newcommand{\Umz}{{}_{\mathcal A}\Um} 
\newcommand{\aff}{{\rm aff}}
\newcommand{\Irr}{\operatorname{Irr}}

\newcommand{\BWzero}{\B_0(\Wtl)}

\begin{document}

\title{Crystal bases and two--sided cells of quantum affine algebras}

\author[J.~Beck]{Jonathan Beck} 
\address{Department of Mathematics,
Bar Ilan University, 52900 Ramat Gan, Israel}
\email{beck@macs.biu.ac.il}

\author[H.~Nakajima]{Hiraku Nakajima}
\address{Department of Mathematics\\
Kyoto University\\
Kyoto 606-8502\\
Japan}
\email{nakajima@kusm.kyoto-u.ac.jp}
\thanks{The work of H.N. is supported by the Grant-in-aid
for Scientific Research (No.13640019), JSPS}

\maketitle

\section{Introduction}

Let $\g$ be an affine Kac-Moody Lie algebra. Let $\U = \U(\g)$ be its
quantum enveloping algebra introduced by Drinfeld and Jimbo, and let
$\Up$ be its positive part. The purpose of this paper is twofold.
First, we define a basis $B$ of $\Up$ which is an analog of a PBW basis
of $\U$ for a finite dimensional simple Lie algebra. It has the
following properties (see \thmref{thm:intbase}):
\begin{enumerate}
\item Each element of $B$ is a product of a monomial in `real root
vectors' and a Schur function in `imaginary root vectors'.
\item The transition matrix between $B$ and Kashiwara--Lusztig's
global crystal basis (or canonical basis) $G(\B(-\infty))$ is
upper-triangular with $1$'s on the diagonal (with respect to a certain
explicitly defined ordering) and with above diagonal entries in
$q_s^{-1}\Z[q_s^{-1}]$.
\end{enumerate}
When $\g$ is symmetric or of type $A_2^{(2)}$, our basis coincides
with the one constructed by Beck-Chari-Pressley \cite{BCP} or Akasaka
\cite{A} respectively, where a property weaker than (2) was
established.

Second, we study the global crystal basis $\B(\Um)$ of the modified
quantum enveloping algebra $\Um$ defined by Lusztig \cite{Lu-mod}. We
obtain a Peter-Weyl like decomposition of the crystal $\B(\Um)$
(\thmref{PWdecomp}), as well as an explicit description of two-sided
cells of $\B(\Um)$ and the limit algebra of $\Um$ at $q=0$
(\thmref{thm:limit}). These results had been conjectured by Nakashima
\cite{Tos}, Kashiwara \cite{Kas00} and Lusztig \cite{Lu-v}
respectively. For type $A^{(1)}_1$, the former was proved in
\cite{Tos}.

Our results are based on the study of ``extremal weight modules''
$V(\lambda)$ for $\lambda\in P$, introduced by Kashiwara \cite{Kas94}.
These $\Um$--modules have global crystal bases $G(\B(\lambda))$, and
for $\lambda \in \bigcup_{w \in W} w(P_+)$ (the Tits cone) are
isomorphic to irreducible highest weight modules [loc.\ cit.]. Outside
of the Tits cone, or equivalently in the affine case for level zero
weights, the structure of these modules has been studied by
Akasaka-Kashiwara \cite{AK} and Kashiwara \cite{Kas00}. In particular,
Kashiwara made conjectures on the crystal $\B(\lambda)$
\cite[\S13]{Kas00}. The present authors independently proved his
conjectures for symmetric $\g$ \cite{Beck3,extrem}. Both proofs used
\cite{BCP}. (The relevance of \cite{BCP} to his conjectures was
already pointed out in \cite{Kas00}.)  In this paper, we first
generalize \cite{BCP} to the non-symmetric cases modulo sign and then
proceed to prove the conjectures on $\B(\lambda)$ modulo sign. Next we
remove the sign ambiguity.  Finally, the above mentioned properties of
$\B(\Um)$ are established.

The construction of the basis $B$ is similar to the previous
construction \cite{BCP,A}, although there are two new ideas worth
remarking on. First, we use the level zero extremal weight modules
$V(\lambda)$ in an essential way to check that after specializing at
$q = \infty$ our basis elements equal canonical basis elements (rather
than equal up to sign, which is easier and is also proved here). Note
that in \cite{BCP} the sign check depends on a positivity result of
Lusztig which is available only for the symmetric case.  Second, we
obtain the fact that our basis $B$ is an integral basis (i.e., a basis
of the Lusztig $\Z[q_s,q_s^{-1}]$--form of $\U^+$) in an interesting
way.  This result follows from the upper triangularity (2) above and
uses the canonical basis and the extremal weight modules in an
essential way.  Our proof is quite different from the proofs
\cite{A,BCP} where the result is obtained by explicitly checking
commutation relations between root vectors.

Furthermore, it should be mentioned that the integrality property of
$B$ also gives (when specializing $q=1$) a construction of an easily
expressed basis of the $\Z$--form of the universal enveloping algebra
of $\g$.  Such a basis has been constructed \cite{Garland, Mitzman},
but these results rely on directly examining all possible commutation
relations between elements of the monomials forming the basis.  Here
we obtain this result as a corollary to the existence of the canonical
basis of $\U^+$.

Let us make one more remark regarding $B$. Recall that Lusztig used a
PBW basis to give an alternative definition of the canonical basis for
finite type $\g$ (an `elementary algebraic definition') \cite{Lu-can}.
Namely, the canonical basis is characterized as 1) an integral basis,
2) invariant under the bar involution, and 3) upper triangular with
respect to the PBW basis.  The existence of such a basis is guaranteed
by the upper triangular property of the bar involution with respect to
the PBW basis. Our definition of elements of $B$ is completely
elementary and
we can prove the upper triangular property without using the global
crystal basis. For symmetric or type $\att$, \cite{A,BCP} proved that
$B$ is a basis of the integral form of $\Up$. The same argument as in
the finite case then gives us the `elementary algebraic definition' of
the global crystal basis. (See \thmref{thm:PBWtoCan}.)
However, when $B$ is not a priori known to be an integral basis (i.e.,
not in the symmetric or $\att$ case), we only show the matrix
expressing the bar involution has entries in $\Q(q_s)$, and we do not
give an alternative algebraic definition of $G(\B(-\infty))$. We hope
that we might avoid the integrality requirement completely in the near
future. In any case, our basis $B$ gives a parametrization of the
global crystal basis and serves us here to prove the above mentioned
conjectures on $\B(\Um)$.

We review the organization of this paper in detail.  In \S2 we
introduce notation and preliminary results from \cite{Kas00},
\cite{extrem} and \cite{Beck3}.  Next, in $\S 3$ we construct the
basis $B$ for the integral form $\Ua^+$ with the properties described
above (up to sign).  In $\S 4$ we consider the crystal structure of
$\Um$ in more detail.  We verify the conjectures of \cite[\S13]{Kas00}
(up to sign) which describe the the crystal structure $\B(\lambda)$ of
$V(\lambda)$. $\B(\lambda)$ decomposes into a product of
$\B_W(\lambda) \times \Irr G_\lambda$, where for $\lambda = \sum_i
\lambda_i \varpi_i$, $\B_W(\lambda)$ denotes the crystal of $W(\lambda) =
\bigotimes_i W(\varpi_i)^{\otimes \lambda_i}$ and $\Irr G_\lambda$ denotes
irreducible representations of $G_\lambda = \prod_i GL_{\lambda_i}(\C)$.
This decomposition is then used to give a $\g \times \g$ bicrystal
decomposition of $\B(\Um) \cong \bigsqcup_{\lambda \in P} \B_0(\lambda)
\times \B(\lambda)/\aW$, where $\B_0(\lambda)$ denotes the connected
component of $\B(\lambda)$ containing the extremal weight vector
$v_\lambda$ and $\aW$ is the affine Weyl group.  In $\S5$ we pause to
remove the sign ambiguity in $\S3$ and $\S4$.

In $\S6$ we study the global basis of the level zero modified quantum
affine algebra $\Um = \bigoplus_{\lambda \in P_\cl^0} \U a_\lambda$.
To each $\lambda \in \cP$ we associate a two sided ideal which is the
intersection of the annihilators of all $V(\lambda')$ for $\lambda'$
outside the cone $\lambda + \cP$ modulo this same ideal further
intersected with the annihilator of $V(\lambda)$.  We show that these
ideals have crystal bases $\Bml$ which have globalizations which
partition the global basis of $\Um$.  We use this partition to
describe the cell structure of $\B(\Um)$ and to verify the conjectures
which appear in \cite{Lu-v}.

Lusztig's conjectures on two-sided cells were based on his conjectures
\cite{Lu-cell} on cells of an affine Hecke algebra, which as far as
the authors know, are still open. In \cite{Lu-cell} Lusztig made a
deep connection between two--sided cells of the affine Hecke algebra
and the geometry of Springer fibers. Our proof is based on extremal
weight modules and is purely algebraic. However geometry is in the
background, since extremal weight modules are isomorphic \cite{extrem}
to universal standard modules, which are defined as $K$-homology
groups of certain quiver varieties introduced by the second author
\cite{Nak00}. For example, values of the $a$-function introduced in
$\S 6$ are equal to the dimensions of the quiver varieties, where the
corresponding result for the affine Hecke algebra was proved in
\cite{Lu-cell}. It is also worthwhile mentioning that the appearance
of $G_\lambda$ is quite natural from quiver varieties.

While the authors were preparing this paper, K. McGerty posted an
article \cite{Mc} to the q-algebra archive where he proves Lusztig's
conjecture for type $A^{(1)}_n$. His proof is completely different
from our proof.

\subsection*{Acknowledgment}
The authors are grateful to Ilaria Damiani who sent us a reprint of
\cite{Da2}.

\section{Preliminaries}

\subsection{Affine Kac--Moody Lie algebras}\label{subsec:affine}

We fix a realization $\g = \g(\xn)$. Here $\xn$ is a diagram from
Table Aff $r$ of \cite[Section 4.8]{Kac}, except in the case of $\xn =
\atnt$ ($n \ge 1$), where we reverse the numbering of the simple
roots.  Let its Cartan subalgebra be $\mathfrak h$.  We denote by $I$
the index set of simple roots. The numbering gives us an
identification $I = \{ 0,1, \dots, n \}$. Let $\{ \alpha_i\}_{i\in I}
\subset \h$ (resp.\ $\{ h_i \}_{i\in I} \subset \h^*$) denote the set
of simple roots (resp.\ simple coroots), where $ \langle h_i,
\alpha_j\rangle = a_{ij},$ where $a_{ij}$ is the Cartan matrix of
$\g$.  Fix $d$ so that $ \langle d, \alpha_j\rangle = \delta_{0j}$.
Denote by $P^* = \bigoplus_{i\in I} \Z h_i \oplus \Z d$ the dual
weight lattice and by $P = \Hom_\Z(P^*, \Z)$ the weight lattice.  Let
$Q = \sum_{i\in I} \Z\alpha_i\subset P$ denote the root lattice,
$\Delta$ the root system and $\Delta^\re = \Delta \setminus \Z\delta$
the set of real roots.  Fix the fundamental weights $\Lambda_i \in P$
defined by $\langle h_i, \Lambda_j \rangle = \delta_{ij}, \langle d,
\Lambda_j\rangle = 0.$ Denote by $Q_+$ the semigroup generated by
positive roots $\sum_{i\in I}\Z_{\ge 0}\alpha_i$; $P_+ = \{ \lambda\in
P\mid \langle h_i, \lambda\rangle \ge 0 \text{ for all } i\in I \}$
the semigroup of integral dominant weights. Let $\Delta^\pm =
\Delta\cap (\pm Q_+)$ be the set of positive and negative roots
respectively.

The center of $\g$ is $1$--dimensional and is spanned by $c =
\sum_{i\in I} a_i^\vee h_i$, where $a_i^\vee$ are the labels of the
dual diagram to $\xn$. $c$ is characterized as the positive
combination of $h_i, i \in I$, for which $\{ h\in P^* \mid \langle h,
\alpha_j \rangle = 0 \text{ for all } j\in I \} = \Z c$.  Let $\delta$
be the unique element $\delta = \sum_{i\in I} a_i \alpha_i $
($a_i\in\Z_{\ge 0}$, where $a_i$ are the numerical labels of $\xn$,
and give a linear dependence between the columns of $a_{ij}$)
satisfying $\{ \lambda \in Q \mid \langle h_i, \lambda\rangle = 0
\text{ for all } i\in I \} = \Z \delta.$ We denote by $h$ the Coxeter
number $\sum_{i\in I} a_i$ and by $h^\vee$ the dual Coxeter number
$\sum_{i\in I} a_i^\vee$.

 Denote the affine Weyl group by $\aW \subset O(\h^*)$ ($=$ the
orthogonal group of $\h^*$ with respect to $(\ ,\ )$) generated by the
simple reflections \( s_i(\lambda) = \lambda - \langle h_i,
\lambda\rangle \alpha_i, \lambda \in \h, i \in I. \) Note that $w
(\delta) = \delta$ for $w \in \aW.$ Denote by $(\ ,\ )$ the
non--degenerate symmetric bilinear form on ${\mathfrak h}^*$ invariant
under the Weyl group action, uniquely characterized by 
\( \langle c, \lambda\rangle = (\delta,\lambda) \), for
$\lambda\in{\mathfrak h}^*$.  Note that $(\alpha_i, \alpha_j) =
a_i^\vee a_i^{-1} a_{ij}$ for $i,j \in I$.

Let $\operatorname{cl}\colon {\mathfrak h}^*\to {\mathfrak
h}^*/\Q\delta$ be the canonical projection.  Let \( {\mathfrak h}^{*0}
= \{ \lambda\in {\mathfrak h}^{*} \mid \langle c,\lambda\rangle =
0\}, \) and define the level zero weight lattice to be \( P^0 =
P\cap{\mathfrak h}^{*0} \).  Let $\h^{*0}_\cl = \cl(\h^{*0})$ and
$P^{0}_\cl = \cl(P^{0})$.  Denote $\Delta_\cl = \cl(\Delta^\re)$.
Since $w(\delta) = \delta$ and $\langle c, \lambda \rangle =
(\delta,\lambda)$, the image of $\aW \subset O(\h^{*0})$ in
$O(\h_\cl^{*0})$ is well-defined and denoted $W_\cl$. Then $W_\cl$ is
the Weyl group of the root system $(\Delta_\cl, \h^{*0}_\cl)$, which 
is reduced, except in type $\atnt$ where it is of type $BC_n$.
The bilinear form $(\ ,\ )$ on $\mathfrak h^*$ descends to a bilinear
form on $\h^{*0}_\cl$, which is denoted also by $(\ ,\ )$. It is
nondegenerate.

We fix $0 \in I$ so that $W_\cl$ is generated by $\{s_i\set i \in I_0
\}$, where $I_0 = I \setminus \{ 0\} = \{1, 2, \dots, n\}$. If $\g$ is
not of type $\atnt$, the choice of $0$ is unique up to a Dynkin
diagram automorphism.  In the case of $\atnt$, there are two choices
of $0$ (either of the two extremal vertices of the Dynkin diagram),
and $(\alpha_{0},\alpha_{0})=1$ or $4$, and accordingly $a_{0}=2$ or
$1$, $a_{0}^\vee=1$ or $2$.  Our choice of $0$ is such that
$(\alpha_{0}, \alpha_{0}) = 4.$ As mentioned above, this is opposite
the numbering convention in \cite{Kac}, but is natural when
constructing $\g = \g(\atnt)$ as a (twisted) loop algebra. We take $\{
\cl(\alpha_i) \mid i\in I_0 \}$ as a set of simple roots of
$\Delta_\cl$, and the corresponding set $\Delta_\cl^+$ of positive
roots.

Let $\alpha^\vee = 2 \alpha/(\alpha,\alpha)$.  Let
$Q^\vee=\sum_{\alpha\in\Delta^\re}\Z\alpha^\vee$.  We set
$Q_\cl=\cl(Q)$, $Q_\cl^\vee=\cl(Q^\vee)$, $\tQ=Q_\cl\cap Q_\cl^\vee$.
We have an exact sequence \beq &&\begin{array}{cccccccccc}
\label{Weylgroupexact} 1&\longrightarrow&\tQ&\xrightarrow{\ t \ }
&\aW&\xrightarrow{{\ \cl \ }}&W_\cl&\longrightarrow&1\,,
\end{array}
\endeq where $t$ is the translation operator given by
\cite[(6.5.2)]{Kac}, and $\cl$ is the above projection $\aW\to W_\cl$.
By abuse of notation we denote $t(\xi)$ simply by $\xi$. For any
$\alpha\in \Delta^\re$, let $\tilde\alpha$ be the element in
$\tQ\cap\Q_{>0}\cl(\alpha)$ with the smallest length.  We set
\[\tD=\{\tilde\alpha\set \alpha\in \Delta^\re\}.\]
Then $\tD$ is a reduced root system, and
$\tQ$ is the root lattice of $\tD$.

An affine Lie algebra $\g$ is either untwisted or the dual of an
untwisted affine Lie algebra or $\atnt$:
\begin{itemize}
\item[(i)] If $\g$ is untwisted, then $2/(\alpha,\alpha)\in\Z$,
  $\tQ=Q_\cl^\vee\subset Q_\cl$, $\tilde\alpha=\cl(\alpha^\vee)$,
  $\tD=\cl(\Delta^{\re\,\vee})$.
\item[(ii)] If $\g$ is the dual of an untwisted algebra, then
  $(\alpha,\alpha)/2\in\Z$, $\tQ=Q_\cl\subset Q_\cl^\vee$,
  $\tilde\alpha=\cl(\alpha)$, $\tD=\cl(\Delta^{\re})$.
\item[(iii)] If $\g=\g(\atnt)$, then $(\alpha,\alpha)/2 = 1/2$, $1$,
  or $2$, $\tQ=Q_\cl=Q_\cl^\vee$, and 
$$\tilde\alpha=
\begin{cases}
\cl(\alpha)&\mbox{ if $(\alpha,\alpha)\not=4$,}\\
\cl(\alpha)/2&\mbox{ if $(\alpha,\alpha)=4$.}
\end{cases}$$
Note that $(\delta - \alpha)/2\in \Delta^\re$ if $(\alpha,\alpha)=4$.
$\tD$ is of type $B_n$.
\end{itemize}
\begin{comment}
For $\atnt$, the first equality of
$\tD=\cl(\Delta^{\re})=\cl(\Delta^{\vee\,\re})$ in \cite{Kas00} is
wrong. $\cl(\Delta^{\re})$ is non-reduced.
\end{comment}

For $\alpha\in\Delta^\re$ or $\alpha\in\Delta_\cl$, we set
\begin{equation} \label{dis}
d_\alpha=\max(1,\dfrac{(\alpha,\alpha)}{2}), 
\end{equation} and
$d_i=d_{\alpha_i}$. We have $m \delta + \alpha \in \Delta^\re \iff
d_\alpha | m$. If $\xn\not=A^{(2)}_{2n}$, then $\tilde\alpha=d_\alpha
\cl(\alpha^\vee)$.

\begin{comment}
  Kashiwara \cite{Kas00} denoted $d_\alpha$ by $c_\alpha$, and defined
  $d_i$ by
\begin{equation*}
    d_i = 
    \begin{cases}
       \max(1,\dfrac{(\alpha,\alpha)}{2}) & \xn\neq\atnt,
\\
       1 & \text{otherwise}.
    \end{cases}
\end{equation*}
But for our choice of $I_0$, $1 = \max(1,
\dfrac{(\alpha_i,\alpha_i)}{2})$ for $i\in I_0$ in $\xn = \atnt$.
\end{comment}

Let $P_\cl^0$ (resp.\ $P_\cl^{0\vee}$) be the dual of $Q_\cl^\vee$
(resp.\ $Q_\cl$), considered as a lattice of $\h^{0*}_\cl$ via $(\ ,\ 
)$. Set $\tP = P_\cl^0\cap P_\cl^{0\vee} = (Q_\cl + Q_\cl^\vee)^*$.
The sets $\Delta$, $\Delta^\vee$ are invariant under the translation
by an element of $\tP$.
We define the {\it extended affine Weyl group\/} by $\eW = \tP \rtimes
W_\cl$. Let $\cal T = \{ w\in \eW \mid w(\Delta^+)\subset \Delta^+\}$.
It is a subgroup of the group of Dynkin diagram automorphisms. We have
$\eW = \cal T\ltimes \aW$. The length function $\ell\colon \aW
\rightarrow \N$ extends to $\ell\colon \eW \rightarrow \N$ where
$\ell(\tau w) = \ell(w)$ for $\tau \in \cal T, w \in \aW.$

Let us denote by $\omega^\vee_i$ ($i\in I_0$) the fundamental
coweights of the root system $(\Delta_\cl, \h_\cl^{*0})$, i.e., 
$(\cl(\alpha_i), \omega^\vee_j) = \delta_{ij}$ for $i,j \in I_0$.  Let
$\tom_i = d_i\omega^\vee_i$. Then $\{\tom_i\}_{i\in I_0}$ is a basis
of $\tP$. We have
\begin{equation}\label{eq:WeylId}
\begin{gathered}
  s_j \tom_i = \tom_i s_j \quad (i\neq j),
   \\
   s_i \tom_i s_i \tom_i^{-1} = t(-d_i\alpha_i^{\vee})
   =
   \begin{cases}
      \tilde\alpha_i^{-1} & \text{if $(\xn,i)\neq(\atnt,n)$},
      \\
      \tilde\alpha_i^{-2} & \text{if $(\xn,i) =(\atnt,n)$}.
   \end{cases}
\end{gathered}
\end{equation}

Set
\begin{equation*}
\begin{gathered}
   \car_> = \{ \alpha\in\Delta^+ \mid \cl(\alpha)\in \Delta_\cl^+ \},
   \quad
   \car_< = \{ \alpha\in\Delta^+ \mid \cl(\alpha)\in -\Delta_\cl^+ \},
\\
   \car_0 = \{ (m\delta,i)\in\Z\delta\times I_0 \mid 
   m > 0, d_i | m \},
\\
   \car = \car_> \sqcup \car_0 \sqcup \car_<.
\end{gathered}
\end{equation*}
These are sets of positive roots counted with multiplicities.


\subsection{Quantum affine algebras}

We define the quantum affine algebra $\U = \U(\g)$ following the
normalization in \cite{AK,Kas00}.  Let $q$ be an indeterminate. For
nonnegative integers $n\ge r$, define
\begin{equation*}
  [n]_q = \frac{q^n - q^{-n}}{q - q^{-1}}, \quad
  [n]_q ! =
  \begin{cases}
   [n]_q [n-1]_q !  &(n > 0),\\
   1 &(n=0),
  \end{cases}
  \quad
  \begin{bmatrix}
  n \\ r
  \end{bmatrix}_q = \frac{[n]_q !}{[r]_q! [n-r]_q!}.
\end{equation*}
We fix the smallest positive
integer $d$ such that $d(\alpha_i,\alpha_i)/2 \in \Z$ for any $i\in
I$. We set $q_s = q^{1/d}$.

Define the quantum affine algebra $\U$ to be the associative algebra
with $1$ over $\Q(q_s)$ generated by elements $E_i$, $F_i$
($i\in I$), $q^h$ ($h\in d^{-1} P^*$),
with defining relations.
{\allowdisplaybreaks
\begin{gather*}
  q^0 = 1, \quad q^h q^{h'} = q^{h+h'},\\
  q^h E_i q^{-h} = q^{\langle h, \alpha_i\rangle} E_i,\quad
  q^h F_i q^{-h} = q^{-\langle h, \alpha_i\rangle} F_i,\\
  E_i F_j - F_j E_i = \delta_{ij}
   \frac{t_i - t_i^{-1}}{q_i - q_i^{-1}},\\
  \sum_{p=0}^{b}(-1)^p E_i^{(p)} E_j E_i^{(b-p)} =
  \sum_{p=0}^{b}(-1)^p F_i^{(p)} F_j F_i^{(b-p)} = 0 \quad
  \text{for $i\ne j$,}
\end{gather*}
where $q_i = q^{(\alpha_i,\alpha_i)/2}$, $t_i =
q^{(\alpha_i,\alpha_i)h_i/2}$,
$b = 1 - \langle h_i,\alpha_j\rangle$,
$E_i^{(p)} = E_i^p/[p]_{q_i}!$,
$F_i^{(p)} = F_i^p/[p]_{q_i}!$.}

Let $\U'$ be the quantized enveloping algebra with
$P_\cl = \cl(P)$ as a weight lattice. It is the
subalgebra of $\U$ generated by $E_i$, $F_i$ ($i\in I$),
$q^h$ ($h\in d^{-1}\bigoplus_i \Z h_i$).

Let $\U^+$ (resp.\ $\U^-$) be the $\Q(q_s)$-subalgebra of $\U$
generated by elements $E_i$'s (resp.\ $F_i$'s).  Let $\U^0$ be the
$\Q(q_s)$-subalgebra generated by elements $q^h$ ($h\in d^{-1}
P^*$). We have the triangular decomposition $\U \cong \U^+\otimes \U^0
\otimes \U^-$.

For $\xi\in  Q$, we define the {\it root space\/} $\U_\xi$ by 
\begin{equation*}
   \U_\xi = \{ x\in \U \mid
     \text{$q^h x q^{-h} = q^{\langle h, \xi\rangle} x$ for all $h\in
      P^*$}\}.
\end{equation*}

Let $\ca = \Z[q_s,q_s^{-1}]$.  Let $\Ui$ be the $\ca$-subalgebra of
$\U$ generated by elements $E_i^{(n)}$, $F_i^{(n)}$, $q^h$ for $i\in
I$, $n\in\Z_{> 0}$, $h\in d^{-1} P^*$.

Let us introduce a $\Q(q_s)$-algebra involutive automorphism $\vee$
and $\Q(q_s)$-algebra involutive anti-automorphisms $*$ and $\psi$ of
$\U$ by
\begin{equation*}
\begin{gathered}
   E_i^\vee = F_i, \quad
   F_i^\vee = E_i, \quad
   (q^h)^\vee = q^{-h},
\\
   E_i^* = E_i, \quad
   F_i^* = F_i, \quad
   (q^h)^* = q^{-h},
\\
   \psi(E_i) = q_i^{-1} t_i^{-1} F_i, \quad
   \psi(F_i) = q_i^{-1} t_i E_i, \quad
   \psi(q^h) = q^h.
\end{gathered}
\end{equation*}
We define a $\Q$-algebra involutive automorphism
$\setbox5=\hbox{A}\overline{\rule{0mm}{\ht5}\hspace*{\wd5}}\,$ of $\U$
by
\begin{equation*}
\begin{gathered}
   \overline{E_i} = E_i, \quad \overline{F_i} = F_i, \quad
   \overline{q^h} = q^{-h},
\\
   \overline{a(q_s)u} = a(q_s^{-1})\overline{u} \quad
   \text{for $a(q_s)\in\Q(q_s)$ and $u\in\U$}.
\end{gathered}
\end{equation*}

We define the coproduct $\Delta$ on $\U$ by
\begin{equation}\label{eq:comul}
\begin{gathered}
   \Delta q^h = q^h \otimes q^h, \quad
   \Delta E_i = E_i\otimes t_i^{-1} + 1 \otimes E_i
\\
   \Delta F_i =  F_i\otimes 1 + t_i \otimes F_i.
\end{gathered}
\end{equation}

\begin{comment}
The coproduct used in \cite{Lu-Book} is equal to
\(
    \tau\circ
    \left(\,\setbox5=\hbox{A}\overline{\rule{0mm}{\ht5}\hspace*{\wd5}}\,
    \otimes
    \setbox5=\hbox{A}\overline{\rule{0mm}{\ht5}\hspace*{\wd5}}\,\right)
    \circ
    \Delta
    \circ
    \,\setbox5=\hbox{A}\overline{\rule{0mm}{\ht5}\hspace*{\wd5}}\,
\),
where $\tau$ is the exchange of factors. We must use this coproduct so 
that the crystal lattice of $V(\lambda)\otimes V(\mu)$ is
$\La(\lambda)\otimes \La(\mu)$. Lusztig's coproduct is OK since he
consider the crystal lattice at $q=\infty$.
\end{comment}

Let us denote by $\varOmega$ the $\Q$-algebra anti-automorphism 
\(
   * \circ \setbox5=\hbox{A}\overline{\rule{0mm}{\ht5}\hspace*{\wd5}}\,
   \circ \vee
\)
of $\U$. We have
\begin{equation*}
   \varOmega(E_i) = F_i, \quad \varOmega(F_i) = E_i, \quad
   \varOmega(q^h) = q^{-h}, \quad \varOmega(q_s) = q_s^{-1}.
\end{equation*}

A $\U$-module $M$ is called {\it integrable\/} if
\begin{enumerate}
\item all $E_i$, $F_i$ ($i\in I$) are locally nilpotent, and
\item it admits a {\it weight space decomposition\/}:
\[
   M = \bigoplus_{\lambda\in P} M_\lambda, \quad
   \text{where $M_\lambda = \{ u\in M\mid \text{$q^h u
     = q^{\langle h,\lambda\rangle} u$ for all $h\in  P^*$}\}$}.
\]
\end{enumerate}
Let $\Um$ be the modified enveloping algebra \cite[Part IV]{Lu-Book}.
It is defined by
\begin{equation*}
    \Um = \bigoplus_{\lambda\in  P} \U a_\lambda,
\quad
    \U a_\lambda = \U \left/
    \sum_{h\in  P^*} \U (q^h - q^{\langle h, \lambda\rangle})
    \right..
\end{equation*}
Here the multiplication is given by
\begin{equation*}
   a_\lambda x = x a_{\lambda-\xi} \quad
   \text{for $\xi\in\U_\xi$},
\qquad
   a_\lambda a_\mu = \delta_{\lambda\mu} a_\lambda,
\end{equation*}
where $a_\lambda$ is considered as the image of $1$ in the above
definition of $\U a_\lambda$.

Let $\lambda$, $\mu\in P_+$. Let $V(\lambda)$ (resp.\
$V(-\mu)$) be the irreducible highest (resp.\ lowest) weight module of 
weight $\lambda$ (resp.\ $-\mu$) \cite[\S3.5]{Lu-Book}. Then there is
a surjective homomorphism
\begin{equation}\label{eq:modify}
   \U a_{\lambda-\mu} \ni u \longmapsto u( u_\lambda\otimes u_{-\mu})
   \in V(\lambda)\otimes V(-\mu),
\end{equation}
where $u_\lambda$ (resp.\ $u_{-\mu}$) is a highest (resp.\ lowest)
weight vector of $V(\lambda)$ (resp.\ $V(-\mu)$).

\subsection{Bilinear Form}

In constructing our crystal base a key component is a variant of a
bilinear form introduced by Drinfeld which characterizes the global
crystal basis of $\Ua^+$.  To introduce the form, first define an
algebra structure on $\U^+ \ot \U^+$ by
\begin{equation*} (x_1\ot x_2)(y_1\ot y_2) 
  =q_s^{(\wt x_2, \wt y_1)}x_1y_1\ot x_2y_2,
\end{equation*}
where $x_t$, $y_t$ ($t=1,2$) are homogeneous. Let $r\colon\U^+ \to\U^+
\ot \U^+$ be the $\Q(q_s)$-algebra homomorphism defined by extending
$r(E_{i}) =E_i \ot 1 + 1\ot E_i$ ($i\in I$). By \cite[1.2.5]{Lu-Book},
the algebra $\Up$ has a unique symmetric bilinear form $(\ , \ )
\colon\Up \times \Up \rightarrow \bq(q_s)$ satisfying $(1,1)=1$ and
\begin{align*}
& (E_i,E_j)=\delta_{ij}\frac{1}{(1-q_i^{-2})}, \\ \ & (x,yy')=(r(x),y\ot
y'), \ \ (xx',y)=(x\ot x',r(y)),
\end{align*} where the form on $\Up \otimes\Up$
is defined by $(x_1\ot y_1, x_2\ot y_2)=(x_1,x_2)(y_1,y_2).$ 

For $i\in I$, introduce the unique $\Q(q_s)$-linear map
$r_i\colon\U^+\to\U^+$ given by $r_i(1)= 0$, $r_i(E_{j})=\delta_{ij}$
for $j\in I$, and satisfying $r_i(xy)=q_s^{(\wt y, \alpha_i)}
r_i(x)y+xr_i(y)$ for all homogeneous $x,y\in\U^+$ ([L3, 1.2.13]).
Similarly, introduce the unique $\Q(q_s)$-linear map
$_ir\colon\U^+\to\U^+$ given by $_ir(1)= 0$, $_ir(E_{j})=\delta_{ij}$,
and satisfying $_ir(xy)= {_ir}(x)y+q_s^{(\wt x, \alpha_i)} x{}
{_ir}(y)$ for all homogeneous $x,y\in\U^+$.

{}From the definition the form satisfies
\begin{align*}{\label{rinner}} 
  (E_{i}y, x) & = (1-q_i^{-2})^{-1}(y, {_ir}(x)), \\ 
  (yE_{i}, x) & = (1-q_i^{-2})^{-1}(y, r_i(x)).
\end{align*}

\subsection{Braid group action}

For each $w\in\aW$, there exists an $\Q(q_s)$-algebra automorphism $T_w$
\cite[\S39]{Lu-Book} (denoted there by $T_{w,1}''$).
Also, for any integrable $\U$-module $M$, there exists $\Q(q)$-linear
map $T_w\colon M\to M$ satisfying $T_w(xu) = T_w(x) T_w(u)$ for $x\in
\U$, $u\in M$ \cite[\S5]{Lu-Book}.
We denote $T_{s_i}$ by $T_i$ hereafter.
By \cite[39.4.5]{Lu-Book} we have
\begin{equation}\label{eq:Omega}
   \varOmega \circ T_w \circ \varOmega = T_w.
\end{equation}

The definition of the automorphism $T_w$ of
$\U$ can be extended to the case $w\in\eW$ by setting
\[
   \tau E_i = E_{\tau(i)}, \quad
   \tau F_i = F_{\tau(i)}, \quad
   \tau q^{h_i} = q^{h_{\tau(i)}}, \quad
   \tau q^d = q^d.
\]

\subsection{Crystal bases}\label{subsec:crystal}

We briefly recall the notion of crystal bases. For the notion of 
(abstract) crystals and more details, we refer to \cite{Kas94,AK}.

For $n\in \Z$ and $i\in I$, let us define an operator acting
on any integrable $\U$-module $M$ by
\begin{gather*}
   \widetilde F_i^{(n)} =
   \sum_{k\ge \max(0,-n)} F_i^{(n+k)} E_i^{(k)} a^n_k(t_i),
\\
   \text{where}\quad
   a^n_k(t_i) = (-1)^k q_i^{k(1-n)} t_i^k
    \prod_{\nu=1}^{k-1} (1 - q_i^{n+2\nu}).
\end{gather*}
And we set $\te_i = F_i^{(-1)}$, $\tf_i = F_i^{(1)}$.

These operators are different from those used for the definition of
crystal bases in \cite{Kas91}, but give us the same crystal bases
by \cite[Proposition~6.1]{Kas00}.

Let
\(
  \A_0 = \{ f(q_s)\in\Q(q_s) \mid \text{$f$ is regular at $q_s=0$}\}.
\)
Let $\Ainfty = \overline{\A_0}$ be the image of $\A_0$ under
$\setbox5=\hbox{A}\overline{\rule{0mm}{\ht5}\hspace*{\wd5}}\,$, that
is, the subring of $\Q(q_s)$ consisting of rational functions regular
at $q_s=\infty$.

\begin{defn}
Let $M$ be an integrable $\U$-module. A pair $(\La,\B)$ is called a
{\it crystal basis\/} of $M$ if it satisfies
\begin{enumerate}
\item $\La$ is a free $\A_0$-submodule of $M$ such that $M\cong
\Q(q_s)\otimes_{\A_0}\La$,
\item $\La = \bigoplus_{\lambda\in P} \La_\lambda$ where
$\La_\lambda = \La\cap M_\lambda$ for $\lambda\in P$,
\item $\B$ is a $\Q$-basis of $\La/q\La \cong \Q\otimes_{\A_0}\La$,
\item $\te_i\La\subset\La$, $\tf_i\La\subset\La$ for all $i\in
I$,
\item if we denote operators on $\La/q\La$ induced by $\te_i$, $\tf_i$ 
by the same symbols, we have
$\te_i\B\subset\B\sqcup\{0\}$, $\tf_i\B\subset\B\sqcup\{0\}$,
\item for any $b,b'\in\B$ and $i\in I$, we have $b' = \tf_i b$ if and
only if $b = \te_i b'$.
\end{enumerate}
\end{defn}

We define functions $\varepsilon_i, \varphi_i\colon\B\to \Z_{\ge 0}$
by
\(
   \varepsilon_i(b) = \max \{ n\ge 0 \mid \te_i^n b\neq 0\},
\)
\(
   \varphi_i(b) = \max \{n \ge 0 \mid \tf_i^n b\neq 0\}.
\)
We set
\(
   \te_i^{\max} b = \te_i^{\varepsilon_i(b)} b,
\)
\(
   \tf_i^{\max} b = \tf_i^{\varphi_i(b)} b.
\)

Let $M$ be an integrable $\U$--module with a crystal basis $(\La,\B)$.
Let $\setbox5=\hbox{A}\overline{\rule{0mm}{\ht5}\hspace*{\wd5}}\,$ be
an involution of an integrable $\U$--module $M$ satisfying \(
\overline{xu} = \overline{x}\,\overline{u} \) for any $x\in\U$, $u\in
M$.
Let ${}_\ca M$ be a $\Ui$--submodule of $M$ such that
$\overline{{}_\ca M} = {{}_\ca M}$, $u - \overline{u} \in (q_s-1){}_\ca
M$ for $u\in {{}_\ca M}$.
We say that {\it $M$ has a global basis\/} $(\La,\B,{{}_\ca M})$ if the
following conditions are satisfied
\begin{enumerate}
\item $M \cong \Q(q_s)\otimes_{\Z[q_s,q_s^{-1}]} {{}_\ca M} \cong
\Q(q_s)\otimes_{\A_0}\La \cong
\Q(q_s)\otimes_{\Ainfty}\overline{\La}$,
\item $\La\cap\overline{\La}\cap {{}_\ca M} \to
  \La\cap {}_\ca M/q_s(\La\cap {}_\ca M)$ is an isomorphism,
\item $\B\subset\La\cap {}_\ca M/q_s(\La\cap {}_\ca M)$.
\end{enumerate}
As a consequence of the definition, natural homomorphisms
\begin{equation*}
\begin{gathered}
  \A_0\otimes_{\Z}\left(\La\cap\overline{\La}\cap {}_\ca M\right) \to \La,
\quad
 \Ainfty\otimes_{\Z}\left(\La\cap\overline{\La}\cap
  {{}_\ca M}\right)\to \overline{\La},
\\
  \Z[q_s,q_s^{-1}]\otimes_{\Z}\left(\La\cap\overline{\La}\cap
  {{}_\ca M}\right) \to {{}_\ca M},
\end{gathered}
\end{equation*}
are isomorphisms.  We call the triple $(\La, \ov{\La}, {{}_\ca M})$ a 
{\it balanced triple\/}. 

\begin{comment}
The previous definition (taken from my earlier paper) was wrong. 
$\La\cap\overline{\La}\cap {{}_\ca M}$ is a $\Z$-module, while
$\La/q_s\La$ is a $\Q$-vector space.
\end{comment}

Let $G$ be the inverse isomorphism $\La\cap {}_\ca M/q_s(\La\cap
{}_\ca M)\to\La\cap\overline{\La}\cap {{}_\ca M}$. Then $\{ G(b) \mid
b\in\B\}$ is a basis of $M$. It is called a {\it global crystal
  basis\/} of $M$. The above conditions imply $\overline{G(b)} =
G(b)$.

$\U^-$ (resp.\ $\U^+$) has a global crystal basis $(\La(\infty),
\B(\infty), \Ua^-)$ (resp.\ $(\La(-\infty), \linebreak[0]\B(-\infty),
\Ua^+)$). (It is not an integrable $\U$-module. But the above
definitions has a modification.)  For a dominant weight $\lambda\in
P_+$, the irreducible highest weight module $V(\lambda)$ has a global
crystal basis $(\La(\lambda), \B(\lambda), {}_\ca V(\lambda))$
\cite{Kas91}. If $\lambda, \mu\in P_+$, then the tensor product
$V(\lambda)\otimes V(-\mu)$ also has a global crystal basis
$(\La(\lambda)\otimes \La(-\mu), \B(\lambda)\otimes\B(-\mu), {}_\ca
V(\lambda)\otimes {}_\ca V(-\mu))$, where the bar involution is
defined by using the quasi $R$-matrix
\cite[Part~IV]{Lu-Book}. Moreover, $\Um$ has a global crystal basis $
( \La(\Um), \B(\Um), {}_\ca \Um)$ such that the homomorphism
\eqref{eq:modify} maps a global basis of $\Um$ to the union of that of
$V(\lambda)\otimes V(-\mu)$ and $0$ (\cite[Theorem~2.1.2]{Kas94} and
\cite[Part~IV]{Lu-Book}).
Note that $\Um$ is {\it not\/} an integrable $\U$-module, and
operators $\te_i$, $\tf_i$ are defined only on
$\La(\Um)/q_s\La(\Um)$. In fact, they are defined so that
$\B(\lambda)\otimes \B(-\mu)\to \B(\U a_{\lambda-\mu}) \subset
\B(\Um)$ is a strict embedding.
Furthermore, the global basis is invariant under $*$
\cite[4.3.2]{Kas94}. The proof given there also shows $\vee$ leaves
the global basis invariant.  We have $\B(\Um) = \bigoplus_{\lambda \in
P} \B(\U a_\lambda),$ the direct sum crystal bases of $\Um a_\lambda.$
We define $\varepsilon^*$, $\varphi^*$, $\te_i^*$, $\tf_i^*$ by
$\varepsilon^*(b) = \varepsilon(b^*)$, $\varphi^*(b) = \varphi(b^*)$,
$\te_i^*(b) = (\te_i(b^*))^*$ and $\tf_i^*(b) = (\tf_i(b^*))^*$. This
another crystal structure on $\Um$ is called the {\it star crystal
structure}. Occasionally we denote $\B(\Um)$ simply by $\Bm$.

\begin{comment}
The $\Um$ is not an integrable $\U$-module, and its global crystal
basis is defined as `limit' of $\lim_{\lambda,\mu\to\infty}
\B(\lambda)\otimes \B(-\mu)$. In particular, the crystal $\B(\Um)$ is
{\it regular}.
\end{comment}

\subsection{Braid group action and global crystal bases}\label{subsec:Lu-bg}

We will recall results of \cite{Lu-bg}. Let
\begin{equation*}
   \U^+[i] = \{ x\in \U^+ \mid T_i(x)\in \U^+ \},
\qquad
   {}^*\U^+[i] = \{ x\in \U^+ \mid T_i^{-1}(x)\in \U^+ \}.
\end{equation*}
By \cite[38.1]{Lu-Book} we have direct sum decompositions of vector spaces
\[
   \U^+ = \U^+[i] \oplus E_i \U^+,
\qquad
   \U^+ = {}^*\U^+[i] \oplus \U^+ E_i.
\]
Let ${}^i\pi\colon \U^+\to \U^+[i]$, $\pi^i\colon \U^+\to {}^*\U^+[i]$
be the natural projections. We have the algebra isomorphism $T_i\colon
\U^+[i]\to {}^*\U^+[i]$. Then \cite[Proposition 1.9]{Lu-bg} says
\begin{equation} \label{eq:prop1.9}
   T_i({}^i\pi(\overline{x})) = \pi^i(\overline{T_i(x)})
\end{equation}
for $x\in \Up[i]$.

\begin{comment}
$T_i$ in \cite{Lu-bg} is $T_i^{-1}$ here.
\end{comment}

By \cite[14.3]{Lu-Book} ${}^i\pi$ (resp.\ $\pi^i$) maps the global
basis to the union of a basis of $\U^+[i]$ (resp.\ ${}^*\U^+[i]$) and
$0$. Then \cite[Theorem~1.2]{Lu-bg} says that $T_i({}^i\pi(G(b)))$, if
it is nonzero, is equal to $\pi^i(G(b'))$ for some $b'$, and the map
$b\mapsto b'$ gives a bijection between \( \{ b \in \B(\infty) \mid
{}^i\pi(G(b)) \neq 0 \} \) and \( \{ b' \in \B(\infty) \mid
\pi^i(G(b')) \neq 0 \}.  \) This result is based in part on an earlier
result \cite{Saito}.

\subsection{Affinization} \label{subsec:aff}
Let $M$ be an integrable $\U'$--module, and let $M=\bigoplus_{\lam\in
P_\cl}M_\lam$ be its weight decomposition.  We define a $\U$--module
$M_\aff$ by
$$M_\aff=\bigoplus_{\lam\in P}M_{\cl(\lam)}.$$ The action of $e_i$ and
$f_i$ are defined by restricting to each summand, so that the
canonical homomorphism $\cl\colon M_\aff\to M$ is $\U'$--linear.  We
define the $\U'$--linear automorphism $z$ of $M_\aff$ with weight
$\delta$ by $(M_\aff)_\lam\isoto
M_{\cl(\lam)}=M_{\cl(\lam+\delta)}\isoto (M_\aff)_{\lam+\delta}$.

Choose a section
$s\colon P_\cl\to P$ of $\cl\colon P\to P_\cl$
such that $s(\cl(\alpha_i))=\alpha_i$ for any $i\in I_0 = I\setminus\{0\}$.
Then $M$ is embedded into $M_\aff$ by $s$ as a vector space.
We have an isomorphism of $\U'$-modules
\begin{equation} \label{emb:aff}
M_\aff\simeq \Q(q_s)[z,z^{-1}]\otimes M.
\end{equation}
Here and $e_i\in\U'$  and $f_i\in\U'$  act on  the right
hand side by $z^{\delta_{i0}}\otimes e_i$ and $z^{-\delta_{i0}}\otimes
f_i$.

Similarly, for a crystal with weights in $P_\cl$,
we can define its affinization
$B_\aff$ by
\beq
&&B_\aff=\bigsqcup_{\lam\in P}B_{\cl(\lam)}.
\eneq
If an integrable $\U'$-module $M$ has a crystal basis 
$(L,B)$, then its affinization
$M_\aff$ has a crystal basis $(L_\aff,B_\aff)$.

For $a\in \Q(q_s)$, we define the $\U'$-module $M_a$ by
\beq
&&M_a=M_\aff/(z-a)M_\aff.
\eneq

\subsection{Extremal weight modules}
\label{subsec:extremal}

A crystal $\B$ over $\U$ is called {\it regular\/} if, for any
$J\subsetneqq I$, $\B$ is isomorphic (as a crystal over
$\U(\g_J)$) to the crystal associated with an integrable
$\U(\g_J)$-module. (It was called {\it normal\/} in \cite{Kas94}.)
Here $\U(\g_J)$ is the subalgebra generated by $E_j$, $F_j$ ($j\in
J$), $q^h$ ($h\in d^{-1}P^*$).
By \cite{Kas94}, the affine Weyl group $\aW$ acts on any regular
crystal. The action $S$ is given by
\begin{equation*}
   S_{s_i} b =
   \begin{cases}
     \tf_i^{\langle h_i, \wt b\rangle} b
       & \text{if $\langle h_i, \wt b\rangle \ge 0$}, \\
     \te_i^{-\langle h_i, \wt b\rangle} b
       & \text{if $\langle h_i, \wt b\rangle \le 0$}
   \end{cases}
\end{equation*}
for the simple reflection $s_i$. We denote $S_{s_i}$ by $S_i$ hereafter.

\begin{defn}\label{def:extremal}
Let $M$ be an integrable $\U$-module. A vector $u\in M$ with weight
$\lambda\in P$ is called {\it extremal\/}, if the following holds for
all $w\in \aW$:
\begin{equation}\label{eq:extremal}
  \begin{cases}
     E_i T_w u = 0 & \text{if $\langle h_i, w\lambda\rangle \ge 0$},
\\
     F_i T_w u = 0 & \text{if $\langle h_i, w\lambda\rangle \le 0$}.
  \end{cases}
\end{equation}
In this case, we define $S_w u$ so that
\begin{equation*}
    S_i S_w u = 
    \begin{cases}
       F_i^{\left(\langle h_i, w\lambda\rangle\right)} S_w u
       & \text{if $\langle h_i, w\lambda\rangle \ge 0$},
       \\
       E_i^{\left(-\langle h_i, w\lambda\rangle\right)} S_w u
       & \text{if $\langle h_i, w\lambda\rangle \le 0$}.
    \end{cases}
\end{equation*}
This is well-defined, i.e., $S_w u$ depends only on $w$.

Similarly, for a vector $b$ of a regular crystal $B$ with weight
$\lambda$, we say that $b$ is {\it extremal\/} if it satisfies
\begin{equation*}
  \begin{cases}
     \te_i S_w b = 0 & \text{if $\langle h_i, w\lambda\rangle \ge 0$},
\\
     \tf_i S_w b = 0 & \text{if $\langle h_i, w\lambda\rangle \le 0$}.
  \end{cases}
\end{equation*}
\end{defn}

\begin{Lemma}[\protect{\cite[Lemma 2.11]{extrem}}] \label{lem:WeylExtremal}
Suppose that an integrable $\U$-module $M$ has a crystal basis
$(\La,\B)$. If $u\in \La\subset M$ is an extremal vector of weight
$\lambda$ satisfying $b= u\bmod q\La \in \B$, then $b$ is an
extremal vector, and we have
\begin{equation*}
  S_w u = (-1)^{N^\vee_+} q^{-N_+} T_w u, \quad
  S_w b = S_w u \bmod q\La \quad\text{for all $w\in \aW$},
\end{equation*}
where
\(
\displaystyle
  N_+ =
  \sum_{\alpha\in \aR_+\cap w^{-1}(\aR_-)}
  \max\left(\left(\alpha, \lambda\right), 0\right),
\)
and $N^\vee_+$ is given by replacing $\alpha$ by $\alpha^\vee$.
\end{Lemma}

For $\lambda\in P$, Kashiwara defined the $\U$-module $V(\lambda)$
generated by $u_\lambda$ with the defining relation that $u_\lambda$
is an extremal vector of weight $\lambda$ \cite{Kas94}\footnote{He
denoted it by $V^{\max}(\lambda)$.}. It has a presentation
\begin{equation*}
   V(\lambda) = \U a_\lambda/I_\lambda, \qquad
   I_\lambda 
   = \bigoplus_{b\in \B(\U a_\lambda)\setminus \B(\lambda)} 
  \Q(q) G(b),
\end{equation*}
where
\(
   \B(\lambda) =  
   \{ b\in \B(\U a_\lambda)\mid \text{$b^*$ is extremal}\}.
\)
(Recall that $\B(\Um)$ is regular \subsecref{subsec:crystal}, so extremal
vectors make sense.)
$I_\lambda$ is a left $\U$--module and $V(\lambda)$ has a crystal
base $(\La(\lambda),\B(\lambda))$ together with a $\Ui$-submodule
${}_\ca V(\lambda)$ with a global crystal basis, naturally induced from
those of $\U a_\lambda$.  We have ${}_\ca V(\lambda) = 
\bigoplus_{b \in \B(\lambda)} \ca G(b)\bmod I_\lambda.$
By the construction of $V(\lambda)$,
$$\{ G(b) u_\lambda \mid b \in \B(\Um)\} \setminus \{ 0\} = \{
G(b)\bmod I_\lambda \mid b \in \B(\lambda) \}.$$ 
By abuse of notation $\B(\lambda)$ is considered both as the crystal 
basis of $V(\lambda)$ and as the subset of the crystal basis of $\B(\U
a_\lambda)$.

For any $w\in W$, $u_\lam\mapsto S_{w^{-1}}u_{w\lam}$ gives an
isomorphism of $\U$-modules:
\[V(\lam)\isoto V(w\lam).\]
This isomorphism sends the global basis to the global basis.
Similarly, we have an isomorphism of crystals
\[S^*_w\colon \B(\lam) \isoto \B(w\lam).\]
Here we regard $\B(\lam)$ as a subcrystal of $\B(\Um)$, and $S^*_w$ is the
Weyl group action on $\B(\Um)$ with respect to the star crystal structure.
Since $\B(\lambda) \subset \B(\U a_\lambda) \cong \B(\infty) \otimes
T_\lambda \otimes \B(-\infty)$ where $T_\lambda$ is the crystal with
one element of weight $\lambda,$ we can consider $\B(\lambda)$ as a
subcrystal of $\B(\infty) \otimes T_\lambda \otimes \B(-\infty)$.

If $\lambda$ is dominant or anti--dominant, then $V(\lambda)$ is
isomorphic to the highest weight module or the lowest weight module of
weight $\lambda$, so in this case the notation is consistent.  For
$\lambda \not \in P^0$, $\lambda$ is in the Tits cone $\bigcup_w w(P_+)$ and
so $V(\lambda)$ is isomorphic to a representation with a dominant or 
anti--dominant weight. 

\begin{Theorem}[\protect{\cite[Theorem 5.1, Corollary 5.2]{Kas00}}]\label{thm:convex}
\textup{(i)} For $\lambda\in P^0$, the weight of any extremal vector
of $\B(\lambda)$ is contained in $\cl^{-1}(\cl(\aW\lambda))$.

\textup{(ii)} For any $\lambda\in P$, the weight of any vector of
$\B(\lambda)$ is contained in the convex hull of $W\lambda$.
\end{Theorem}

When $\xn \neq \atnt$, we define the ``fundamental weights of level
zero'' by setting
$$\varpi_i = \Lambda_i - a_i^\vee \Lambda_0 \in P^0, \ i \in I_0.$$
When $\xn = \atnt$ we set
\begin{align*}
& \varpi_n = 2 \Lambda_n - \Lambda_0  \in P^0, \\
& \varpi_i = \Lambda_i - \Lambda_0 \in P^0,\  i = 1, \dots, n-1.
\end{align*}
We have $P_\cl^0 = \cl(P^0) = \bigoplus_{i \in I_0} \Z \cl(\varpi_i)$.
We say that $\lambda$ is a ``basic weight'' if $\cl(\lambda)$ is
$W_\cl$ conjugate to $\cl(\varpi_i)$ for some $i \in I_0$.

In the $\atnt$ case, our choice of fundamental level zero weights is
different than that of \cite{Kas00}.  It is a simple check that both
choices span the same $\Z$--lattice in $\h$.  It follows that the
image of this lattice under $\cl$ is independent of the choice.
Choosing a basis of $P^0_\cl$ amounts to fixing a Weyl chamber of
$(\Delta_\cl, \h_\cl^{*0})$. Since any two Weyl chambers are conjugate
under the Weyl group action, it follows that our choice of $\varpi_i,
i \in I_0$ give the same set of basic weights as that in
\cite{Kas00}.

\begin{comment}
Kashiwara \cite[\S4.2]{Kas00} defined $M_{\rm aff}$ after choosing a
simple root $\alpha_0$ with $a_0 = 1$. (For $\atnt$, $\alpha_0$ is the
long root and coincides with our choice.) On the other hand, he choose
$0^\vee$ with $a_{0^\vee}^\vee = 1$ to define fundamental
representations \cite[\S5.2]{Kas00}. In our notation $0^\vee = n$ for
$\atnt$. Here we use $0$ to define fundamental representations.
\end{comment}

In \cite[\S5]{Kas00}, Kashiwara describes the structure of level zero
fundamental representations corresponding the basic weights $\varpi_i$ 
($i\in I_0$). Let us note
\begin{equation*}
   (\varpi_i, \tal_j) = \delta_{ij} d_i, \qquad
   \{ n\in \Z \mid \varpi_i + n \delta\in \aW\varpi_i \}
   = \Z d_i,
\end{equation*}
where $d_i$ is as in \eqref{dis}. 
We obtain a $\U'$-linear automorphism $z_i$ of $V(\varpi_i)$ of weight
$d_i\delta$, which sends $u_{\varpi_i}$ to $u_{\varpi_i+d_i\delta}$.
We define the ``fundamental level zero representation'' $\U'$--module
$W(\varpi_i)$ by:
\begin{equation*}
W(\varpi_i)=V(\varpi_i)/(z_i-1)V(\varpi_i).
\end{equation*}

\begin{comment}
If $\xn\neq\atnt$, we have
\begin{equation*}
   (\varpi_i, \tal_j)
   = (\Lambda_i - a_i^\vee\Lambda_0, d_j\alpha_j^\vee)
   = d_j \delta_{ij}.
\end{equation*}
If $\xn = \atnt$ and $i\neq n$, we have
\begin{equation*}
   (\varpi_i, \tal_j)
   = (\Lambda_i - \Lambda_0, \alpha_j)
   = \delta_{ij}.
\end{equation*}
Finally if $(\xn,i) = (\atnt,n)$, we have
\begin{equation*}
   (\varpi_n, \tal_j)
   = (2\Lambda_n-\Lambda_0, \alpha_j)
   = 
   \begin{cases}
     (2\Lambda_n, \alpha_j^\vee) & \text{if $j\neq n$}
     \\
     (2\Lambda_n, \frac12 \alpha_n^\vee) & \text{if $j=n$}
   \end{cases}
   = \delta_{jn}
\end{equation*}
\end{comment}

\begin{Theorem}[\protect{\cite[Theorem 5.15]{Kas00}}]\label{thm:fund}
\textup{(i)} $W(\varpi_i)$ is a finite--dimensional irre\-duc\-ible
integrable $\U'$-module.

\textup{(ii)} $W(\varpi_i)$ has a global crystal basis with a simple
crystal, i.e., the weight of any extremal vector of $\B(W(\varpi_i))$
is contained in $W_\cl \cdot\cl(\varpi_i)$ and $\#
\B(W(\varpi_i))_{\cl(\varpi_i)} = 1$.

\textup{(iii)}
For any $\mu\in \wt(V(\varpi_i))$,
$$W(\varpi_i)_{\cl(\mu)}\simeq V(\varpi_i)_\mu.$$

\textup{(iv)} $\dim W(\varpi_i)_{\cl(\varpi_i)}=1$.

\textup{(v)} The weight of any extremal vector of $ W(\varpi_i)$
belongs to $\aW \cdot \cl(\varpi_i)$.

\textup{(vi)} $\wt(W(\varpi_i))$ is the intersection of
$\cl(\varpi_i)+Q_\cl$ and the convex hull of $\aW\cdot \cl(\varpi_i)$.

\textup{(vii)}
$\Q(q_s)[z_i^{1/d_i}]\otimes_{\Q(q_s)[z_i]}V(\varpi_i)\simeq
W(\varpi_i)_{\rm aff}$, where the action of $z_i^{1/d_i}$ on the left
hand side corresponds to the action of $z$ on the right hand side as
defined in \S\ref{subsec:aff}.

\textup{(viii)} $V(\varpi_i)$ is isomorphic to the submodule
$\Q(q_s)[z^{d_i},z^{-d_i}]\otimes W(\varpi_i)$ of $W(\varpi_i)_{\rm
aff}$ as a $\U$-module.  Here we identify $W(\varpi_i)_{\rm aff}$ with
$\Q(q_s)[z,z^{-1}]\otimes W(\varpi_i)$ as in \eqref{emb:aff}.

\textup{(ix)} Any irreducible finite-dimensional integrable
$\U'$-module with $\cl(\varpi_i)$ as an extremal weight is isomorphic
to $W(\varpi_i)_a$ for some $a\in \Q(q_s)\setminus\{0\}$.
\end{Theorem}

\subsection{Bilinear form characterizing $V(\lambda)$}

We recall the following properties of $V(\lambda)$ with respect to a
natural bilinear form.

\begin{Proposition}[\protect{\cite[Proposition 4.1]{extrem}}]\label{prop:bilinear}
The extremal weight module $V(\lambda)$ has a unique bilinear form $(\ 
,\ )$ satisfying
\begin{gather}
   (u_\lambda, G(b)) = 
   \begin{cases}
      1 & \text{if $G(b) = u_\lambda$}, \\
      0 & \text{otherwise}
   \end{cases}\label{eq:bil1}
\\
   (x u, v) = (u, \psi(x)v) \quad\text{for $x\in\Ua$, $u,v\in V(\lambda)$}.
   \label{eq:bil2}
\end{gather}
\end{Proposition}

The following theorem gives a characterization of the global basis of
$V(\lambda)$ with respect to the form:
\begin{Theorem}[\protect{\cite[Theorem 3]{extrem}}] \label{thm:main2}
\textup{(i)} $\{ G(b) \mid b\in \B(\lambda)\}$ is almost orthonormal
for $(\ ,\ )$, that is, $(G(b), G(b')) \equiv \delta_{bb'} \mod
q\Z[q]$.

\textup{(ii)}
\(
   \{ \pm G(b) \mid b\in\B(\lambda)\}
   = \left\{ u\in {}_\ca V(\lambda) \left|\, \overline{u} = u, \;
       (u,u) \equiv 1 \mod q\Z[q] \right\}\right..
\)
\end{Theorem}

For fundamental representations this result is due to \cite{VV2}.

\begin{Remark} By the ensuing discussion in \S3, 
\thmref{thm:main2} holds in the non--symmetric case also. 
\end{Remark}

\section{An integral crystal base for $\Ua^+$} \label{sec:intbase}

In this section we construct an integral basis $B$ of $\Ua^+$ such
that $B \subset \ov{\La(-\infty)}$ and under the natural map $\pi:
\ov{\La(-\infty)} \rightarrow \ov{\La(-\infty)}/q_s^{-1}
\ov{\La(-\infty)}$, $\pi(B) = \B(-\infty)$. (In fact, we construct 
a series of bases parametrized by $p\in\Z$ with these properties.)

\subsection{Root Vectors}
We introduce root vectors in $\Ua^+$. Recall that we chose $\tom_i
\in \tP \vartriangleleft \eW$ for each $i\in I_0$ (see
\subsecref{subsec:affine}). We choose $\tau_i\in\cal T$ so that
$\tom_i\tau_i^{-1}\in \aW$. Choosing a reduced expression for
$\tom_i\tau_i^{-1}$ for each $i\in I_0$ we fix a reduced expression of
$\tom_n \tom_{n-1} \dots \tom_1$:
\begin{equation}\label{eq:tomorder}
    \tom_n \tom_{n-1} \dots \tom_1
    = s_{i_1} s_{i_2} \dots s_{i_N} \tau, \qquad
    (\tau = \tau_n\dots \tau_1).
\end{equation}
We define a doubly infinite sequence 
\begin{equation*}
   \boh =(\ldots, i_{-1}, i_{0}, i_{1},\ldots )
\end{equation*}
by setting $i_{k+N} = \tau(i_k)$ for $k\in\Z$. Note that for any
integers $m<p$, the product $s_{i_m}s_{i_{m+1}}\ldots s_{i_p} \in
\aW$ is a reduced expression. We have
\begin{equation*}
\begin{split} \car_> = \{ \alpha_{i_0}, s_{i_0} (\alpha_{i_{-1}}), 
s_{i_0} s_{i_{-1}} (\alpha_{i_{-2}}), \dots \}, \\
  \car_< = \{ \alpha_{i_1}, s_{i_1} (\alpha_{i_{2}}), 
s_{i_1} s_{i_{2}} (\alpha_{i_3}), \dots \}.
\end{split}
\end{equation*}

\begin{Remark}
  Our definition of the PBW basis will depend on the sequence $\boh$.
  In particular, it depends on the choice of the numbering $I_0 = \{
  1,2,\dots, n\}$. Almost all of the results in \cite{Da2} which we
  will use are independent of the numbering. But when $\xn = \atnt$,
  \cite[Corollary~4.2.6]{Da2} depends on our choice such that $\tom_n$
  corresponding to the short root $\alpha_n$ appears first in
  \eqref{eq:tomorder}. We choose our $\boh$ to agree with that in
  \cite{Da2}. (Our vertex $n$ is labeled by $1$ there.)
\end{Remark}

Set
\begin{align} 
  \beta_{k}&=
\begin{cases}
s_{i_0}s_{i_{-1}}\ldots s_{i_{k+1}}(\alpha_{i_k}), & \quad {\text{if
 $k\le 0$},}\\ s_{i_1}s_{i_2}\ldots s_{i_{k-1}}(\alpha_{i_k}), & \quad
 {\text{if $k>0$}}.
\end{cases}
\end{align} Define a total order on
  $\car$ by setting
\begin{equation} \label{rootorder}
  \beta_0 < \beta_{-1} < \beta_{-2} \dots < \delta^{(1)} < \dots <
  \delta^{(n)} < 2 \delta^{(1)} < \dots < \beta_3 < \beta_2 < \beta_1,
\end{equation}
where $k\delta^{(i)}$ denotes $(k\delta,i)\in\car_0$.

We now define root vectors for each element of $\car_>\sqcup \car_<$.
\begin{align} \label{realroots}
E_{\beta_k}  = 
\begin{cases}    T_{i_0}^{-1} T_{i_{-1}}^{-1} \dots
  T_{i_{k+1}}^{-1} (E_{i_k}) & \text{if $k\le 0$}, \\ T_{i_1}
 T_{i_{2}} \dots T_{i_{k-1}} (E_{i_k}) &\text{if $k>0$}.
\end{cases}
\end{align}
By \cite[40.1.3]{Lu-Book} these are in $\Up$.
As usual, set $F_\beta = \Omega(E_\beta)$ for all $\beta \in
\car_>\sqcup \car_<$.

\begin{Remark} \label{rem:loopgen} We note that the root 
vectors $E_{d_i k \delta \pm \alpha_i}$
are described explicitly by:
$$E_{k d_i \delta + \alpha_i} = T_{\tom_i}^{-k} E_i\ (k \ge 0), E_{k
  d_i \delta - \alpha_i} = T_{\tom_i}^{k} T_i^{-1} E_i\ (k > 0).$$
These are the Drinfeld generators for $\U$.
\end{Remark}

Having defined real root vectors, we define the imaginary root vectors.
For $k>0$, $i \in I_0$, set
\begin{equation*}
\tilde\psi_{i,k d_i} = 
E_{k d_i \delta - \alpha_i}E_{\alpha_i} -q_i^{-2}E_{\alpha_i}E_{k d_i
  \delta-\alpha_i},
\end{equation*} 
and define elements $E_{i, k d_i \delta}\in\Up$ by the functional equation
\begin{equation*} 
(q_i-q_i^{-1})\sum_{k=1}^\infty \ed[i]{k d_i} u^k = \log\left(1
+ \sum_{k=1}^\infty (q_i-q_i^{-1}) {\tilde{\psi}}_{i, k d_i} u^k\right).
\end{equation*}
We have
\begin{equation*}
   \left[ E_{i, k d_i \delta}, E_{j, l d_j\delta} \right] = 0.
\end{equation*}
(For the untwisted case see \cite{Beck}. For general case see
\cite[Theorem 5.3.2]{Da2}.)

\begin{comment}
{\bf Notations}: $\tilde\psi_{i,k d_i} = -\tilde E_{(kd_i\delta,i)}$
(defined in \cite[p.101]{Da2}), $E_{i, kd_i\delta} =
E_{(kd_i\delta,i)}$ (defined in \cite[p.103, p.110, middle of the
proof of Prop.~4.4.10]{Da2}). The normalization is different: $q$ in
\cite{Da2} is our $q_n$ for $\atnt$. For $\att$: $\tilde\psi_{1, kd_1}
= \tilde\psi_k$ in \cite[Def.~3.3, Prop.~3.26(2)]{A},
$E_{1,kd_i\delta} = E_{k\delta}$ in \cite[Prop.~B.3]{A}.
\end{comment}

For each $i \in I_0$ we introduce the ``integral'' imaginary root
vectors (\cite{CP1}, \cite{A}) by
\begin{align}  \label{tobespecialized}
\sum_{k \ge 0} \tilde P_{i, k d_i}u^{k}=
\begin{cases}
\exp\left(\sum_{k \geq
1} \frac{\ed[n]{k}}{[2k]_{n}} u^{k}\right) 
 & \text{if } (\xn,i) = (\atnt,n),\\ 
\exp\left(\sum_{k \geq
1} \frac{\ed[i]{k d_i}}{[k]_{i}} u^{k}\right) 
& \text{otherwise.}
\end{cases}
\end{align} 

They satisfy the following recursive identity:
\begin{align} \label{integralimag} 
\tilde P_{i, k d_i} = 
\begin{cases}
 \frac{1}{[2k]_n}
  \sum_{s=1}^k q_n^{2(s-k)} \tilde{\psi}_{n,s}\tilde P_{n,k-s}
 & \text{if } (\xn,i) = (\atnt,n),\\ \\
 \frac{1}{[k]_i}
  \sum_{s=1}^k q_i^{s-k} \tilde{\psi}_{i,s d_i}\tilde P_{i, (k-s) d_i}
 & \text{otherwise.} 
\end{cases} 
\end{align}
These definitions are based on the definition of imaginary root
vectors which appeared for type $A_1^{(1)}$ in \cite{CP1}. The $\atnt$
case is generalized from \cite{A} where it appears in the $\att$ case.

Let $\Up(>)$ (resp.\ $\Up(<)$, $\Up(0)$) be the
$\Q(q_s)$-subalgebra of $\Up$ generated by the $E_{\beta_k}$ for
$k\le 0$ (resp.\ $E_{\beta_k}$ for $k>0$, $E_{i,kd_i\delta}$ for $k>0$).

The $\tilde P_{i,k d_i}$ ($i\in I_0$, $k>0$) are used to construct
a basis of the imaginary part of $\Up$ as follows. Let $\boc_0 =
(\rho^{(1)}, \rho^{(2)}, \dots, \rho^{(n)})$ be an $n$--tuple of
partitions where each $\rho^{(i)} = (\rho^{(i)}_1 \ge \rho^{(i)}_2 \ge
\dots )$.  For a partition $\rho$, we denote its transpose by
$\rho^t$.  For each $\rho^{(i)}$, define the Schur
function in the $\tilde P_{i,k d_i}$ by
\begin{equation*}
S_{\rho^{(i)}} = \det(\tilde{P}_{i, ({\rho^{(i)}}^t_k-k+m) d_i})_{1\le
k,m\le t},
\end{equation*}
where $t\ge l({\rho^{(i)}}^t)$. This puts the $\tilde P_{i,k d_i}$ in
the role of elementary symmetric functions.  Note that in \cite{BCP,
A} the transpose of $S_{\rho^{(i)}}$ is considered, but this make no
difference.  Denote the product over $i\in I_0$ of $S_{\rho^{(i)}}$ by
\begin{equation*} 
S_{\boc_0} = \prod_{i=1}^n S_{\rho^{(i)}}. 
\end{equation*}

\begin{Definition} \label{precp} 
Let $\boc_+ \in \N^{\Z_{\le 0}}$ and $\boc_- \in \N^{\Z_{>0}}$ be
functions which are almost everywhere $0$.  Let $\boc_0$ as above. Let
$\aC$ denote the set of all such triples $\boc = (\boc_+, \boc_0,
\boc_-).$

For each $p$ and $\boc \in \aC$ we define:
\begin{equation*}
\begin{split}
  &\boc_{+_p} = \left( \boc(p), \boc({p-1}), \boc({p-2}),
   \dots\right)
\\
  &\text{and} \quad
   \boc_{-_p} = \left( \boc({p+1}), \boc({p+2}), \boc({p+3}),
   \dots\right).
\end{split}
\end{equation*}
where the components $\boc(j)$ are actually $\boc_+(j)$
(resp.\ $\boc_-(j)$) when $j \le 0$ (resp.\ $j > 0$).

Define 
\begin{equation*}
\begin{split}
  & E_{\boc_{+_p}} = E_{i_p}^{\boc(p)} T_{i_p}^{-1}
(E_{i_{p-1}}^{\boc(p-1)}) T_{i_p}^{-1} T_{i_{p-1}}^{-1}
(E_{i_{p-2}}^{\boc(p-2)}) \cdots
\\ 
  & E_{\boc_{-_p}} = \cdots
T_{i_{p+1}} T_{i_{p+2}} (E_{i_{p+3}}^{(\boc(p+3))})
T_{i_{p+1}}(E_{i_{p+2}}^{(\boc(p+2))}) E_{i_{p+1}}^{(\boc(p+1))},
\end{split}
\end{equation*}
where the exponents above are written $\boc(j)$ when they should be
$\boc_+(j)$ or $\boc_-(j)$ for $j \le 0$ or $j > 0$ respectively.

Note that when $p=0$, the monomials $E_{\boc_{\pm_0}}$ are formed by
multiplying the $(E_{\beta_k})^{(\boc_\pm(k))}$ for real roots
$\beta_k\in\aR_\gtrless$ in the order \eqref{rootorder}.  Also, in
this case we will omit the subindex and write $\boc_\pm$ instead of
$\boc_{\pm_0}$---there should be no confusion with the $\boc_\pm$
above.

For $\boc \in \aC, p \in \Z$ we define 
(cf.\ \cite[40.2.3]{Lu-Book}):
\begin{equation*}
\begin{split}
L(\boc,p) = & \left(E_{i_{p}}^{(\boc(p))}
    T_{i_{p}}^{-1}(E_{i_{p-1}}^{(\boc(p-1))})
    T_{i_{p}}^{-1}T_{i_{p-1}}^{-1}(E_{i_{p-2}}^{(\boc(p-2))})
    \cdots\right) \\ & \times T_{i_{p+1}}T_{i_{p+2}} \cdots T_{i_0}
    (S_{\boc_0}) \\ & \times \left(\cdots T_{i_{p+1}} T_{i_{p+2}}
    (E_{i_{p+3}}^{(\boc(p+3))}) T_{i_{p+1}}(E_{i_{p+2}}^{(\boc(p+2))})
    E_{i_{p+1}}^{(\boc(p+1))}\right) \\ = & \; E_{\boc_{+_p}} \big(
    T_{i_{p+1}}T_{i_{p+2}} \cdots T_{i_0} (S_{\boc_0})\big)
    E_{\boc_{-_p}},
\end{split}
\end{equation*}
for $p\le 0$. For $p\ge 1$, we replace the middle part by
$T_{i_p}^{-1} T_{i_{p-1}}^{-1} \dots T_{i_2}^{-1}
    T_{i_1}^{-1}(S_{\boc_0}).$

Note that 
\begin{equation} \label{eq:trans} 
\begin{split}
 & L(\boc,p-1) = T_{i_{p}}L(\boc,p) \text{ if
$\boc(p) = 0$},  \\ & L(\boc,p+1) =
T_{i_{p+1}}^{-1}L(\boc,p) \text{ if $\boc({p+1}) = 0$},
\end{split}
\end{equation}

When $p=0$ we will denote $L(\boc,0)$ by $B_\boc$.
In this case, each $\boc \in \aC$ indexes an element which we will
call of ``PBW--type'' in $\Up$:
\begin{equation}   \label{def:PBW}    
 B_{\boc} = E_{\boc_+} \cdot S_{\boc_0} \cdot E_{\boc_-}.
\end{equation} 
We will call $B_\boc$ for $\boc \in \aC$ {\it purely imaginary\/} if
$\boc_+ = \boc_- = 0$.  In this case we may write $\boc_0$ instead of
$\boc = (0,\boc_0,0)$.

For each $p \in \Z$, we define a partial ordering $\prec_p$ on $\{
\boc \ | \ \boc \in \aC \}$ by letting $\boc \prec_p \boc'$ if and
only if
\begin{equation*}
  \boc_{+_p}   \le \boc_{+_p}' \quad \text{and}\quad
   \boc_{-_p} \le \boc_{-_p}' \quad \text{and one of these is strict.}
\end{equation*}
Here both $\le$ are the lexicographic ordering from left to right.
For example, the first inequality means either $\boc = \boc'$ or
$\boc(p) = \boc'(p), \,
\boc({p-1}) = \boc'({p-1}), \,  \dots, \,$ $
\boc({p-k+1}) = \boc'({p-k+1})$ and 
$\boc({p-k}) < \boc'({p-k})$ for some $k \ge 0$. 
Note that
\begin{equation} \label{eq:ord_trans}
\text{if $\boc(p) = \boc'(p) = 0$, then
\(
   \boc \prec_p \boc' \Leftrightarrow \boc \prec_{p-1} \boc'.
\)
}
\end{equation}
\end{Definition}

We now state the main the theorem of this section.

\begin{Theorem} \label{thm:intbase}
\textup{(i)} For each $p\in \Z$, $\{ L(\boc,p) \mid \boc \in \aC \}$  
is an almost orthonormal basis of the $\Q(q_s)$-vector space $\Up$,
i.e., $(L(\boc,p), L(\boc',p)) \in \delta_{\boc,\boc'} +
q_s^{-1}\Z[[q_s^{-1}]]\cap \Q(q_s)$.

\textup{(ii)} Let $p \in \Z$. The transition matrix between $\{
L(\boc,p) \mid \boc \in \aC \}$ and the global crystal basis of $\U^+$
is upper-triangular with $1$'s on the diagonal and with above diagonal
entries in $q_s^{-1}\Z[q_s^{-1}]$.
\end{Theorem}

The property (i) follows from (ii) and the almost orthonormality of
the global crystal base \cite[14.2.3]{Lu-Book}. However, we will use
(i) during the proof of (ii) and we prove it independently. Let us
also remark that (ii) implies that $\{ L(\boc,p) \mid \boc \in \aC \}$
is a basis of the integral form $\Ua^+$.  The proof of (i) will be
given in \subsecref{subsec:alonb}. We postpone that of (ii) until
\secref{sec:sign}.

\subsection{Vertex subalgebras}

A key part of describing the PBW basis is its reduction to
``vertex subalgebras''.  The following proposition describes
the $n$ ``vertex'' subalgebras of $\U = \U(X_N^{(k)})$:

\begin{prop} \label{ivertex}
Let $i\in I_0.$ Let $\Uiv$ be the $\Q(q_i)$-subalgebra of $\U =
\U(\xn)$ generated by
\begin{equation*}
\{E_i, E_{d_i \delta -\alpha_i}, F_i, F_{d_i \delta -\alpha_i},
K_i^{\pm 1}, K_{d_i \delta}^{\pm 1} \}.
\end{equation*}  Then $\Uiv$ is isomorphic as a $\Q(q)$-algebra to
$\U(A_{1}^{(1)})$ in all cases except when $\xn = \atnt$ and $i = n$,
i.e., $(\alpha_i, \alpha_i) = 1$.  The isomorphism $\U(A_{1}^{(1)})
\rightarrow \Uiv$ is given by $E_1 \mapsto E_i, E_0 \mapsto E_{d_i
\delta - \alpha_i}$, $F_1 \mapsto F_i, F_0 \mapsto F_{d_i \delta -
\alpha_i}$, $q \mapsto q_i$. For $(\xn, i) = (\atnt, n)$, $\Uiv$ is
isomorphic as a $\Q(q^{1/2})$-algebra to $\att$, and the isomorphism
$\att\to \Uiv$ is given by $E_0 \mapsto E_{\delta - 2 \alpha_n}$, $F_0
\mapsto F_{\delta - 2\alpha_n}$, $q^{1/2} \rightarrow q^{1/2}$.
\textup(In particular, $E_{\delta - 2 \alpha_n}$, $F_{\delta -
2\alpha_n}\in \Uiv$.\textup)
\end{prop}

\begin{proof}  For the untwisted  case see \cite{Beck}.  
In the twisted case, the result is due to \cite{Da2}.  Note that in
\cite{Da2} the quantum algebra of type $A_{2n}^{(2)}$ is normalized
differently (the invariant bilinear form on $\h^*$ is $2$ times the
one here).
\end{proof}

Next we have:
\begin{prop}  \label{tildePinA}
For $i \in I_0, k > 0$ we have $\tilde P_{i, kd_i} \in \Ua^+.$
\end{prop}
\begin{proof}
This is proved as in \cite[Corollary 2.2]{BCP} for the symmetric case.  In
the $(\xn,i) = (A_{2n}^{(2)},n)$ case this follows from Proposition
\ref{ivertex} and \cite[Corollary 8.6]{A}.
\end{proof}

\begin{prop}\label{prop:int}  Let $\boc \in \aC$. We have
$L(\boc,p) \in\Ua^+$.
\end{prop}
\begin{pf}
By \cite[6.3.2]{Da2} (see \cite{Beck2} for the untwisted case) we have
$T_{i_{p+1}}T_{i_{p+2}} \cdots T_{i_0} (S_{\boc_0})$, $T_{i_p}^{-1}
T_{i_{p-1}}^{-1} \dots T_{i_2}^{-1} T_{i_1}^{-1}(S_{\boc_0})\in \Up$.
Since $T_w$ preserves $\Ua$, they are contained in $\Ua\cap\Up =
\Ua^+$ by \propref{tildePinA}. We also have $E_{\boc_{\pm p}}\in\Ua^+$ 
by \cite[41.1.3]{Lu-Book}.
\end{pf}

Next we cite a key property of the $\tilde P_{i, k d_i}$: 
\begin{prop}  \label{keyidentity} Let $k > 0$. 
\begin{align}
\tilde P_{i, k d_i} = 
\begin{cases}
 E_{d_i \delta - \alpha_i}^{(k)} E_{i}^{(k)} + q_s^{-1} x
 & \text{if } (\xn,i) \neq (A_{2n}^{(2)},n), \\
 E_{\delta -  2 \alpha_n}^{(k)} E_{n}^{(2k)} + q_s^{-1} x
 &  \text{if } (\xn,i) = (A_{2n}^{(2)},n), 
\end{cases} 
\end{align}
where $x$ is a sum of terms $B_\boc$ with coefficients in
$\Z[q_s^{-1}]$ where for each term $\boc_+, \boc_- \neq 0$.
Furthermore, for each such $\boc \in \aC$, only imaginary root vectors
$P_{i, t d_i}$ with $0 < t < k$ appear in $S_{\boc_0}$.
\end{prop}  
\begin{proof}
This is derived as in \cite[Proposition 2.2 and eq.~(4.9)]{BCP}, where
it appears for the untwisted affine case. For the $\atnt$ this
appears as a special case of \cite[Theorem 8.5]{A}.
\end{proof}

\subsection{Almost orthonormality}\label{subsec:alonb}

The following proposition is central to our calculations.
\begin{prop}{\label{specialinner}} 
For $i,j\in I_0$, $k, k' >0$ we have 
\begin{equation*}
(\tilde{P}_{i,kd_i},\tilde{P}_{j,k'd_j}) \equiv \delta_{k,k'}\delta_{ij}
\pmod{q_s^{-1}\Ainfty}.\end{equation*}
\end{prop}
\begin{proof} 
This result appears in the symmetric case in \cite{BCP} and in the
$\att$ case it appears in \cite{A}.  In general, the proof is
analogous to one of the previous cases.  In the non--symmetric case
where $i \neq j$, we may assume $a_{ij} = -1$ since the condition is
symmetric in $i$ and $j$.  The result then follows from the following
identity:
\begin{align}
r_j(E_{i,kd_i\delta}) =
\begin{cases} 
(-1)^{k+1} q_j^{-1} (1 - q_j^{-2}) E_{k d_i \delta -
\alpha_j} & \text{ if }d_j | k d_i  \\ 0 & \text{ otherwise.}
\end{cases}
\end{align}
which is derived from the following special case 
 of \cite[Theorem 5.3.2]{Da2} (just as in the symmetric case
the previous identity is derived in \cite{BCP}):
\begin{align}
[E_{i, k d_i \delta}, F_j] =
\begin{cases} 
(-1)^{k+1} \frac{[r]_{q_i}}{r} K_{\alpha_j} E_{k d_i \delta -
\alpha_j} & \text{ if }d_j | k d_i  \\ 0 & \text{ otherwise.}
\end{cases}
\end{align}
\end{proof} 

Next we need the following result regarding the coproduct formula for
the imaginary root vectors.  For any algebra $A$, let $A_+$ denote
its augmentation ideal.

\begin{prop} \label{prop:coprod} Let $k > 0, i \in I_0$.  Then 
\begin{align*}
r(\tilde{P}_{i,k d_i}) =\sum_{s=0}^{k}& \tilde{P}_{i,s d_i}\ot
\tilde{P}_{i,(k-s) d_i} \\ & + {\text{terms in }}
\Up(<)_+\Up(0)\otimes \Up(0)\Up(>)_+.
\end{align*}
\end{prop}
\begin{proof}
The proof of this follows from the relation between the braid group
action and the coproduct (see \cite[3.1.5, 37.3.2]{Lu-Book}) by using
\remref{rem:loopgen}.  In the untwisted case, the argument is given in
\cite{Da}.  In general, the argument is identical.
\end{proof}

We have the following result which is proved as in \cite{BCP} (see
\cite{A} for the $\att$ case):
\begin{prop} \label{schurorth} Let $\boc_0, \boc'_0$ be two $n$--tuples of 
partitions as above. Then
\begin{equation} (S_{\boc_0}, S_{\boc'_0}) \equiv
  \delta_{\bf{c_0},\bf{c'_0}} \pmod{q_s^{-1}\Ainfty}.
\end{equation} 
\end{prop}

\begin{proof}[Proof of \thmref{thm:intbase}(i)]
Let $\boc,\boc' \in \aC.$ Suppose $p \ge 1$.
We have
\begin{align*}
   T_{i_p}^{-1} T_{i_{p-1}}^{-1} \cdots T_{i_2}^{-1}
    T_{i_1}^{-1}(S_{\boc_0}) \in \Up\cap T_{i_p}^{-1}(\Up). 
\end{align*}
as we already remarked in the proof of \propref{prop:int}.
By \cite[38.2.1]{Lu-Book} we have
\begin{align*}
   \begin{gathered}
       \left(T_{i_p}^{-1} T_{i_{p-1}}^{-1} \cdots T_{i_2}^{-1}
    T_{i_1}^{-1}(S_{\boc_0}),
    T_{i_p}^{-1} T_{i_{p-1}}^{-1} \cdots T_{i_2}^{-1}
    T_{i_1}^{-1}(S_{\boc'_0})\right)
\\
 = \left(T_{i_{p-1}}^{-1} \cdots T_{i_2}^{-1} T_{i_1}^{-1}(S_{\boc_0}),
    T_{i_{p-1}}^{-1} \cdots T_{i_2}^{-1} T_{i_1}^{-1}(S_{\boc'_0})\right).
   \end{gathered}
\end{align*}
By the induction and Proposition~\ref{schurorth} these are equal to
$\delta_{\boc, \boc'}$ modulo $q_s^{-1}\Ainfty$.

We have by \cite[40.2.4]{Lu-Book}
\begin{equation} \label{eq:inner}
  (L(\boc,p), L(\boc',p)) =
  (S_{\boc_0},S_{\boc'_0}) \prod_{s\in\bz}
  (E_{\alpha_{i_s}}^{(\boc(s))},E_{\alpha_{i_s}}^{(\boc'(s))})
  \equiv \delta_{\boc, \boc'} \mod{q_s^{-1}\Ainfty}.
\end{equation}
Moreover, we know $(L(\boc,p), L(\boc',p))\in \Z[[q_s^{-1}]]$ since
$L(\boc,p)$, $L(\boc',p)\in\Ua^+$ by \cite[14.2.6]{Lu-Book}. Hence we
get the almost orthonormality.

For $p<0$, we have the same result thanks to
\begin{equation*}
T_{i_{p+1}}T_{i_{p+2}} \cdots T_{i_0}(S_{\boc_0})
\in \Up\cap T_{i_p}^{-1}(\Up).
\end{equation*}

By \eqref{eq:inner} $\{ L(\boc,p) \}$ is linearly independent. The PBW
theorem says that the dimension of a weight space of $\Up$ is equal to
the number of $\boc$'s with the given weight. Thus it is a basis.
\end{proof}

\begin{Corollary} \label{cor:qbasis}
Let $\La(-\infty) = \{ x\in \Up \mid (x,x)\in \A_0
\}$. Then it is an $\A_0$-submodule of $\Up$, and $\{ \ov{L(\boc,p)}
\mid \boc\in \aC \}$ is its $\A_0$-basis for each $p\in\Z$.
\end{Corollary}

\begin{proof}
The assertion follows from \cite[14.2.2]{Lu-Book}.
\end{proof}

Note that our proof of this statement, as well as that of
\thmref{thm:intbase}(i), is independent of the existence of the
global crystal basis. The above definition of $\La(-\infty)$ coincides
with one in \cite[17.3.3]{Lu-Book}, while it was a characterizing
property in \cite[5.1.4]{Kas91}.

\begin{prop} \label{crystalbase}
For every $\boc \in \aC, p \in \Z$, there exist $b = b(\boc,p)\in
\B(-\infty)$ and $\sgn(\boc,p)=\pm$ such that
\begin{equation}
  L(\boc,p) \equiv \sgn(\boc,p) G(b(\boc,p))
  \mod q_s^{-1}\ov{\La(-\infty)}.
\end{equation}
\end{prop}
\begin{pf}
By \cite[14.2.2]{Lu-Book}, we know that if $x\in \Ua^+$ satisfies
$(x,x)\in 1 + q_s^{-1}\Ainfty$, then $x\in\ov{\La(-\infty)}$ and there
exists $b\in\B(-\infty)$ such that $x \equiv \pm G(b)\mod q^{-1}_s
\ov{\La(-\infty)}$.
The assertion follows from \thmref{thm:intbase}(i).
\end{pf}

In \secref{sec:sign} we will show $\sgn(\boc,p) = 1$.

Note that the map $\aC\to \B(-\infty)$ given by $\boc\mapsto
b(\boc,p)$ is bijective for any $p$, since both are bases.

\begin{Remark}
For any $p$, $q\in\Z$ there exists a bijection
$\boc\in\aC\leftrightarrow\boc'\in\aC$ given by
\(
   b(\boc,p) = b(\boc',q)
\)
by \propref{crystalbase}.  It is extremely interesting problem to give
an explicit description of the bijection. For finite type $\g$, the
same construction gives the piecewise linear bijections
\cite[42.1.3]{Lu-Book}, which have been studied by various peoples.
Let $b = b(\boc,p-1)$ with $\boc(p) = 0$. By \cite{Saito} we have
\(
   b(\boc,p) = \te_i^{*\varepsilon_i(b)} \tf_i^{\varphi_i(b)} b
\),
where $i = i_p$. This is a first step towards this problem.
\end{Remark}

\subsection{Upper triangular property of the bar involution}

This subsection is a small detour. We give a proof that the bar
involution is upper triangular with respect to our basis $\{ L(\boc,p)
\}$. For symmetric or type $\att$, this together with \cite{A,BCP}
gives us the elementary algebraic definition of the global crystal
basis as explained in the Introduction. The reader in a hurry may skip
the rest of this section.

We will need a ``reordering lemma'' to prove the upper-triangularity
property of the bar action with respect to the ordering $\prec_p$.

\begin{Lemma} \label{lem:reordering} Let $p \in \Z$. 
Let $\boc, \boc' \in \aC$. Write 
\begin{equation}
L(\boc,p) L(\boc',p) = \sum_{\boc''} a_{\boc,\boc'}^{\boc''}
L(\boc'',p),
\end{equation}
where $a_{\boc,\boc'}^{\boc''} \in \Q(q_s)$.  

\textup{(i)} Then for each $\boc''$ in the above sum, $\boc_{+_p}''
\ge \boc_{+_p}$ and $\boc_{-_p}'' \ge \boc'_{-_p}$ with respect to the
lexicographical ordering.  

\textup{(ii)} Further, if $L(\boc,p) = E_{\boc_{+_p}}$
(resp.\ $E_{\boc_{-_p}}$) and $L(\boc',p) = E_{\boc_{+_p}'}$
(resp.\ $E_{\boc_{-_p}'}$), then $\boc_{+_p}'' > \boc_{+_p}$
(resp.\ $\boc_{-_p}'' > \boc_{-_p}'$) for each summand.
\end{Lemma}
\begin{proof}
We prove this for $p = 0$, noting that for $p \neq 0$ the proof is
identical. Assume $L(\boc,0) = E_{\boc_+} S_{\boc_0} E_{\boc_-}$ and
$L(\boc',0) = E_{\boc_+'} S_{\boc_0'} E_{\boc'_-}$.  Write the
expression 
\begin{equation}
S_{\boc_0} E_{\boc_-} E_{\boc_+'} S_{\boc_0'} E_{\boc'_-}
=\sum_{\bob \in \aC} a_{\boc,\boc'}^{\bob} E_{\bob_+} S_{\bob_0}
E_{\bob_-}.
\end{equation}
Here and for the remainder of the section any structure constants
(such as $a_{\boc,\boc'}^{\bob}$) are assumed to be in $\Q(q_s)$
unless otherwise stated. We have
\begin{equation} \label{eq:reord0}
L(\boc,0) L(\boc',0) = \sum_{\bob} a_{\boc,\boc'}^{\bob}
E_{\boc_+} E_{\bob_+} S_{\bob_0} E_{\bob_-}.
\end{equation}

For a given summand, if $\boc_+ = 0$ then clearly $\boc_+'' = \bob_+
\ge \boc_+$ in the Lemma. Assume $\boc_+ > 0$.  For any given $\bob$
in the sum, we may assume $\bob_+ > 0$ or else that summand also
fulfills the requirement of the Lemma. Under the assumption that
$\bob_+ > 0$ and $\boc_+ > 0$, the $\boc_+'' \ge \boc_{+}$ part in (i)
and the $\boc_+'' > \boc_{+}$ part in (ii) follow if we check that
for a fixed $\bob_+$ in \eqref{eq:reord0},
\begin{equation} \label{eq:reord1} 
E_{\boc_+} E_{\bob_+}  = \sum_{\substack{\bod_+ \in \N^{\car_>} \\ 
\bod_+ > \boc_+}} a_{\boc_+,\bob_+}^{\bod_+} E_{\bod_+}.
\end{equation}
Let $k > k' > 0$. We have the following useful identity \cite{LS} (it
was proved there in the finite type case, but the same proof works
here):
\begin{equation} \label{eq:LSid}
E_{\beta_{-k}} E_{\beta_{-k'}} = q^{(\beta_k,\beta_{k'})}
E_{\beta_{-k'}} E_{\beta_{-k}} + \sum_{\bod_+'} a_{\bod_+'}
E_{\bod_+'},
\end{equation}
where $\bod_+'(j) = 0$ if $j \ge k$ or $j \le k'$.  Note that this
condition is equivalent to saying that $e_{-k'} > \bod_+' > e_{-k}$,
where $e_{k}$ is the tuple whose $j$-th position is $\delta_{jk}$.  

Consider the set $S$ of all monomials of weight $\gamma = \wt
(E_{\boc_+} E_{\bob_+}) \in Q_+$ formed from the real root vectors
$E_{\beta_{-k}}, k \ge 0$.  $S$ is a finite set.  We will order this
set by a lexicographic order on the monomials, where $E_{\beta_{-k'}}
> E_{\beta_{-k}}$ if $k > k' \ge 0$.  In this ordering, a monomial $M$ is
in the PBW order if and only if it is maximal among all monomials
where $E_{\beta_{-k}}$ appears for all $k$ the same number of times as
in $M$.

On the left hand side of \eqref{eq:reord1}, moving from left to right
take the first root vector in $E_{\boc_+} E_{\bob_+}$ which is out of
PBW order, i.e., the first root vector which is larger than in it's
immediate predecessor.  Use \eqref{eq:LSid} to reorder these two root
vectors.  By \eqref{eq:LSid}, we obtain a linear combination of
monomials from $S$, each greater in lexicographic order than
$E_{\boc_+} E_{\bob_+}$.  Repeating this for each summand, and taking
into account that $S$ is a finite set, we ultimately obtain a linear
combination of elements of $S$, each of which is maximal among
monomials formed from the same root vectors.  Thus, each is in PBW
form and also larger than $E_{\boc_+}$ in lexicographic order, so that
we obtain $E_{\boc_+} E_{\bob_+}$ in the form of \eqref{eq:reord1} as
required.
\end{proof}

\begin{Proposition} \label{bartriangular}  
Let $\boc \in \aC, p \in \Z$.  Then 
\begin{equation} \label{eq:barupt}
\overline{L(\boc,p)} = L(\boc,p) + \sum_{\boc \prec_p \boc'} a_{\boc,
\boc'} L(\boc',p),
\end{equation}
where $a_{\boc, \boc'} \in \Q(q_s)$. 
\end{Proposition}

\begin{proof}
We first prove
\begin{description}
\item[(a)] Fix $\boc_0$ and $q$. Suppose that \eqref{eq:barupt} is
true for
\(
   T_{i_{q+1}}T_{i_{q+2}} \cdots T_{i_0} (S_{\boc_0})
\)
($q \le 0$), or
\(
  T_{i_q}^{-1} T_{i_{q-1}}^{-1} \cdots T_{i_1}^{-1}(S_{\boc_0})
\)
($q \ge 1$). Let $p \ge q$. If $\boc$ satisfies $\boc_{-_p} = 0$,
$\boc_{+_q} = 0$ and $\boc_0$ is the given one, \eqref{eq:barupt} is
also true for $L(\boc,p)$. Furthermore, the condition
$\boc\prec_p\boc'$ can be replaced by the stronger condition
$\boc_{+_p} < \boc'_{+_p}$. (The other condition $(0 =) \boc_{-_p} \le
\boc'_{-_p}$ is trivially satisfied.)
\end{description}

We are considering
\begin{multline*}
    L(\boc,p) = E_{i_{p}}^{(\boc(p))}\times
    T_{i_{p}}^{-1}(E_{i_{p-1}}^{(\boc({p-1}))})\times
    T_{i_{p}}^{-1}T_{i_{p-1}}^{-1}(E_{i_{p-2}}^{(\boc({p-2}))})\times
    \cdots 
\\
    \cdots 
    \times T_{i_p}^{-1}\cdots T_{i_{q+2}}^{-1}(E_{i_{q+1}}^{(\boc(q+1))})
    \times T_{i_{p+1}}T_{i_{p+2}} \cdots T_{i_0} (S_{\boc_0}).
\end{multline*}
(When $p \ge 1$, the last part is 
$T_{i_p}^{-1} T_{i_{p-1}}^{-1} \dots T_{i_2}^{-1}
    T_{i_1}^{-1}(S_{\boc_0}).$)
We prove the assertion by the induction on $p$. When $p=q$,
\eqref{eq:barupt} is true by the assumption.

First assume $\boc(p) = 0$. We consider
\(
   L(\boc,p-1) = T_{i_p}L(\boc,p)
\)
(see \eqref{eq:trans}). By the induction hypothesis, we have
\begin{equation*}
   \ov{T_{i_p}L(\boc,p)} = \ov{L(\boc,p-1)} =
   L(\boc,p-1) + 
   \sum_{\boc_{+_{p-1}} < \boc'_{+_{p-1}}}
   a_{\boc,\boc'}^{p-1} L(\boc',p-1).
\end{equation*}
(We put the superscript $p-1$ in order to clarify its dependence on
$p-1$ in this part. For the other part, $p$ will be fixed, and there
will be no confusion.)
We apply the composition $T_{i_p}^{-1} \circ \pi^{i_p}$ to both
sides. By \eqref{eq:prop1.9}, the left hand side becomes
${}^{i_p}\pi(\ov{L(\boc,p)})$, which is equal to $\ov{L(\boc,p)}$ modulo
$E_{i_p}\Up$. For the right hand side, we use \eqref{eq:trans}. We get
\begin{equation*}
   \ov{L(\boc,p)} \in L(\boc,p) + 
   \sum_{\substack{\boc_{+_{p-1}} < \boc'_{+_{p-1}} \\ \boc'(p) = 0}}
   a_{\boc,\boc'}^{p-1} L(\boc',p-1) + E_{i_p}\Up.
\end{equation*}
The condition $\boc'(p) = 0$ comes from $\pi^{i_p} (L(\boc',p-1))
\neq 0$. (Otherwise, $L(\boc',p-1)\in \Up E_{i_p}$ and 
$\pi^{i_p} (L(\boc',p-1)) = 0$.)
The part in $E_{i_p} \U^+$ is a linear combination of $L(\boc'',p)$'s
with $\boc''(p) > 0$. They satisfy $\boc_{+,p} < \boc''_{+,p}$ since
$\boc(p) = 0$. The summation in the second term can be replaced as
$\sum_{\boc_{+_p} < \boc'_{+_p}, \boc'(p) = 0}$ by \eqref{eq:ord_trans}.
Thus we have the assertion under the assumption $\boc(p) = 0$.

Next we assume $\boc(p) > 0$. Let us define $\tilde\boc$
by setting $\tilde\boc(p) = 0$ and all other entries are the
same as $\boc$. We have
\begin{equation*}
    L(\boc,p) = E_{i_p}^{(\boc(p))}L(\tilde\boc,p).
\end{equation*}
Since $\tilde\boc(p) = 0$, we have just proved
\begin{equation*}
    \ov{L(\tilde\boc,p)} = L(\tilde\boc,p) + 
    \sum_{\tilde\boc_{+_p} < \tilde\boc'_{+_p}} 
    a_{\tilde\boc,\tilde\boc'}^p L(\tilde\boc',p).
\end{equation*}
Therefore
\begin{equation*}
    \ov{L(\boc,p)} = E_{i_p}^{(\boc(p))}\ov{L(\tilde\boc,p)}
    = L(\boc,p) + 
    \sum_{\tilde\boc_{+_p} < \tilde\boc'_{+_p}}
    a_{\tilde\boc,\tilde\boc'}^p 
    \begin{bmatrix}
      \boc(p) + \tilde\boc'(p) \\ \boc(p)
    \end{bmatrix}_{q_i}
    L(\boc',p),
\end{equation*}
where $\boc'$ is defined by setting $\boc'(p) =
\tilde\boc'(p)+\boc(p)$ and other entries are the same as
$\tilde\boc'$. We have $\tilde\boc_{+_p} < \tilde\boc'_{+_p}
\Longleftrightarrow\boc'_{+_p} < \boc_{+_p}$, and therefore the
assertion.

Similarly we have
\begin{description}
\item[(a')] Fix $\boc_0$ and $q$. Suppose that \eqref{eq:barupt} is
true for
\(
   T_{i_{q+1}}T_{i_{q+2}} \cdots T_{i_0} (S_{\boc_0})
\)
($q \le 0$), or
\(
  T_{i_q}^{-1} T_{i_{q-1}}^{-1} \cdots T_{i_1}^{-1}(S_{\boc_0})
\)
($q \ge 1$). Let $p \le q$. If $\boc$ satisfies $\boc_{+_p} = 0$,
$\boc_{-_q} = 0$ and $\boc_0$ is the given one, \eqref{eq:barupt} is
also true for $L(\boc,p)$. Furthermore, the condition
$\boc\prec_p\boc'$ can be replaced by the stronger condition
$\boc_{-_p} < \boc'_{-_p}$. (The other condition $(0 =) \boc_{+_p} \le
\boc'_{+_p}$ is trivially satisfied.)
\end{description}

Our next task is to show
\begin{description}
\item[(b)] Fix $\boc_0$. Suppose that \eqref{eq:barupt} is true for
$S_{\boc_0}$. Then it is also true for $L(\boc,p)$ with $\boc_0$ is
the given one.
\end{description}

By {\bf (a)}, {\bf (a')}, if \eqref{eq:barupt} is true for
$S_{\boc_0}$, then it is also true for \( T_{i_{p+1}}T_{i_{p+2}}
\cdots T_{i_0} (S_{\boc_0}) \) ($p \le 0$), or \( T_{i_p}^{-1}
T_{i_{p-1}}^{-1} \cdots T_{i_1}^{-1}(S_{\boc_0}) \) ($p \ge 1$). (The
conditions on $\boc_{\pm,p}$, $\boc_{\pm,q}$ are vacuous since
$\boc_{\pm} = 0$.) This is a special case of {\bf (b)}.
We return back to general $L(\boc,p)$ as in {\bf (b)}. We decompose it
as $L(\boc,p) = L(\boc_{+_p},p) L(\boc_0,p) L(\boc_{-_p}, p)$ where
$\boc_{\pm p}$, $\boc_0$ are understood as elements of $\aC$ by
setting other entries as $0$. By {\bf (a)} and the above special case
of {\bf (b)}, we have
\begin{equation*}
   \ov{L(\boc_{+_p},p) L(\boc_0,p)}
   = L(\boc_{+_p},p) L(\boc_0,p) 
   + \sum_{\boc_{+_p} < \boc'_{+,p}} a_{\boc_{+_p},\boc'} L(\boc',p).
\end{equation*}
By {\bf (a')} we have
\begin{equation*}
   \ov{L(\boc_{-_p},p)}
   = L(\boc_{-_p},p) 
   + \sum_{\boc_{-_p} < \boc''_{-_p}} a_{\boc_{-_p},\boc''} L(\boc'',p).
\end{equation*}
Note that the assumption of {\bf (a')} is trivially satisfied since
$\boc_0 = 0$ in this case. Therefore
\begin{multline*}
   \ov{L(\boc,p)}
   = L(\boc,p) + 
   \sum_{\boc_{+_p} < \boc'_{+,p}} a_{\boc_{+_p},\boc'}
   L(\boc',p) L(\boc_{-_p},p)
\\
   + \sum_{\boc_{-_p} < \boc''_{-_p}} a_{\boc_{-_p},\boc''}
   L(\boc_{+_p},p) L(\boc_0,p) L(\boc'',p)
\\
   + \sum_{\boc_{+_p} < \boc'_{+,p}, \boc_{-_p} < \boc''_{-_p}}
   a_{\boc_{+_p},\boc'} a_{\boc_{-_p},\boc''} L(\boc',p) L(\boc'',p).
\end{multline*}
By \lemref{lem:reordering}, $L(\boc',p) L(\boc_{-_p},p) = \sum_{\bod}
a_{\boc',\boc_{-_p}}^{\bod} L(\bod, p)$ where the summation is
over $\bod$ satisfying $\boc'_{+_p} \le \bod_{+_p}$, $\boc_{-_p}\le
\bod_{-_p}$. Since $\boc_{+_p} < \boc'_{+,p}$, we have
$\boc_{+_p} < \bod_{+_p}$, and hence $\boc \prec_p \bod$. (Recall that 
one of the inequalities must be strict in the definition of $\prec_p$.)
The other two summations can be handled in the same way, and we have
the assertion {\bf (b)}.

Next we replace the inequality $\boc\prec_p\boc'$ by a further stronger
inequality in a special case.
\begin{description}
\item[(c)] Consider the case $p=0$. Let $\boc_k \in \aC$ be such that
\(
  L(\boc_k,0) 
\)
is equal to
\(
  E^{(k)}_{d_i \delta - \alpha_i}
\)
if $(\xn,i)\neq (\atnt, n)$, and
\(
  E^{(k)}_{\delta - 2\alpha_n}
\)
if $(\xn,i) = (\atnt, n)$.
Then \eqref{eq:barupt} is true for $L(\boc_k,0)$ with
$\boc_k\prec_p\boc'$ replaced by $(0 =) \boc_{k,+} < \boc'_+$
and $\boc_{k,-} < \boc'_-$.
\end{description}
Here we have written $\boc_{k,-}$ instead of $(\boc_k)_-$,
etc.

We check
\begin{equation} \label{kdelta}
  \ov{L(\boc_k,0)} = 
   L(\boc_k,0) + \sum_{
   \substack{\boc_{k,+} < \boc'_+\\ \boc_{k,-} < \boc'_-}}
   a_{\boc_k,\boc'} L(\boc',0).
\end{equation}
By {\bf (a')}, we already know $\boc_{k,-} < \boc'_-$. Thus we only
need to show $\boc'_+\neq 0$. Suppose that $\boc'_+ = 0$,
and hence $L(\boc',0) = S_{\boc'_0} E_{\boc'_-}$. Then
\(
   \wt(L(\boc',0)) = 
   \wt(S_{\boc'_0}) + \wt(E_{\boc'_-}),
\)
which is equal to $k (d_i\delta - \alpha_i)$ if $(\xn,i) \neq (\atnt,n)$,
and $k (\delta - 2\alpha_n)$ if $(\xn,i) = (\atnt,n)$.
Since $\wt(S_{\boc'_0})\in \N\delta$, $E_{\boc'_-}$ must be a product of
\(
   E_{l d_i\delta-\alpha_i}^{(c_l)}
\)
when $(\xn,i) \neq (\atnt,n)$, and
\(
   E_{(2l-1)\delta-2\alpha_n}^{(c_l)},
\)
\(
   E_{m\delta-\alpha_n}^{(d_m)}
\)
when $(\xn,i) = (\atnt,n)$ (in an appropriate order).
When $(\xn,i) \neq (\atnt,n)$, we have $\sum c_l = k$ and
$\wt(S_{\boc'_0}) + \sum l c_l = k$. These equations simultaneously
hold only if $c_1 = k$, $c_l = 0$ ($l\neq 1$), and $\wt(S_{\boc'_0}) =
0$ (i.e., $\boc'_0 = 0$). When $(\xn,i) = (\atnt,n)$, we have $\sum
2c_l + \sum d_m = 2k$, $\sum (2l-1)c_l + \sum md_m + \wt(S_{\boc'_0})
= k$. Thus $\sum (2l-2)c_l + \sum (m-\frac12) d_m + \wt(S_{\boc'_0}) =
0$, and hence $c_l = 0$ ($l\ge 2$), $d_m = 0$, $\wt(S_{\boc'_0}) = 0$.
We get $c_1 = k$. In both cases, we have $L(\boc',0) =
E_{\boc_{k,-}}$, which is impossible.

The following together with {\bf (b)} completes the proof.
\begin{description}
\item[(d)] \eqref{eq:barupt} is true for $S_{\boc_0}$ for any
  $\boc_0$.  Namely $\ov{S_{\boc_0}} = S_{\boc_0} + \sum_{0 <
    \boc'_\pm} a_{\boc_0,\boc'} L(\boc',0)$
\end{description}
(The inequality $\boc_0 \prec_p \boc'$ is equivalent to the above
inequalities since $\wt(L(\boc',0))\in\Z\delta$.)

We prove {\bf (d)} by the induction on the length of $\boc_0$. When
$\ell(\boc_0) = 0$, we understand $S_{\boc_0} = 1$, and the assertion
is trivial. We assume $(\xn,i)\neq (\atnt,n)$ now, but the argument
works even when $(\xn,i) = (\atnt,n)$ with obvious modifications.
By \propref{keyidentity} we have
\begin{equation*}
   \tilde P_{i, k d_i} = 
   E^{(k)}_{d_i \delta - \alpha_i}E^{(k)}_{\alpha_i} 
   + \sum_{0 < \boc_\pm'} a_{\boc'} L(\boc', 0).
\end{equation*}
Taking 
$\setbox5=\hbox{A}\overline{\rule{0mm}{\ht5}\hspace*{\wd5}}\,$ 
of both sides of this equality we have
\begin{align*}
   \ov{\tilde P_{i, k d_i}} 
   &=
   \ov{E^{(k)}_{d_i \delta - \alpha_i}}E^{(k)}_{\alpha_i}
   + \sum_{0 < \boc'_\pm} \ov{a_{\boc'} L(\boc',0)}
\\
  & = E^{(k)}_{d_i \delta - \alpha_i}E^{(k)}_{\alpha_i}
  + \sum_{\boc_{k,\pm} < \boc''_{\pm}} a_{\boc,\boc''} 
  L(\boc'',0) E^{(k)}_{\alpha_i} +
  \sum_{0 < \boc'_\pm} \ov{a_{\boc'} L(\boc',0)},
\end{align*}
where we have used {\bf (c)} in the second equality. By
\lemref{lem:reordering}, $L(\boc'',0) E^{(k)}_{\alpha_i}$ is a sum of
$L(\bod,0)$ with $\boc''_+ \le \bod_+$. (The other inequality $0 \le
\bod_-$ is trivially satisfied.) Since $(0 =) \boc_{k,+} < \boc''_+$,
we have $0 < \bod_+$.

Recall that by \propref{keyidentity} the purely imaginary part of
$L(\boc',0)$ is a polynomial in $\tilde P_{i,td_i}$ with $0 < t < k$.
By the induction hypothesis, \eqref{eq:barupt} is true for this
polynomial. Then it is also true for $L(\boc',0)$ by (b). Thus
\eqref{eq:barupt} is true for $\tilde P_{i,kd_i}$. By
\lemref{lem:reordering}, the assertion is also true for $S_{\boc_0}$
if it is a polynomial in $\tilde P_{i,td_i}$ with $t\le k$.
\end{proof}

When $\g$ is symmetric or type $\att$ we know that the set
$\{L(\boc,0) \mid \boc \in \aC\}$ is an $\ca$--basis of $\Ua^+$ from
\cite{BCP} and \cite{A}.  Using this result, we can obtain the more
general
\begin{Lemma}\label{lem:LpInt}
For each $p \in \Z$, the set $\{L(\boc,p) \mid \boc \in \aC \}$ is
an $\ca$--basis of $\Ua^+$.  
\end{Lemma}
\begin{proof}
First note that 
\begin{equation} \label{lzero}
\big\{ L(\boc,0) \mid \boc \in \aC, \boc(0) = 0 \big\}
\end{equation}
is an $\ca$--basis of $\Ua^+ \cap \U^+[i_0]$. To see this take any $x
\in \Ua^+ \cap \U^+[i_0]$.  Since the lemma holds for $p = 0$, we have
$x = \sum_{\boc' \in \aC} a_{\boc'} L(\boc', 0)$ with $a_{\boc'} \in
\ca$.  Now $x \in \U^+[i_0] \Leftrightarrow$ each $L(\boc',0) \in
\U^+[i_0]$ since $\U^+ = \U^+[i_0] \oplus E_{i_0} \U^+$ and clearly
each $L(\boc',0)$ is in one of these direct summands.  Consider the
image of \eqref{lzero} under the two maps $T_{i_0}$ and $*$
respectively.  We have
\begin{equation}
       T_{i_0} : \big\{ L(\boc,0) \mid \boc \in \aC, \boc(0) = 0 \big\}
\xrightarrow{1-1} \big\{ L(\boc,-1) \mid \boc \in \aC, \boc(0) = 0 \big\},
\end{equation}
as well as 
\begin{equation}
       * : \big\{ L(\boc,0) \mid \boc \in \aC, \boc(0) = 0 \big\}
\xrightarrow{1-1} \big\{ L(\boc,0)^* \mid \boc \in \aC, \boc(0) = 0
\big\}.
\end{equation}
Since both $*$ and $T_{i_0}$ leave $\Ua^+$ invariant, and both take
$\U^+[i_0]$ isomorphically to ${}^*\U^+[i_0]$ we obtain that the two
sets
\begin{equation}
\big\{ L(\boc,-1) \mid \boc \in \aC, \boc(0) = 0 \big\}, \quad 
\big\{ L(\boc, 0)^* \mid \boc \in \aC, \boc(0) = 0 \big\},
\end{equation}
are both $\ca$--bases of ${}^*\U^+[i_0] \cap \Ua^+$. 
This implies that the two sets
\begin{equation}
\big\{ L(\boc,-1) E_{i_0}^{(k)} \mid \boc \in \aC, \boc(0) = 0, k \in
\N \big\}, \quad \big\{ L(\boc, 0)^* E_{i_0}^{(k)} \mid \boc \in \aC,
\boc(0) = 0, k \in \N \big\},
\end{equation}
span the same $\ca$--submodules of $\Ua^+$. The right hand set is a
$\ca$--basis of $\Ua^+$ since it's the image under $*$ of a
$\ca$--basis of $\Ua^+$.  Therefore the left hand set is a
$\ca$--basis of $\Ua^+$, but this set is exactly $\big\{ L(\boc, -1)
\mid \boc \in \aC \big\}$.  With the same reasoning, an induction
gives the lemma for all $p < 0$.  A similar argument works for $p >
0$.
\end{proof}

{}From \propref{bartriangular} and \lemref{lem:LpInt} we deduce using
\cite[24.2.1]{Lu-Book}:

\begin{Theorem}\label{thm:PBWtoCan}  Let $\g$ be affine of symmetric
or of type $\att$.

\textup{(i)} For any $p\in \Z, \boc \in \aC$ there is a unique
$b(\boc,p) \in \Ua^+$ such that

\begin{aenume}
\item $\ov{b(\boc,p)} = b(\boc,p),$
\item $b(\boc,p) = L(\boc,p) + \sum_{\boc \prec_p \boc', \boc'
\in \aC} a_{\boc,\boc'} L(\boc',p),$ where $a_{\boc,\boc'} \in
q_s^{-1} \Z[q_s^{-1}]$.
\end{aenume}

\textup{(ii)} The $\Z$--homomorphism 
\[
   \La(-\infty) \cap \ov{\La(-\infty)}\cap\Ua^+ \rightarrow 
   \La(-\infty)\cap\Ua^+/q_s^{-1}(\La(-\infty)\cap\Ua^+)
\]
is an isomorphism. 
\end{Theorem}

For $p,q\in\Z$, there exists a bijection
$\boc\in\aC\leftrightarrow\boc'\in\aC$ such that
\(
   b(\boc,p) = \pm b(\boc',q)
\)
by the proof of \cite[14.2.3]{Lu-Book}. We will show $b(\boc,p) =
G(b)$ for some $b\in\B(-\infty)$ in \secref{sec:sign}, but it is
desirable to have a proof of $+$-sign in the above equality,
independent of the existence of the global crystal basis.

\section{Crystal structure of $\Um$} \label{sec:crystal}

Let $P_+^0 = \sum_{i\in I_0} \N \varpi_i$. Let $\lambda = \sum_i
\lambda_i \varpi_i \in P_+^0$. Let $G_\lambda = \prod_i
\GL_{\lambda_i}(\C)$ and $\Irr G_\lambda$ the set of irreducible
representations of $G_\lambda.$
Let $\breve B(\lambda)$ be the crystal of $\bigotimes_i
V(\varpi_i)^{\otimes \lambda_i}$.  In \cite[\S13]{Kas00} Kashiwara
conjectures a description of the crystal structure of $V(\lambda)$ in
terms of $\Irr G_\lambda$ and $\breve B(\lambda).$ This conjecture was
proved in \cite{Beck3,extrem} in the symmetric affine case. It can now
be checked for arbitrary type modulo sign using the results of the
previous section. The modification is straightforward, but we recall
the proof for the sake of reader. We also give a Peter--Weyl type
description of $\B(\Um)$ which is conjectured in \cite[\S13]{Kas00}
(see also \cite{Tos}), but not proved in \cite{Beck3,extrem}.

The sign ambiguity will be removed in \secref{sec:sign} based on
results in this section. Thus returning back to this section again, we 
get Kashiwara's conjecture.

\begin{rem} In the previous section we constructed an integral basis 
of $\Ua^+$.  In this section, for the purposes of calculation we
replace the basis elements $B_\boc, \boc \in \aC$ with $B_\boc^+ =
\ov{B^*_\boc}, \boc \in \aC$.  We also replace the $F_\alpha =
\Omega(E_\alpha), \alpha \in \aR^+$ which appeared in the previous
section, with $F^-_\alpha = \ov{E_\alpha^\vee}, \alpha \in \aR^+.$
Additionally, we define the integral basis of $\Ua^-$ given by
$B_\boc^- = \ov{B_\boc^\vee}$.
Let $\tilde P_{i,-kd_i} = \ov{\tilde P_{i,kd_i}^\vee}$. 
There are two reasons to do this.  First, applying the
$\setbox5=\hbox{A}\overline{\rule{0mm}{\ht5}\hspace*{\wd5}}\,$
operator allows us to work in $\La(\pm \infty)$ rather than
$\ov{\La(\pm \infty)}.$ Second, the operators $*$ and $\vee$ reverse
the root orderings in $B_\boc^\pm$ and so we are able to work with
highest weight (relative to the map $\cl$) level zero representations
instead of lowest weight level zero representations.
\end{rem}

\begin{comment}
The definition of $\tilde P_{i,-kd_i}$ does not fit with our conventions.
More consistent might be $\tilde P_{i,-kd_i}^-$.
\end{comment}

\begin{Definition} For $\lambda = \sum_i \lambda_i \varpi_i\in P^0_+$ let
$$\bocind = \{ \boc_0 = (\rho^{(1)}, \dots, \rho^{(n)}) \mid
\rho^{(i)} \text{ a
 partition, } \ell(\rho^{(i)}) \le \lambda_i, i = 1,
 \dots, n \}.$$
This is identified with the set of irreducible {\it polynomial\/}
representations of $G_\lambda$, and the set of $n$-tuples
$(s_{\rho^{(1)}}, \dots, s_{\rho^{(n)}})$, where $s_{\rho^{(i)}}$
is a Schur function in variables $z_{i,1}$, \dots, $z_{i, \lambda_i}$.
\end{Definition}

\subsection{Preparatory results}

\begin{prop} \label{iextremal}
\textup{(i)} For any $\lam \in P_+^0$, any vector in $\B(\lambda)$ is
connected to an extremal weight vector of the form $b_1 \otimes
t_\lambda \otimes u_{-\infty},$ where $b_1$ is purely imaginary with
respect to the crystal base. 

\textup{(ii)} Furthermore, all such possible $b_1 \in \B(\infty)$ are
given by $\sgn(\boc_0,0)S^-_{\boc_0} u_\infty \bmod q_s\La(\infty)$
where $\boc_0 \in \bocind.$
\end{prop}

\begin{proof} (i) By \cite[Proof of Theorem 5.1]{Kas00} any vector
is connected to an extremal weight vector of the form $b_1 \otimes
t_\lambda \otimes u_{-\infty},$ where $\wt(b_1) = -k\delta.$ Using the
basis of \S3 to express $b_1,$ we take $\sgn(\boc,0)B^-_\boc = b_1
\bmod q_s\La(\infty),$ for some $\boc$.  Assume that $B^-_\boc$ isn't
purely imaginary.  Since $\wt(B^-_\boc) = -k\delta$, and
$B^-_{\boc_+}$ (resp.\ $B^-_{\boc_-}$) consists only of terms in root
vectors with positive real part (resp.\ negative real part), it
follows $\boc_- \neq 0.$ By \thmref{thm:convex} we have
$B^-_\boc u_\lam = 0$.  However, by assumption $B^-_\boc u_\lam \in
\La(\lam)$ such that $B^-_\boc u_\lam \neq 0 \bmod
q_s\La(\lambda)$. This is a contradiction.

(ii) From Proposition~\ref{keyidentity} we have 
\begin{equation}\label{eq:tP}
   \tilde P_{i,-kd_i} u_{\lambda}
   = 
   \begin{cases}
    F_{d_i \delta - \alpha_i}^{-\ (k)} F_{i}^{(k)}u_{\lambda}
    & \text{if } (\xn,i) \neq (\atnt,n), \\
    F_{\delta -  2 \alpha_n}^{-\ (k)} F_{n}^{(2k)}u_{\lambda}
    &  \text{if } (\xn,i) = (\atnt,n).
   \end{cases}
\end{equation}
Here $u_\lambda$ generates $V(\lambda)$. Since the weights of
$V(\lambda)$ are in the convex hull of $W \lambda$ (see
\thmref{thm:convex}), this implies that $\tilde P_{i,-kd_i} u_\lambda
= 0$ for $k > \lambda_i$.  Note that for any $i$, $\ell(\rho^{(i)})
\le \lambda_i \iff {\rho^{(i)}}^t_1 \le \lambda_i.$ Since the $\tilde
P_{i,-kd_i}$ all commute, considering the top row of the determinant
$S^-_{\boc_0}$, we have $S^-_{\boc_0} u_\lambda = 0$ for $\boc_0
\notin \bocind.$
\end{proof} 

Let $\tal_i\in\tQ$ as in \subsecref{subsec:affine} and let
$S_{\tal_i}$ be the corresponding operator in
Definition~\ref{def:extremal}.

\begin{Lemma}\label{lem:Stal}
Let $\lambda = \sum_i \lambda_i\varpi_i \in P_+^0$.  Let $V(\lambda)$
be the extremal weight module generated by $u_\lambda$.
\begin{equation*}
\begin{split}
S_{\tal_i} u_{\lambda} &= 
\begin{cases}
 F_{d_i \delta - \alpha_i}^{-\ (\lambda_i)} F_{\alpha_i}^{(\lambda_i)}
 u_{\lambda} & \text{if } (\xn,i) \neq (A_{2n}^{(2)},n)
\\
 F_{\delta - 2 \alpha_n}^{-\ (\lambda_n)} F_{\alpha_n}^{(2\lambda_n)}
 u_{\lambda} & \text{if } (\xn,i) = (A_{2n}^{(2)},n)
\end{cases}
\\
 &= \tilde P_{i,-\lambda_i d_i} u_{\lambda}
\end{split}
\end{equation*}
\end{Lemma}

\begin{proof}
The second equality follows from the first equality and \eqref{eq:tP}.
We prove the first equality after applying
\(
\setbox5=\hbox{A}\overline{\rule{0mm}{\ht5}\hspace*{\wd5}}\,
   \circ \vee.
\)
(See the above remark.)

First consider the case $(\xn,i) \neq (A_{2n}^{(2)},n)$.
We have an identity $\tal_i = s_{d_i\delta-\alpha_i} s_i = \tom_i s_i
\tom_i^{-1} s_i$ in the affine Weyl group (see \eqref{eq:WeylId}),
where $s_{d_i\delta-\alpha_i}$ is the reflection with respect to
$d_i\delta-\alpha_i$. Then
\begin{align*} 
   S_{\tal_i} u_{-\lambda} &
   = S_{\tom_i} S_i S_{\tom_i^{-1}} S_i u_{-\lambda}
   = S_{\tomp_i} S_i S_{{\tilde \omega^{\prime -1}_i}}
      E_{\alpha_i}^{(\lambda_i)} u_{-\lambda}
      \text{ where $\tomp_i = \tom_i s_i$ }
\\
   & = S_{\tomp_i} E_{\alpha_i}^{(\lambda_i)}
   S_{{\tilde \omega^{\prime -1}_i }}
   E_{\alpha_i}^{(\lambda_i)} u_{-\lambda}
  = q^{-N_+} (-1)^{N^\vee_+} T_{\tomp_i} (E_{\alpha_i}^{(\lambda_i)})
    T_{\tomp_i} (S_{{\tilde \omega^{\prime -1}_i}}
    E_{\alpha_i}^{(\lambda_i)} u_{-\lambda})
\\
   & = q^{M_+} (-1)^{M^\vee_+} q^{-N_+} (-1)^{N^\vee_+} T_{\tomp_i}
  (E_{\alpha_i}^{(\lambda_i)}) S_{\tomp_i}
  S_{{\tilde \omega^{\prime -1}_i}} E_{\alpha_i}^{(\lambda_i)} u_{-\lambda}
\\
   & = q^{M_+} (-1)^{M^\vee_+} q^{-N_+} (-1)^{N^\vee_+}
   E_{d_i \delta - \alpha_i}^{(\lambda_i)}
   E_{\alpha_i}^{(\lambda_i)} u_{-\lambda}
\end{align*}
by \lemref{lem:WeylExtremal}. Here
\begin{align*}
\textstyle M_+  = \sum_{\alpha
\in \Delta^+ \cap \tom_i^{\prime -1}(\Delta^-)} & \max(-(\alpha,
\tom_i^{\prime -1}s_i(\lambda)), 0)
\\
 & \text{ and } N_+ = \textstyle \sum_{\alpha
\in \Delta^+ \cap \tom_i^{\prime -1}(\Delta^-)} \max(-(\alpha,
s_i\tom_i^{\prime -1}s_i(\lambda)), 0)
\end{align*}
and $M^\vee_+$, $N_+^\vee$ are defined by replacing $\alpha$ by
$\alpha^\vee$. We have used Remark~\ref{rem:loopgen} in the last
equality.

We set $\alpha' = -\tom_i' \alpha$. Then $\alpha'\in\Delta^+\cap
\tom_i'\Delta^-$ and
\begin{equation*}
\begin{gathered}
   -(\alpha, s_i \tom_i^{\prime -1}s_i(\lambda)) = (\alpha', \tal_i\lambda)
   = (\alpha', \lambda),
\\
   -(\alpha, \tom_i^{\prime -1}s_i(\lambda)) = (\alpha', s_i \lambda).
\end{gathered}
\end{equation*}

For $\beta\in\Delta_\cl$, let us denote by $\beta'$ the unique element
of $\Delta^+$ such that $\cl(\beta') = \beta$ and
$\beta'-n\delta\notin\Delta^+$ for any $n > 0$. We have
\begin{equation*}
    \Delta^+ \cap \tom_i \Delta^-
    = \{ \beta' + n d_\beta \delta \mid \beta\in\Delta_{\cl},
    n\in\Z, 0\le n < -(\tom_i, \beta)/d_\beta \}.
\end{equation*}
Therefore,
\begin{align*}
 \Delta^+ \cap \tom'_i(\Delta^-)
  &
= \Delta^+ \cap \tom_i(\Delta^-) \setminus \{ \tom_i'(\alpha_i) \}
\\
& = \{\beta' + n d_\beta \delta \mid
   \beta \in \Delta_\cl,\ n\in\Z,\ 
   0 \le n < -(\tom_i, \beta)/d_\beta \} \setminus
   \{ d_i\delta - \alpha_i \}
,
\end{align*}
where we have used $(\tom_i, \cl(\alpha_i)) = d_i$. If
$\beta\in\Delta_\cl^+$, then $(\tom_i, \beta)\ge 0$ and there are no
corresponding terms in the summation of $M_+$, $N_+$. If
$\beta\in\Delta_\cl^-$, then $(\alpha',\lambda)\le 0$ for $\alpha' =
\beta'+nd_\beta\delta$ and $(\alpha', s_i\lambda)\le 0$ except
possibly when $\beta=-\cl(\alpha_i)$, i.e., $\alpha' = -\alpha_i +
nd_i\delta$. However, there are no such roots in $\Delta^+ \cap
\tom'_i(\Delta^-)$. Therefore $M_+ = N_+ = 0$. By the same reason, we
also have $M^\vee_+ = N^\vee_+ = 0$.

Next consider the case $(\xn,i) = (A_{2n}^{(2)},n)$. We have an identity
$\tal_n = s_{\delta-2\alpha_n} s_n$. Following \cite[\S4.2]{Da2}, we set
\begin{equation*}
    w = s_0 s_1 \dots s_n \in \aW.
\end{equation*}
(Our numbering is different from one in [loc.~cit.]) We have
$w^{n-1}(\alpha_0) = \delta-2\alpha_n$. We can repeat the above
calculation replacing the identity by
\(
  \tal_n = w^{n-1} s_0 w^{1-n} s_n.
\)
We have
\begin{equation*}
   S_{\tal_n} u_{-\lambda} = 
   q^{M_+} (-1)^{M^\vee_+} q^{-N_+} (-1)^{N^\vee_+} 
   T_{w^{n-1}}(E_{\alpha_0}^{(\lambda_n)})
   E_{\alpha_n}^{(2\lambda_n)} u_{-\lambda},
\end{equation*}
where
\begin{align*}
   & M_+  = \sum_{\alpha'\in \Delta^+ \cap w^{n-1}(\Delta^-)}
   \max((\alpha', s_n(\lambda)), 0),
\\
   & N_+ = \sum_{\alpha'\in \Delta^+ \cap w^{n-1}(\Delta^-)}
     \max((\alpha', \lambda), 0),
\end{align*}
and $N_+^\vee$, $M_+^\vee$ are similar.

We have $\Delta^+\cap w^{n-1}(\Delta^-)\cap \cl^{-1}(\Delta_\cl^+) =
\emptyset$. Therefore the above summations are over
$\alpha'\in\Delta^+\cap w^{n-1}(\Delta^-)\cap \cl^{-1}(\Delta_\cl^-)$.
Then $(\alpha',\lambda)\le 0$ for any such $\alpha'$, and
$(\alpha',s_n(\lambda)) \le 0$ except possibly when
$\cl(\alpha')=-\cl(\alpha_n)$ or $-2\cl(\alpha_n)$, i.e., $\alpha' =
m\delta-\alpha_n$ or $(2m-1)\delta-2\alpha_n$ for some
$m\in\Z_{>0}$. However we have
\begin{equation*}
   w^{1-n}(\alpha_n) = \alpha_1 + \dots + \alpha_n
\end{equation*}
by \cite[4.2.3]{Da2}, which means that such $\alpha'$ cannot be in
$\Delta^+\cap w^{n-1}(\Delta^-)$. Therefore $M_+ = M_+^\vee = N_+ =
N_+^\vee = 0$.

Finally we have $T_{w^{n-1}}(E_{\alpha_0}) = E_{\delta-2\alpha_n}$ by
\cite[4.2.6]{Da2}.
\end{proof}

Let $z_i$ be the $\U'$-module automorphism of $V(\varpi_i)$ defined in
\S2.6.
\begin{lem}\label{piaction} Let $i\in I_0$. Then on $V(\varpi_i)$:
\begin{equation*}
\tilde P_{i,-d_i} u_{\varpi_i} = \overline{\tilde P_{i,-d_i} u_{\varpi_i}}
= z_i^{-1} u_{\varpi_i}, \qquad
\tilde P_{i,-kd_i} u_{\varpi_i} = 0, \quad \text{for $k > 1$}.
\end{equation*}
\end{lem}

\begin{proof}
The statement for $k > 1$ was already observed in the proof of
\propref{iextremal}(ii). \lemref{lem:Stal} (with $\lambda = \varpi_i$,
$\lambda_i = 1$) says $\tilde P_{i,-d_i} u_{\varpi_i} = S_{\tal_i}
u_{\varpi_i}$. In particular, $\tilde P_{i, - d_i} u_{\varpi_i}$ is an
element of the global basis.
On the other hand, $V(\varpi_i)$ has a unique global basis element of
weight $\varpi_i - d_i\delta$ (see \thmref{thm:fund}), which by
definition equals $S_w u_{\varpi_i} = z_i^{-1} u_{\varpi_i},$ where $w
= \tal_i.$ The assertion follows.
\end{proof}

\subsection{The map $\Phi_\lambda$ and Kashiwara's conjecture}

Let $\lam = \sum_{i\in I_0} \lambda_i\varpi_i \in P^0_+$.  The module
$$\Vtl = \bigotimes_{i\in I_0}V(\varpi_i)^{\otimes \lambda_i}$$ has a
crystal base $(\bigotimes_i \La(\varpi_i)^{\otimes \lambda_i},
\bigotimes_i \B(\varpi_i)^{\otimes \lambda_i}).$ Let $\utl =
\mathop\bigotimes_i u_{\varpi_i}^{\otimes \lambda_i}.$
 
For each $i$, and each of the $\nu = 1, \dots, \lambda_i$ factors of
$V(\varpi_i)^{\otimes \lambda_i}$, we let $z_{i,\nu}$ be the commuting
automorphisms defined in \cite[\S 4.2]{Kas00}. By
\cite[Theorem~8.5]{Kas00}, the submodule
\begin{equation} \label{defW}
   \Wtl = \Ua[z_{i,\nu}^{\pm 1}]_{1 \le i \le n, 1 \le
    \nu \le \lambda_i} \utl \subset \Vtl
\end{equation}
has a global crystal basis $(\breve \La, \breve \B)$ such
that 
\(
   \breve \La \subset \bigotimes_i \La(\varpi_i)^{\otimes \lambda_i},
\)
\(
  \breve \B = \bigotimes_i \B(\varpi_i)^{\otimes \lambda_i}.
\)
Since $\Wtl$ contains the extremal vector $\utl$ of weight $\lam$ we
have a unique $\U$-linear morphism
\begin{equation} \label{mainmap} \Phi_\lam\colon V(\lam)\to \Wtl, 
\end{equation} 
sending $u_\lam$ to $\utl,$ and which commutes with the crystal
operators $\tilde e_i, \tilde f_i$.

\begin{comment}
In earlier versions, $\Wtl$ (defined by the macro \verb+\Wtl+) was
denoted by $\breve W(\lambda)$. We hope there are no remaining
mistakes caused by the changes.
\end{comment}

For each $n$--tuple of partitions $\boc_0 = (\rho^{(1)},
\rho^{(2)}, \dots, \rho^{(n)})$ we consider the product of Schur
functions in the variables $z^{\pm 1}_{i,\nu}$ (see \cite{Mac}):
\begin{equation}
s_{\boc_0}(z_{i,\nu}^{\pm 1}) = \prod_{i\in I_0}
s_{\rho^{(i)}}(z_{i,1}^{\pm 1}, \dots, z_{i,\lambda_i}^{\pm 1}).
\end{equation}
Note that for each $i$, $s_{\rho^{(i)}}(z^{\pm 1}_{i,\nu})$ acts as
the $0$ map if $\lambda_i < \ell(\rho^{(i)}).$ We will omit the indices
$i,\nu$ and write $s_{\boc_0}(z^{\pm 1}).$

Using Lemma~\ref{piaction} we have:
\begin{prop} \label{schurcor} 
Let $\boc_0 = (\rho^{(1)},\rho^{(2)}, \dots, \rho^{(n)})$ be an
$n$--tuple of partitions:
\begin{equation} \label{schur}
\Phi_\lam(S^{-}_{\boc_0}\, u_\lambda) =
            s_{\boc_0}(z^{-1})\, \utl.
\end{equation}
\end{prop}

\begin{proof} 
Note that $\sigma \circ (^\vee \times ^\vee) \circ \Delta(a) =
\Delta(a^\vee)$ for $a \in \U$. Here $\sigma$ is the exchange of two
factors of the tensor product. Since our $\tilde{P}_{i,-kd_i}$ $(k >
0)$ are those given in \secref{sec:intbase} after applying $- \circ
\vee$ we have by \propref{prop:coprod} and \cite[3.1.5]{Lu-Book}
(after noting the difference between our coproduct and the one there)
\begin{multline*}
\Delta (\tilde{P}_{i,-kd_i}) =\sum_{s=0}^{k}
\tilde{P}_{i,-sd_i} \ot \tilde{P}_{i,(s-k)d_i}
\\
 + {\text{terms acting as }} 0 \text{ on } v_{\varpi_{j_1}} \otimes
v_{\varpi_{j_2}} \text{ for } j_1, j_2 \in I_0.
\end{multline*}
This and Lemma~\ref{piaction} imply that
$\Delta^{\lambda_i}(\tilde P_{i,-kd_i})$ acts as
$e_k(z_{i,1}^{-1}, \dots, z_{i,\lambda_i}^{-1})$ on $\Vtl$ where $e_k$
is the $k$--th elementary symmetric function.  Since polynomials in
the $\tilde P_{i,-kd_i}$ (resp.\ elementary symmetric functions)
generate the Schur functions $S^{-}_{\boc_0}$ (resp.\
$s_{\boc_0}(z^{-1})$), we have $\Phi_\lam(S^{-}_{\boc_0}\, u_\lambda) = 
s_{\boc_0}(z^{-1})\, \utl.$
\end{proof}

Next we consider the image of $\B(\lambda)$ under $\Phi_\lambda.$ By
\propref{iextremal} every element of $\B(\lam)$ is connected to
an extremal vector of the form
\(
   \sgn(\boc_0,0)S^-_{\boc_0}u_\infty
   \otimes t_\lambda \otimes u_{-\infty} \bmod q_s \La(\lam).
\)
Therefore we have,
\begin{multline} \label{mvtoimaginary}
  \B(\lambda) = \{X_1\, X_2\, \cdots\, X_n 
   (\sgn(\boc_0,0)S^-_{\boc_0} 
   u_\lambda\bmod q_s \La(\lambda)) \mid 
\\
   X_i = \tilde e_i \text{ or } \tilde f_i, \boc_0
   \in \bocind \} \setminus \{0\}.
\end{multline}
Since $\Phi_{\lambda}$ commutes with crystal operators, and the
$z_{i,\nu}$ induce automorphisms of the $\U'$--crystal of
$V(\varpi_i)$, we have $\Phi_\lambda(\La(\lambda)) \subset \breve \La$
by \propref{schurcor}. Denote by $\Phi_{\lambda | q=0}$ the induced
map $\La(\lam)/q_s\La(\lam) \rightarrow \breve \La/q_s\breve \La$.

\begin{prop} \label{crystalcorresp} Let $\BWzero$ be the connected
component of $\breve \B$ containing $\utl.$ Then
$$
   \Phi_{\lambda | q=0}\colon \{ b \mid b \in \B(\lambda) \}
   \rightarrow
   \{\sgn(\boc_0,0) s_{\boc_0}(z^{-1}) b'\mid
    \boc_0 \in \bocind,\, b' \in \BWzero \}
$$
is a bijection.
\end{prop}

\begin{proof}
It is clear that $\Phi_{\lambda | q=0}(\B(\lambda)) \setminus \{0\}$
is equal to the right hand side of the above.
Using \eqref{mvtoimaginary}, we check that $\Phi_{\lambda |
q=0}(\B(\lambda))$ does not contain $0$ and is injective.  Let $b \in
\B(\lambda)$ such that $\Phi_{\lambda | q=0} (b) = 0$. Since $b$ is
connected by crystal operators to $b_1 \otimes t_\lambda \otimes
u_{-\infty},$ where $b_1 = \sgn(\boc_0,0) S^-_{\boc_0} u_\infty \bmod
q_s\La(\infty)$, $\boc_0 \in \bocind$, this implies $\Phi_{\lambda |
q=0} (S^-_{\boc_0} u_\lambda \mod q_s\La(\lambda)) = 0$. This
contradicts Proposition \ref{schurcor}.  Now let $b_1, b_2 \in
\B(\lambda), b_1 \neq b_2$ and assume $\Phi_{\lambda | q=0}
(b_1)=\Phi_{\lambda | q=0} (b_2)$.  By applying crystal operators we
may assume $b_1 = \sgn(\boc_0,0)S^-_{\boc_0} u_\lambda \bmod q_s
\La(\lambda)$ is extremal and purely imaginary, and where $b_2$ is of
the same weight as $b_1$.  For $w \in W$, we have
$$\Phi_\lambda|_{q=0}(X_i S_w b_1) = X_i S_w \Phi_\lambda|_{q=0}(b_1)
= X_i S_w \Phi_\lambda|_{q=0}(b_2) = \Phi_\lambda|_{q=0}(X_i S_w b_2)
= 0,$$ where $X_i$ is $\te_i$ or $\tf_i$.  Since
$\Ker\,\Phi_\lambda|_{q=0}\cap \B(\lambda) = \emptyset$, we have $X_i
S_w b_2 = 0$. Therefore, $b_2$ is also extremal.  Applying a sequence
of ${\tf_i}^{\max}$'s (since $b_2$ is extremal, this is equivalent to
the Weyl group action) we may assume $b_2 =
\sgn(\boc',0)S^-_{\boc^\prime_0} u_\lambda \bmod q_s \La(\lambda)$ (as
in \cite[Theorem 5.1]{Kas00}). Then we have that $\Phi_\lambda(b_1) =
\Phi_\lambda(b_2) \neq 0$ and also that $b_1 =
\sgn(\boc'',0)S_{\boc''_0}^- u_\lambda \mod q_s \La(\lambda)$ for some
$\boc''_0$ purely imaginary in $\aC$ by \propref{iextremal} (ii).  But
by Proposition \ref{schurcor} it is clear that
$\Phi_\lambda(\sgn(\boc'_0,0)S^-_{\boc'_0} u_\lambda) =
\Phi_\lambda(\sgn(\boc''_0,0)S^-_{\boc''_0} u_\lambda)$
only if $S^-_{\boc'_0} u_\lambda = S^-_{\boc''_0} u_\lambda.$
\end{proof}

\begin{Lemma}\label{lem:Wlam}
\textup{(i)} Let $\boc_0$, $\boc_0'\in \bocind$ and $b$, $b'\in
\BWzero$. Then $s_{\boc_0}(z^{-1}) b =\linebreak[0] s_{\boc'_0}(z^{-1}) b'$
if and only if $s_{\boc'_0}(z^{-1}) = s_{\boc_0}(z^{-1}) p(z)$ and $b'
= p(z)^{-1} b$ for some $p(z) = \prod_i (z_{i,1}\cdots\linebreak[0]
z_{i,\lambda_i})^{r_i}$ \textup($r_i\in\Z$\textup).

\textup{(ii)} For $b\in\BWzero$ and $p(z)$ as above, we have
$p(z)b\in\BWzero$.
\end{Lemma}

\begin{proof}
(i) Recall that for fundamental representations $V(\varpi_i)$, we have
an isomorphism $V(\varpi_i)\cong \Q(q_s)[z_i^\pm]\otimes W(\varpi_i)$,
where $W(\varpi_i)$ is the finite dimensional representation (see
\subsecref{subsec:extremal}).
For general $\lambda = \sum_i \lambda_i \varpi_i\in P^0_+$, we
put
\[
   W(\lambda) = \bigotimes_{i\in I_0} W(\varpi_i)^{\otimes \lambda_i}.
\]
This is irreducible and has a global crystal basis $\B_W(\lambda)$. We 
have
\begin{equation*}
\begin{gathered}
  \Vtl \cong
   \Q(q_s)[z_{i,\nu}^\pm]_{1 \le i \le n, 1 \le \nu \le \lambda_i}
   \otimes W(\lambda),
\\
  \breve \B \cong
   \{ \text{monomials in $z_{i,\nu}^\pm$} \} \times \B_W(\lambda).
\end{gathered}
\end{equation*}
Therefore
\begin{equation*}
   s_{\boc_0}(z^{-1})b = s_{\boc_0}(z^{-1}) m(z) b_W,
\quad
   s_{\boc'_0}(z^{-1})b' = s_{\boc'_0}(z^{-1}) m'(z) b'_W,
\end{equation*}
for some monomials $m$, $m'$ and $b_W$, $b'_W\in \B_W(\lambda)$. These 
two are equal if and only if $b_W = b'_W$ and $s_{\boc_0}(z^{-1})m(z)
= s_{\boc'_0}(z^{-1}) m'(z)$. Then $p(z) = m(z)/m'(z)$ is a symmetric
function which is also a monomial. Therefore it must be of the above
form.

(ii) Since $p(z)$ commutes with $\te_i$, $\tf_i$, we may assume $b =
\utl$. For $w \in \aW$, let $S_w$ denote the corresponding crystal
operator. Then we have
\begin{equation*}
S_{w} \utl = \bigotimes_i (S_w u_{\varpi_i})^{\otimes \lambda_i}
\end{equation*}
by \cite[Lemma 1.6]{AK} since $\B(V(\varpi_i))$ is the affinization of
a simple crystal. Let us take $w = t(-\sum_i r_i \tal_i)$. By
Lemmas~\ref{lem:Stal}, \ref{piaction} we get
\begin{equation*}
   S_w \utl = \bigotimes_{i,\nu} (z_{i,\nu})^{r_i} \utl
   = p(z) \utl.
\end{equation*}
Thus $p(z)\utl$ is connected to $\utl$ in the crystal graph.
\end{proof}

Let
\[
  \bocindp
  = \{(\rho^{(1)}, \dots, \rho^{(n)}) \mid \rho^{(i)} \text{ a
 partition, } \ell(\rho^{(i)}) < \lambda_i, i = 1,
 \dots, n \}.
\]
This can be identified with the set of irreducible representations of
$\prod_i \SL_{\lambda_i}(\C)$. Any $\boc_0\in\bocind$ decomposes
uniquely as
\(
   s_{\boc_0}(z^{-1}) = s_{\boc'_0}(z^{-1}) p(z)
\)
with $\boc'_0\in\bocindp$ and a monomial $p(z)$ as in the above
lemma. Therefore \propref{crystalcorresp} can be strengthened as
\begin{equation*}
\begin{split}
  \Phi_{\lambda | q=0}\colon \B(\lambda)
   \isoto
   \{ & \sgn(\boc_0,0) s_{\boc_0}(z^{-1}) b'\mid
    \boc_0 \in \bocind',\, b' \in \BWzero \}
\\
   & \cong \bocind'\times \BWzero.
\end{split}
\end{equation*}
A connected component of $\B(\lambda)$ is mapped to $\{ \boc_0
\}\times \BWzero$ for some $\boc_0\in\bocind'$.  In particular, each
connected component is isomorphic to each other as a $P_\cl$-crystal.

\begin{cor} \label{cor:inj} The map $\Phi_\lambda\colon 
V(\lambda) \rightarrow \Wtl$ is injective.
\end{cor}
\begin{proof} 
Since $\Phi_{\lam | q=0}\colon \La(\lam)/q_s\La(\lam) \rightarrow
\breve \La/q_s \breve \La$ maps the crystal basis $\B(\lam)$
bijectively to $\{ \sgn(\boc_0,0) s_{\boc_0}(z^{-1}) b'\mid
    \boc_0 \in \bocind',\, b' \in \BWzero \}$, which is linearly
independent. Therefore $\Phi_{\lambda | q=0}$ is injective.
Write an element $v \in \ker \Phi_\lambda, v \neq
0,$ in terms of the global basis $\{G(b)\mid b \in \B(\lambda)\}$ as
$v = \sum_b c_b(q_s) G(b).$ Multiplying by a power of $q$ we may assume
that each $c_b(q_s)$ is regular at $q=0,$ so that $v \bmod q_s \La(\lambda)
\neq 0$. This implies $\Phi_{\lam | q=0} (v \bmod q_s \La(\lam)) \neq 0$,
which is a contradiction.
\end{proof}

\begin{comment}
The linear independence is clear now: If
\begin{equation*}
   \sum_{(\boc_0,b')\in\bocind'\times \BWzero}
   a_{\boc_0,b'} s_{\boc_0}(z^{-1}) b' = 0,
\end{equation*}
then
\begin{equation*}
   \sum_{b'} a_{\boc_0,b'} b' = 0
\end{equation*}
for every $\boc_0$. The linear independence of $\BWzero$ implies
$a_{\boc_0,b'} = 0$.
\end{comment}

We now state a main result in this subsection:
\begin{Theorem}\label{thm:crystaldesc}
\textup{(i)} We have an isomorphism of crystals
\begin{equation*}
\begin{split}
  \Phi_{\lambda | q=0}\colon \B(\lambda)
   \isoto
   \{ & \sgn(\boc_0,0) s_{\boc_0}(z^{-1}) b'\mid
    \boc_0 \in \bocind',\, b' \in \BWzero \}
\\
   & \cong \bocind'\times \BWzero.
\end{split}
\end{equation*}
A connected component of $\B(\lambda)$ is mapped to $\{ \boc_0
\}\times \BWzero$ for some $\boc_0\in\bocind'$. Also, any two connected
components are isomorphic to each other as $P_\cl$-crystals.

\textup{(ii)} $\Phi_\lam$ induces a bijection between the sets
\begin{equation*}
  \Phi_\lam \colon G(\B(\lam))
   \rightarrow \{ \sgn(\boc_0,0) s_{\boc_0}(z^{-1}) G(b')
    \mid \boc_0 \in \bocind', \, b' \in \BWzero \}.
\end{equation*}

\textup{(iii)} A vector $b\in \B(\lambda)$ is extremal if and only if
\[
   \Phi_{\lambda | q=0}(b) = 
   \sgn(\boc_0,0) s_{\boc_0}(z^{-1})
   S_{w} \utl
\]
for some $\boc_0\in\bocindp$, $w\in\aW$.
\end{Theorem}

\begin{proof}
(i) is proved already. Let us show (ii). By Proposition
\ref{crystalcorresp}, for each $b \in \B(\lam)$ there exist $b' \in
\BWzero$ and $s_{\boc_0}(z^{-1})$ such that $\Phi_\lam (b) =
s_{\boc_0}(z^{-1}) b' \bmod q_s \breve \La$.  Let $G(b)$, $G(b')$ be
the respective globalizations of $b$ and $b'$.  Then $\Phi_\lam (G(b))
\equiv s_{\boc_0}(z^{-1}) G(b') \mod q_s \breve \La$.  Since
$\Phi_\lam$ commutes with the bar involution, we have $\Phi_\lam
(G(b)) \equiv s_{\boc_0}(z^{-1}) G(b') \mod q_s^{-1} \ov{\breve \La}$.
We get the assertion.

Let us prove (iii). It is enough to consider the case when
$\Phi_{\lambda | q=0}(b)\in \BWzero$. As in the proof of
\lemref{lem:Wlam}, we can write $\Phi_{\lambda | q=0}(b) = p(z) b_W$
where $p(z) = \prod_i (z_{i,1}\cdots z_{i,\lambda_i})^{r_i}$
($r_i\in\Z$) and $b_W\in \B_W(\lambda)$. Then $b$ is extremal if and
only if $b_W$ is so. Furthermore, $b_W$ is extremal
if and only if $b_W = S_{w_0} \utl$ for $w_0 \in W_\cl.$
This follows from \cite[Lemma 1.6]{AK}, after noting that
(\cite[Theorem~5.15]{Kas00}) $W(\varpi_i)$ has a simple crystal such
that the weight of any extremal vector belongs to
$\aW\,\cl(\varpi_i)$.
Now we have $\Phi_{\lambda | q=0}(b) = p(z) b_W = S_w \utl$ for some
$w\in\aW$ as in the proof of \lemref{lem:Wlam}.
\end{proof}

\begin{rem} Taken together the results of this section and
$\sgn(\boc_0,0) = 1$, which will be proved in the next section, give the
conjectures \cite[13.1, 13.2]{Kas00}. To obtain 13.1 (iii), consider that
the crystal $\bigotimes_{i\in I_0} \B(\lambda_i\varpi_i)$ is by Proposition
\ref{crystalcorresp} in bijective correspondence with
$\{s_{\boc_0}(z^{-1}) \BWzero\}$, and note that $\Phi_\lambda$ factors
through $\bigotimes_{i\in I_0} V(\lambda_i\varpi_i)$.
\end{rem}

\subsection{A Peter--Weyl type decomposition theorem}

Let $\B_0(\lambda)$ be the connected component of $\B(\lam)$ containing
$u_\lambda$. 
Consider $\bigsqcup_{\lam\in P}\B_0(\lam)\times \B(-\lam)$ as a
crystal over $\g\oplus\g$. Here for $u\ot v\in \B_0(\lambda)\otimes
\B(-\lambda)$, $X_i(u\ot v)=X_i u\ot v$ and $X_i^*(u\ot v)=u\ot
X_i v$ where ($X_i=\tilde e_i,\tilde f_i$).  The Weyl group $\aW$
acts on $\bigsqcup_{\lam\in P}\B_0(\lam)\times \B(-\lam)$ by
$S^*_w\times S^*_w\colon \B_0(\lam)\times \B(-\lam)\to
\B_0(w\lam)\times \B(-w\lam)$.

Consider the map $\B_0(\lam)\times \B(-\lam) \rightarrow
\B(\U a_\lambda)$, which sends $u_\lam\otimes b\in \B_0(\lam)\times
\B(-\lam)$ to $b^*\in \B(\U a_\lambda)$. We will show that this is
well-defined later. It is a map between crystals
over $\g\oplus\g$ where the usual crystal structure on $\B(\tU)$
corresponds to the one on $\B_0(\lam)$ and the star crystal structure
on $\B(\tU)$ corresponds to the one on $\B(-\lam)$. We have

\begin{thm} \label{PWdecomp} 
$\Bigl(\bigsqcup_{\lam\in P}\B_0(\lam)\times \B(-\lam)\Bigr)/\aW
\isoto \ \B(\tU)$ as crystals over $\g \oplus \g$.
\end{thm}

\begin{proof}
When $\langle c, \lambda \rangle \neq 0$, the $\B(\U a_\lambda)$ part
of the crystal decomposition in Theorem \ref{PWdecomp} appears as
\cite[Proposition 10.2.2]{Kas94}.  Therefore, it is sufficient to check the
map in Theorem \ref{PWdecomp} for $\lambda \in P^0$, where the
image on the right hand side is in $\B(\U a_\lambda).$

The following proof is a modification of a result of Nakashima
\cite[Proposition~4.4]{Tos}. We give it for the sake of the reader.

We first prove
\begin{description}
\item[(C1)] For any extremal vector $b\in \B(\tU)$, there exists
a crystal embedding:
$$
\B_0(\wt(b))\hookrightarrow \B(\tU),
$$
given by $u_{\wt(b)}\mapsto b$.
\end{description}

Let $b \in \B(\tU)$ be extremal of weight $\lambda$ and $B'$ be
connected component containing $b$. We want to prove $B'\cong
\B_0(\wt(b))$. By \cite[Corollary 9.3.4]{Kas94} the connected
component $B'$ of $b$ can be embedded into some $\B(\mu)$ for some
$\mu \in P^0.$ For each $w \in \aW$, we have $S_w^* (\B(\mu)) \isoto \ 
\B(w\mu)$.  Using this Weyl group action, and noting that for
$\cl(\mu) = \cl(\mu'), \B(\mu) \isoto \ \B(\mu')$ as $P_\cl$ crystals,
we may assume $\mu = \sum_{i=1}^{n} m_i \varpi_i \in P_0^+$. By
\thmref{thm:crystaldesc}(iii), (i), we have that $\Phi_\mu(b) =
\sgn(\boc_0,0) s_{\boc_0}(z^{-1}) S_w \tilde u_\mu$ for fixed $\boc_0
\in \aC, w \in \aW$, and $B'$ is isomorphic to $B_0(\mu)$ as a
$P_\cl$--crystal, so that $b$ is mapped to $S_w u_\mu$. This is
further isomorphic to $B_0(\lambda)$ so that $b$ is mapped to
$u_\lambda$. But since $\wt(b) = \wt(u_\lambda) = \lambda$ this is
also an isomorphism of $P$--crystals. This completes the proof of
$\bold{(C1)}$.

We now define a map
\(
   \varphi\colon \bigsqcup_\lambda \B_0(\lambda)\times \B(-\lambda) \to 
   \B(\tU)
\)
by
\begin{equation*}
    \varphi(X_{1} X_{2} \cdots X_{N} u_\lambda\otimes b)
    = X_{1} X_{2} \cdots X_{N} b^*,
\end{equation*}
where $X_j = \te_i$ or $\tf_i$. By $\bold{(C1)}$ this is well-defined,
i.e., (a) if $X_{1} X_{2} \cdots X_{N} u_\lambda = 0$, then $X_{1}
X_{2} \cdots X_{N} b^* = 0$, and (b) if $X_{1} X_{2} \cdots X_{N}
u_\lambda = X'_{1} X'_{2} \cdots X'_{N'} u_\lambda$, then $X_{1} X_{2}
\cdots X_{N} b^* = X'_{1} X'_{2} \cdots X'_{N'} b^*$. $\bold{(C1)}$ is
applicable since $b^*$ is extremal of weight $\lambda$ for $b\in
\B(-\lambda)$.

Since $\te_i^*$, $\tf_i^*$ commute with $\te_j$, $\tf_j$ on
$\B(\tU)$ (\cite[Theorem~5.1.1]{Kas94}), $\varphi$ is a morphism of
bi-crystal. Since any connected component of $\B(\tU)$ contains an extremal vector (\cite[Corollary~9.3.3]{Kas94}), the map $\varphi$ is surjective.

Finally we show that $\varphi$ becomes injective if we divide
\(
  \bigsqcup_\lambda \B_0(\lambda)\times \B(-\lambda)
\)
by $\aW$. Suppose
\(
   X_{1} X_{2} \cdots X_{N} b^* = X'_{1} X'_{2} \cdots X'_{N'} b^{\prime *}
\)
where each $X_j, X'_{j'} = \te_i$ or $\tf_i$, $b\in \B(-\lambda)$,
$b'\in\B(-\lambda')$. Then $b^*$ and $b^{\prime *}$ are in the same
connected component of $\B(\tU)$, which is isomorphic to $\B_0(\lambda)$
by $\bold{(C1)}$. By \thmref{thm:crystaldesc}(iii), we have
$b^{\prime *} = S_w b^*$ for some $w\in\aW$. In particular, $\lambda' =
w \lambda$.
We have
\(
   X_{1} X_{2} \cdots X_{N} b^* = 
   X'_{1} X'_{2} \cdots X'_{N'} S_w b^*.
\)
$\bold{(C1)}$ implies
\(
   X_{1} X_{2} \cdots X_{N} u_\lambda = 
   X'_{1} X'_{2} \cdots X'_{N'} S_w u_\lambda.
\)
The isomorphism $S_w^*\colon \B_0(\lambda)\to \B_0(w\lambda)$ sends
$u_\lambda$ to $S_w^{-1} u_{w\lambda} = S_w^{-1} u_{\lambda'}$.
Therefore
\begin{equation*}
  X'_{1} X'_{2} \cdots X'_{N'} u_{\lambda'}\otimes b'
  = (S_w^*\times S_w^*)(X_{1} X_{2} \cdots X_{N} u_\lambda\otimes b).
\end{equation*}
\end{proof}

\begin{Remark}
In \cite{Tos} a condition labeled $\bold{(C2)}$ requiring
$\B_0(\lambda)_\lambda = \{u_\lambda\}$ is used. This is false in
general. A counter-example for type $A_2^{(1)}$ can be found in
\cite[5.10]{Kas00}.
\end{Remark}

\begin{comment}
It seems to be false for all affine $\g$ of rank $> 1$.
\end{comment}

\section{\protect{The proof of \thmref{thm:intbase}(ii)}}\label{sec:sign}

In this section we will complete the proof of
\thmref{thm:intbase}. The proof uses results in the previous section,
in particular extremal weight modules.

\begin{rem}
Since we are working with generators in $\U^+$ we will need to consider the
extremal weight modules $V(-\lambda)$, $\lambda \in P_+^0$. 
\end{rem}

\begin{Lemma} \label{globalbaseofV}
Let $\lambda \in P_+^0$ and $\boc_0 \in \bocind$. Then 
\(
    S_{\boc_0} u_{-\lambda} \in G(\B(-\lambda)).
\)
\end{Lemma}

\begin{proof}
We have \( \sgn(\boc_0,0) S_{\boc_0} u_{-\lambda} \in G(\B(-\lambda))
\) by \propref{schurcor} and Theorem \ref{thm:crystaldesc}.  So we
only need to show $\sgn(\boc_0,0) = 1$. Therefore it is enough to show
$S_{\boc_0} u_{-\mu} \in G(\B(-\mu))$ for {\it some\/} $\mu$. We check
this by induction on $\sum_i |\rho^{(i)}|$. Take $\mu = \sum_i
\linebreak[0]\ell(\rho^{(i)}) \varpi_i$. Fix an $i$ and consider
$\rho^{\prime (i)}$ obtained by removing the last box in each row of
$\rho^{(i)}$. Let $\boc'_0$ be the set of partitions we get by
replacing $\rho^{(i)}$ by $\rho^{\prime (i)}.$

By \lemref{lem:Stal} and the Pieri rules for symmetric functions
\cite[p.~73]{Mac},
\begin{equation*}
    S_{\boc'_0} S_{\tal_i} u_{-\mu}
  = S_{\boc'_0} \tilde P_{i,\mu_i d_i} u_{-\mu}
  = S_{\boc_0} u_{-\mu},
\end{equation*}
where $\mu_i = \ell(\rho^{(i)})$. A priori, the right hand side is a
sum over all $S_{\boc_0} u_{-\mu}$ such that $\boc_0$ is an $n$--tuple
of partitions which are identical to $\boc'_0$ in factors different
than $i$, and for $i$ we have $\rho^{(i)} \setminus \rho^{\prime(i)}$
is a vertical $\mu_i$ strip. But there is only one such summand.
(Under the identification between Schur functions and irreducible
representations of $G_\mu$, $\tilde P_{i,\mu_i d_i}$ is identified
with a power of the determinant representation. So the tensor product
of it with an irreducible representation remains irreducible.)

Let us note the following fact:
If $\lambda' \in P^0$ and $\lambda = \sum_i \lambda_i \varpi_i \in
P_+^0$ satisfy $\cl(\lambda') = \cl(\lambda)$ then
$u_\infty \otimes t_{-\lambda} \otimes b_2 \in \B(\lambda) \iff
u_\infty \otimes t_{-\lambda'} \otimes b_2 \in \B(\lambda')$ (see
\cite[Appendix]{Kas00}).  In other words, $S_{\boc_0}
u_{-\lambda} \bmod q_s \La(\lambda) \in \B(\lambda) \iff
S_{\boc_0} u_{-\lambda'} \bmod q_s \La(\lambda') \in
\B(\lambda')$.
Then by the inductive assumption, $S_{\boc'_0} u_{\mu_i d_i
\delta-\mu}$ is an element of the global crystal basis. Therefore
$S_{\boc_0} u_{-\mu} = S_{\boc'_0} S_{\tal_i} u_{-\mu}$ is also in the
global basis, since the isomorphism $V(\mu_i d_i \delta - \mu)
\xrightarrow{\sim} V(-\mu)$ sends the global basis to the global basis
(see \cite[Prop. 8.2.2 (iv)]{Kas94}).
\end{proof}

\begin{proof}[Proof of $\sgn(\boc,p) = 1$ in Proposition \ref{crystalbase}]
Let $\boc_0 \in \aC$ be purely imaginary. In the proof of
\lemref{globalbaseofV} we have shown $\sgn(\boc_0,0) = 1.$ From 
\cite[Proposition~8.3(a)]{Lu-bg} we have $\sgn(\boc,0) = 1$ for an
arbitrary $\boc$. The general case now follows from
\cite[Theorem~1.2]{Lu-bg} (see \subsecref{subsec:Lu-bg}).
\end{proof}

Let $b(\boc,p)$ as in \propref{crystalbase}. We denote
$G(b(\boc,p))$ by $G(\boc,p)$.
By the same argument in the proof of \cite[42.1.12]{Lu-Book}, we have
\begin{equation}\label{eq:ve}
\begin{gathered}
  \varphi_{i_p}(b(\boc,p)) = \boc(p)
  \quad\text{(hence $G(\boc,p)\in E_{i_p}^{\boc(p)}\U^+$)},
\\
  \te_{i_p}(b(\boc,p)) = b(\boc^+,p),
\\
  \tf_{i_p}(b(\boc,p)) = b(\boc^-,p) \quad\text{if $ \boc(p) > 0$},
\end{gathered}
\end{equation}
where $\boc^+$ (resp.\ $\boc^-$) is defined by setting $\boc^+(p) =
\boc(p)+1$ (resp.\ $\boc^-(p) = \boc(p)-1$) and other entries are the
same as $\boc$. We have similar formulas for $\varphi_{i_{p+1}}^*$,
$\tf_{i_{p+1}}^*$, $\te_{i_{p+1}}^*$.

Let $\Up[i]$, ${}^*\Up[i]$, ${}^i\pi$, $\pi^i$ be as in
\subsecref{subsec:Lu-bg}.
By \eqref{eq:ve} we have
\begin{equation}\label{eq:nonzero_trans}
   {}^{i_p}\pi(G(\boc,p)) \neq 0 \Longleftrightarrow
   \boc(p) = 0  \Longleftrightarrow
   \pi^{i_{p}}(G(\boc,p-1)) \neq 0,
\end{equation}
and
\begin{equation}\label{eq:flip}
    T_{i_p}({}^{i_p}\pi(G(\boc,p))) = \pi^{i_{p}}(G(\boc,p-1)).
\end{equation}
The latter follows from \eqref{eq:trans} (see \cite[\S6.2]{Lu-bg} for
detail).

Let us define $a_{\boc,\boc'}^p \in \Q(q_s)$ by
\begin{equation*}
  L(\boc,p) = \sum_{\boc'} a_{\boc,\boc'}^p G(\boc',p).
\end{equation*}
The sum is over a finite set of $\boc' \in \aC$ indexing all elements
$G(\boc,p) \in \Up$ having the same weight with $L(\boc,p)$. Since the
global basis is an integral base of $\Ua^+$ and $\La(-\infty)$ we have
$a_{\boc,\boc'}\in \ca \cap \A_\infty = \Z[q_s^{-1}]$. We also have
$a_{\boc,\boc'}|_{q=\infty} = \delta_{\boc,\boc'}$ by
\propref{crystalbase} together with $\sgn(\boc,p) = 1$.
\thmref{thm:intbase}(ii) follows from
\begin{Lemma}\label{lem:upper}
  The transition matrix $(a_{\boc,\boc'})$ between the global basis
and the $L(\boc,p)$, $\boc \in \aC, p \in \Z$ is upper triangular with
respect to the ordering $\prec_p$, and the diagonal entries are $1$.
\end{Lemma}

\begin{pf}
We first consider the case $\boc_+ = 0 = \boc_-$, so that
\begin{equation*}
    L(\boc,p) = T_{i_{p+1}}T_{i_{p+2}} \cdots T_{i_0} (S_{\boc_0})
\quad\text{or}\quad T_{i_p}^{-1} T_{i_{p-1}}^{-1} \cdots T_{i_2}^{-1}
T_{i_1}^{-1}(S_{\boc_0}).
\end{equation*}
In this case, by Lemma \ref{globalbaseofV} we know that
$L(\boc,p)u_{-\lambda}$ is in the global basis of $V(-\lambda)$ for any
sufficiently large $\lambda \in P_+^0$. Therefore $(L(\boc,p) -
G(\boc,p))u_{-\lambda} = 0$ for any sufficiently large $\lambda \in
P_+^0$, i.e.,
\begin{equation*}
    \sum_{\boc'\neq\boc} (a_{\boc,\boc'}^p - \delta_{\boc,\boc'})
    G(\boc',p) u_{-\lambda} = 0.
\end{equation*}
By the construction of $V(-\lambda)$, $\{ G(\boc',p)u_{-\lambda} \mid
\boc'\in \aC\}$ is mapped to the union of the global basis of
$V(-\lambda)$ and $0$. Hence
\begin{equation*}
    a_{\boc,\boc'}^p = \delta_{\boc,\boc'} \quad\text{or}\quad
    G(\boc',p)u_{-\lambda} = 0.
\end{equation*}
If $\boc'_+ = 0 = \boc'_-$, then
\begin{equation*}
  G(\boc',p)u_{-\lambda} =
  T_{i_{p+1}}T_{i_{p+2}} \cdots T_{i_0} (S_{\boc'_0}) u_{-\lambda}
\quad\text{or}\quad
    T_{i_p}^{-1} T_{i_{p-1}}^{-1} \dots T_{i_2}^{-1} T_{i_1}^{-1}(S_{\boc'_0})
    u_{-\lambda}
\end{equation*}
is nonzero for sufficiently large $\lambda$. This means that
$a_{\boc,\boc'}^p = \delta_{\boc,\boc'}$ for such $\boc'$. By the
definition of the ordering, the remaining $\boc'$'s have $\boc'
\succ_p \boc$. We have the assertion in this case.

Next consider the case $\boc({p+1}) = \boc({p+2}) = \cdots
= 0$, i.e.,
\begin{equation*}
  \begin{split}
    L(\boc,p) = & \left(E_{i_{p}}^{(\boc(p))}
    T_{i_{p}}^{-1}(E_{i_{p-1}}^{(\boc({p-1}))})
    T_{i_{p}}^{-1}T_{i_{p-1}}^{-1}(E_{i_{p-2}}^{(\boc({p-2}))})
    \cdots\right)
\\
    & \times T_{i_{p+1}}T_{i_{p+2}} \cdots T_{i_0} (S_{\boc_0}).
  \end{split}
\end{equation*}
We prove the assertion by the induction on $q$ such that
$\boc({p-q}) = \boc({p-q-1}) = \cdots = 0$. When $q=0$, we have
$\boc_+ = 0 = \boc_-$, which we have checked already.

First assume $\boc(p) = 0$. We consider
\(
   L(\boc,p-1) = T_{i_p}L(\boc,p)
\)
(see \eqref{eq:trans}). By the induction hypothesis, we have
\begin{equation*}
   L(\boc,p-1) = 
   G(\boc,p-1) + 
   \sum_{\boc \prec_{p-1}\boc'}a_{\boc,\boc'}^{p-1} G(\boc',p-1).
\end{equation*}
We apply the composition $T_{i_p}^{-1} \circ \pi^{i_p}$ to both
sides. We have $L(\boc,p-1)\in {}^* \U^+[i_p]$, so the left hand side
becomes $L(\boc,p)$. Therefore
\begin{equation*}
   L(\boc,p) = T_{i_p}^{-1}(\pi^{i_p} (G(\boc,p-1))) + 
   \sum_{\boc \prec_{p-1}\boc'}a_{\boc,\boc'}^{p-1}
   T_{i_p}^{-1}(\pi^{i_p} (G(\boc',p-1))).
\end{equation*}
By \eqref{eq:flip} the right hand side is contained in
\begin{equation*}
   G(\boc,p) 
   + \sum_{\substack{\boc \prec_{p-1}\boc' \\ \boc'(p) = 0}}
     a_{\boc,\boc'}^{p-1}
   G(\boc',p) + E_{i_p} \U^+.
\end{equation*}
The condition $\boc'(p) = 0$ comes from ${}^{i_p}\pi
(G(\boc',p-1))\neq 0$ (see \eqref{eq:nonzero_trans}).
The part in $E_{i_p} \U^+$ is a linear combination of $G(\boc'',p)$'s
with $\boc''(p) > 0$. Each such $\boc''$ is greater than $\boc$ with
respect to $\prec_p$. The summation in the second term can be replaced
as $\sum_{\boc \prec_p\boc', \boc'(p) = 0}$ by \eqref{eq:ord_trans}.
Thus we have the assertion under the assumption $\boc(p) = 0$.

Next we assume $\boc(p) > 0$. Let us define $\tilde\boc$
by setting $\tilde\boc(p) = 0$ and all other entries are the
same as $\boc$. We have
\begin{equation*}
    L(\boc,p) = E_{i_p}^{(\boc(p))}L(\tilde\boc,p).
\end{equation*}
Since $\tilde\boc(p) = 0$, we have just proved
\begin{equation} \label{eq:indstep}
    L(\tilde\boc,p) = G(\tilde\boc,p) + 
    \sum_{\tilde\boc' \succ_p \tilde\boc} 
    a_{\tilde\boc,\tilde\boc'}^p G(\tilde\boc',p).
\end{equation}
By \cite[14.3]{Lu-Book},
\(
  E_{i_p}^{(\boc(p))}G(\tilde\boc,p) = G(\boc,p)
\)
plus an element in $E_{i_p}^{\boc(p)+1}\U^+$. If we write the
element in a linear combination of $G(\boc'',p)$'s, all elements appearing
satisfy $\boc'' \succ_p \boc$ by \eqref{eq:ve}.
Next consider \( E_{i_p}^{(\boc(p))}G(\tilde\boc',p), \)
obtained from \eqref{eq:indstep} by multiplication by
$E_{i_p}^{(\boc(p))}$.  If $\tilde\boc'(p) > 0$, then this
is an element in $E_{i_p}^{\boc(p)+1}\U^+$. Otherwise, it is
equal to \( G(\boc',p) \) plus an element in
$E_{i_p}^{\boc(p)+1}\U^+$, where $\boc'$ is defined by setting
$\boc'(p) = \boc(p)$ and all other entries are the same as
$\tilde\boc'$. In either cases, it is a linear combination of elements
$G(\boc'',p)$ with $\boc''\succ_p\boc$.  Thus we have the assertion
for $G(\boc,p)$.

Finally we consider the general case. We first remark
\begin{equation}\label{eq:lexineq}
  \left( \boc(p), \boc({p-1}), \boc({p-2}),
   \dots\right)
  \le
   \left( \boc'(p), \boc'({p-1}), \boc'({p-2}),
   \dots\right).
\end{equation}
This is proved by the induction on $q$ such that $\boc({p-q}) =
\boc({p-q-1}) = \cdots = 0$. When $q=0$, the left hand side is
$(0,0,\dots)$, so the inequality trivially holds. The remaining
argument of the induction is exactly the same as above.

Now we can prove the assertion of the lemma by the induction on $q$
such that $\boc({p+q}) = \boc({p+q+1}) = \cdots = 0$. When
$q=1$, we are reduced to the case studied above, and have the
assertion.
The remaining argument is almost the same as above. When
$\boc({p+1}) = 0$, we apply the induction hypothesis to
$L(\boc,p+1)$ and consider $T_{i_{p+1}}{}^{i_{p+1}}\pi L(\boc,p+1)$. We
use \eqref{eq:flip} to get
\begin{equation*}
   L(\boc,p) \in G(\boc,p) 
   + \sum_{\substack{\boc \prec_{p+1}\boc' \\ \boc'(\beta_{p+1}) = 0}}
     a_{\boc,\boc'}^{p+1} G(\boc',p)
   + \U^+E_{i_{p+1}} .
\end{equation*}
The summation in the second part can be replaced as
$\sum_{\boc \prec_p\boc', \boc'({p+1}) = 0}$ by \eqref{eq:ord_trans}.
The part in $\U^+E_{i_{p+1}} $ is a linear combination of $G(\boc'',p)$'s
with $\boc''({p+1}) > 0$. Since we already have
\eqref{eq:lexineq} (with $\boc'$ replaced by $\boc''$), we have
$\boc \prec_p \boc''$. We have the assertion in the case
$\boc({p+1}) = 0$.
The case $\boc({p+1}) > 0$ can be reduced to the case $\boc({p+1}) =
0$ as above with use of \eqref{eq:lexineq}.
\end{pf}

\section{Cell structure of $\Um$.}

In this section we prove Lusztig's conjecture \cite{Lu-v} on
two--sided cells of the modified quantum affine algebra of level
$0$. Our strategy follows the same line as his proof of the
corresponding conjecture for finite type cases. However, here we need
to show that a certain bi-module $\Um[\lambda]$ has a global crystal
basis, whereas for finite type cases this assertion is a direct
consequence of the definition. In defining $\Uml$, our earlier results
on extremal weight modules play a key role. Their role is analogous to
the role that dominant highest weight representations play in the
finite type case.

Let $\Um$ be the level $0$ modified quantum affine algebra, i.e., \(
\Um = \bigoplus_{\lambda\in P^0_\cl} \U a_\lambda.  \) Unfortunately
this contradicts the notation in \secref{sec:crystal}, where $\Um$ was
$\bigoplus_{\lambda\in P} \U a_\lambda$, but we do not want to
introduce new notation. Let $\Lm$ be its crystal lattice, $\Bm$ be its
crystal base.

Throughout this section, $\lambda$ is a dominant classical level $0$ 
weight, i.e., $\lambda \in \cP = \sum_{i\in I_0} \Z_{\ge 0}
\cl(\varpi_i)$, except in the proof of \thmref{thm:Umcry}.
Let $V(\lambda)$ be the extremal weight module of weight
$\lambda$. A priori, it is defined for $\lambda\in P^0$, but its
$\Um$-module structure depends only on $\cl(\lambda)$. 

\subsection{A bi-module $\Uml$}

\begin{Definition}
For
$b\in\B(\lambda)$, we denote by $G_\lambda(b)$ the corresponding
element in the global basis of $V(\lambda)$. If we consider $b$ as an
element of $\Bm$, we denote by $G(b)$ the corresponding element in
$\Um$.
\end{Definition}

Denote $\# = * \circ \vee$. It is an involutive anti-automorphism of
$\Um$. By \cite[4.3.2]{Kas94} it leaves the global crystal basis of
$\Um$ invariant. (The result was proved for $*$, but the same proof
works for $\vee$.)
We have 

\begin{Lemma}[\protect{\cite[19.1.1]{Lu-Book}, \cite[3.7]{Lu-v}}]
\label{lem:psi}
Let $x\in a_\lambda \Um a_{\lambda'}$ and $\nu =
\lambda - \lambda'$. Then
\begin{equation}\label{eq:psi}
   \psi(x) = q^{-(\nu, \lambda+\lambda')/2} x^\#.
\end{equation}
\end{Lemma}

We denote by $\le$ the dominance partial order on the level $0$
classical weights relative to the fundamental level zero weights
$\cl(\varpi_i) \in \cP$ defined in $\S 2.6$.

Let us give a slightly different parametrization of $\B(\lambda)$ from
that of \secref{sec:crystal}. Let $W(\lambda)$, $\B_W(\lambda)$ be as
in the proof of \lemref{lem:Wlam}. We have maps
\begin{equation*}
    \B_0(\lambda) \xrightarrow{\Phi_{\lam | q=0}}
    \{ \text{monomials in $z_{i,\nu}^\pm$} \} \times \B_W(\lambda)
    \xrightarrow{\text{projection}} \B_W(\lambda).
\end{equation*}
The composition is surjective since $\B_W(\lambda)$ is connected by
\cite[Lemmas~1.9, 1.10]{AK}. By \lemref{lem:Wlam} each fiber is
identified with the set of monomials $p(z) = \prod_i (z_{i,1}\cdots
\linebreak[0]z_{i,\lambda_i})^{r_i}$ ($r_i\in\Z$). We choose and fix a
section $\B_W(\lambda)\to \B_0(\lambda)$. (We do not require that it
respect the crystal structure.) Then we have an identification (of
sets, not of crystals)
\begin{equation}\label{eq:ident}
   \B(\lambda) \cong \Irr G_\lambda \times \B_W(\lambda),
\end{equation}
where $G_\lambda = \prod_i \GL_{\lambda_i}(\C)$, and $\Irr G_\lambda$
is the set of irreducible representations of $G_\lambda$. We identify
$\Irr G_\lambda$ with the set of Schur Laurent polynomials in
$\{z_{i,\nu}^\pm \}$. We denote by $s(z)$ the polynomial corresponding
to $s$.

\begin{Definition}\label{def:pureim}
For $s\in\Irr G_\lambda$ we define $S\in G(\B(\lambda))$ by
$\Phi_\lambda(S u_\lambda) = s(z) \tilde u_\lambda$. The existence and
uniqueness are guaranteed in \thmref{thm:crystaldesc}(ii). If $s$ is
the dual of a polynomial representation and corresponds to $\boc_0$,
we have $S = S^-_{\boc_0} u_\lambda$ by \propref{schurcor}. We also
regard $S$ as an element in $G(\B(\U a_\lambda))$.
\end{Definition}

We will use the correspondence between $s\leftrightarrow S$ hereafter.

\begin{Definition}
Let $\Um[{}^\ge\lambda]$ (resp.\ $\Um[{}^>\lambda]$) be the two--sided
ideal of $\Um$ consisting of all elements $x\in\Um$ acting on
$V(\lambda')$ by $0$ for any $\lambda'\ngeq\lambda$ (resp.\
$\lambda'\ngtr\lambda$).  Let
$\Um[\lambda] = \Um[{}^\ge\lambda]/\Um[{}^>\lambda]$. We define
$\Umz[{}^\ge\lambda]$, $\Umz[{}^>\lambda]$ and $\Umz[\lambda]$ in an
analogous manner.

For $\xi$, $\xi'\in P_+$ let $\Um[{}^\ge\lambda]_{\xi,\xi'}$ (resp.\ 
$\Um[{}^>\lambda]_{\xi,\xi'}$) be the image of $\Um[{}^\ge\lambda]$
(resp.\ $\Um[{}^>\lambda]$) under the natural map $\Um\to
V(\xi)\otimes V(-\xi')$ given by $x\mapsto x(u_\xi\otimes u_{-\xi'})$.
Let $\Um[\lambda]_{\xi,\xi'} =
\Um[{}^\ge\lambda]_{\xi,\xi'}/\Um[{}^>\lambda]_{\xi,\xi'}$.
\end{Definition}

\begin{Lemma}\label{lem:inv}
\textup{(i)} $\Um[{}^\ge\lambda]$ is invariant under
$\setbox5=\hbox{A}\overline{\rule{0mm}{\ht5}\hspace*{\wd5}}\,$.

\textup{(ii)} For any $\xi$, $\xi'\in P_+$, the $\U$-submodules
$\Um[{}^\ge\lambda]_{\xi,\xi'}$ and $\Um[{}^>\lambda]_{\xi,\xi'}$ of
$V(\xi)\otimes V(-\xi')$ are invariant under $\widetilde F_i^{(n)}$
defined in \subsecref{subsec:crystal}.
\end{Lemma}

\begin{proof}
(i) The bar involution
$\setbox5=\hbox{A}\overline{\rule{0mm}{\ht5}\hspace*{\wd5}}\,$ on
$V(\lambda)$ satisfies $\overline{xu} = \overline{x}\, \overline{u}$
for $x\in\Um$, $u\in V(\lambda)$. The assertion follows.

(ii) Let $x\in\Um[{}^\ge\lambda]$. We have
\[
  \widetilde F_i^{(n)} x (u_\xi\otimes u_{-\xi'}) 
  = \sum_k F_i^{(n+k)} E_i^{(k)} a_k^n(t_i) x(u_\xi\otimes u_{-\xi'}),
\]
where all but finitely many terms of the sum are $0$. Since each
$F_i^{(n+k)} E_i^{(k)} a_k^n(t_i) x$ is contained in
$\Um[{}^\ge\lambda]$, we have the assertion.
\end{proof}

\begin{comment}
In earlier version, we stated

\textup{(ii)} $\Um[{}^\ge\lambda]$ is invariant under $\te_i$, $\tf_i$.

\begin{proof}
(ii) The assertion is a consequence of $(\te_i x)u = \te_i(xu)$,
$(\tf_i x)u = \tf_i(xu)$ for $x\in\Um$, $u\in V(\lambda)$.
\end{proof}

But this is not a correct statement since $\te_i$, $\tf_i$ are not
{\it well-defined\/} on $\Um$.

Another way to correct is :

Let $(\operatorname{forget})$ be the functor from the category of
integrable $\U$-modules to the category of $\Q(q)$-vector spaces,
forgetting the $\U$-module structure. Let $R$ be its endomorphism
ring. Note that $\U$ is a subalgebra of $R$. The element $\widetilde
F_i^{(n)}$ in \subsecref{subsec:crystal} belongs to $R$. We define
$R[{}^\ge\lambda]$, $R[{}^>\lambda]$ in the same way as above.

\textup{(ii)} $R[{}^\ge\lambda]$ and $R[{}^>\lambda]$ are invariant
under $\widetilde F_i^{(n)}$.

\begin{proof}
The assertion is a consequence of $(\te_i x)u = \te_i(xu)$,
$(\tf_i x)u = \tf_i(xu)$ for $x\in R$, $u\in V(\lambda)$.
\end{proof}
\end{comment}

\begin{Lemma}[\protect{cf.\ \cite[4.5]{Lu-v}}]\label{lem:a}
\textup{(i)} $\U a_\lambda, a_\lambda \U\subset\Um[{}^\ge\lambda]$.

\textup{(ii)} For $b\in\B(\U a_\lambda)$, $G(b) \in \Um[{}^> \lambda]$
if and only if $b\notin \B(\lambda)$.
\end{Lemma}

\begin{proof}
(i) Let $x a_\lambda\in \U a_\lambda$. If $x a_\lambda$ acts on
$V(\lambda')$ by a nonzero map, then so does $a_\lambda$. Since
$a_\lambda$ is a projector to the weight space with weight $\lambda$,
$\lambda$ is a weight of $V(\lambda')$. Therefore we have
$\lambda\le\lambda'$. Similarly, we have
$a_\lambda\U\subset\Um[{}^\ge\lambda]$.

(ii) Let us first remark that $V(\lambda)_\lambda$ has a basis $\{ S
u_\lambda \mid s\in\Irr G_\lambda \}$. This follows from results in
\secref{sec:crystal}. Note that $\lambda\in P^0_\cl$ in this section,
and the weight space is a direct sum of weight spaces $V(\lambda)_\mu$
with $\cl(\mu) = \lambda$.

Let $b\in \B(\U a_\lambda)$. Then $G(b) = 0$ as an operator on
$V(\lambda)$ if and only if $G(b) S u_\lambda = 0$ for any $s \in\Irr
G_\lambda$. This is equivalent to $\Phi_\lambda(G(b) S u_\lambda) =
s(z)\Phi_\lambda(G(b)u_\lambda) = 0$, and in particular $G(b)u_\lambda
= 0$. If $b\in\B(\lambda)$, then $G(b)u_\lambda$ is nothing but
$G_\lambda(b)$.  In particular, it is nonzero. If
$b\notin\B(\lambda)$, then $G(b)u_\lambda = 0$ by the definition of
$V(\lambda)$.
\end{proof}

Let $(\ ,\ )_\lambda$ be the bilinear form on $V(\lambda)$ satisfying
$(xu, v)_\lambda = (u, \psi(x)v)_\lambda$ for $u,v\in V(\lambda)$,
$x\in \Um$ and $(u_\lambda, G_\lambda(b))_\lambda = 1$ or $0$
according to $G_\lambda(b) = u_\lambda$ or not (see
\propref{prop:bilinear}).

\begin{Proposition} \label{prop:sharpinv}
$\Um[{}^\ge\lambda]$ is invariant under $\#$.
\end{Proposition}

\begin{proof}
Since $\Um[{}^\ge\lambda]$ is a direct sum of its intersection with
$a_{\lambda'} \U a_{\lambda''}$, by \eqref{eq:psi} it is enough to show
that $x\in\Um[{}^\ge\lambda]$ if and only if
$\psi(x)\in\Um[{}^\ge\lambda]$.
The global crystal basis is almost orthonormal, that is \(
(G_\lambda(b), G_\lambda(b'))_\lambda \in \delta_{bb'} + q_s\A_0.  \)
(see \thmref{thm:main2}). In
particular, it is non-degenerate on $V(\lambda)$.  Therefore, $x \in
\operatorname{Ann} V(\lambda) \iff (xG_\lambda(b), G_\lambda(b'))$ for
all $b, b' \in \B(\lambda) \iff (G_\lambda(b), \psi(x) G_\lambda(b'))
= 0 \iff \psi(x) \in \operatorname{Ann} V(\lambda).$ It follows
$x\in\Um[{}^\ge\lambda]$ if and only if
$\psi(x)\in\Um[{}^\ge\lambda]$.
\end{proof}

\begin{rem} Since \propref{prop:sharpinv} implies $x \in \operatorname{Ann}
V(\lambda)$ if and only if $x^\# \in \operatorname{Ann} V(\lambda)$,
 \lemref{lem:a} implies that $b \in \B(\lambda)$ if and only if 
 $G(b)^\#\in\Um[{}^\ge\lambda]\setminus\Um[{}^>\lambda]$ for
 $b\in\B(\U a_\lambda)$.
\end{rem}

Note that $\Um[\lambda] = \Umgel/\Umgl$ is integrable as a
$\Um$-module by \thmref{thm:convex}. In particular, the operators
$T_w$ are defined.

\begin{Lemma}\label{lem:extrem}
Let $b\in\B(\lambda)$. Then $G(b)^\#\bmod\Um[{}^>\lambda]$ is an
extremal vector in the $\Um$-module $\Um[\lambda]$.
\end{Lemma}

\begin{proof}
We show that $S_{i_1}\cdots S_{i_N} G(b)^\#\bmod\Um[{}^>\lambda]$ is
$i$-extremal for all $i\in I$ by induction on $N$.

If $N=0$, we need to show that $e_i G(b)^\#$ or $f_i G(b)^\#$ acts on
$V(\lambda)$ by $0$ for each $i\in I$. In other words if and only if
$e_i G(b)^\#$ or $f_i G(b)^\# \in \Umgl$.  Note that $G(b)^\# =
a_\lambda G(b)^\#$ and $a_\lambda$ is a projector to the weight space
of weight $\lambda$. Since $V(\lambda)_\lambda$ has a basis $\{ S
u_\lambda \mid s\in\Irr G_\lambda \}$ consisting of extremal vectors,
we are done.
Moreover, we have \( S_i G(b)^\#\bmod\Um[{}^>\lambda] = a_{s_i\lambda}
S_i G(b)^\#\bmod\Um[{}^>\lambda] \).

We show the statement for $N$ assuming the statement for $N-1$.
By the induction hypothesis and the last part of the above argument,
$S_{i_2}\cdots S_{i_N} G(b)^\#\bmod\Um[{}^>\lambda]$ is $i$-extremal 
for $i \in I$ and
\[
  S_{i_1}\cdots S_{i_N} G(b)^\#\bmod\Um[{}^>\lambda]
   = a_{s_{i_1}\cdots s_{i_N}\lambda}  
S_{i_1} \cdots S_{i_N}  G(b)^\#\bmod\Um[{}^>\lambda].
\]
Since $V(\lambda)_{s_{i_1}\cdots s_{i_N}\lambda}$ is spanned by
extremal vectors, we are done.
\end{proof}

By \lemref{lem:extrem}, for each $b\in \B(\lambda)$ we have a
unique $\Um$-homomorphism
\begin{equation*}
   \Psi_{b} \colon V(\lambda) \to \Um[\lambda],
\end{equation*}
which sends $u_\lambda$ to $G(b)^\#\bmod\Um[{}^>\lambda]$.
By definition, we have
\[
   \Psi_{b}(G_\lambda(b_1)) = G(b_1)G(b)^\#\bmod\Um[{}^>\lambda].
\]

\begin{Lemma}\label{lem:mul1}
Let $b_1, b_2\in\B(\lambda)$. Then
\begin{equation*}
   G(b_1)^\# G(b_2) \equiv q^n \sum_{s\in\Irr G_\lambda}
   (G_\lambda(b_2), G(b_1)S u_\lambda)_\lambda
   S a_\lambda \mod \Um[{}^>\lambda],
\end{equation*}
where $n = (\wt b_1, 2\lambda + \wt b_1)/2$ and $S$ is an element
corresponding to $s$ by Definition~\ref{def:pureim}.
\end{Lemma}

\begin{proof}
Since $\Um[{}^\ge\lambda]$ is a two-sided ideal of $\Um$, both sides
are in $\Um[{}^\ge\lambda]$. Therefore it is enough to show that both
sides define the same operator on $V(\lambda)$.  Since we have
$G(b_1)^\# G(b_2)u_\lambda = a_\lambda G(b_1)^\# G(b_2)u_\lambda$, it
is contained in the weight space $V(\lambda)_\lambda$. We have a basis
$\{ S u_\lambda \}$ which is orthonormal by
\cite[Proposition~4.10]{extrem}. Hence we have
\begin{equation*}
   G(b_1)^\# G(b_2) u_\lambda
   = \sum_{s\in\Irr G_\lambda}
   (G(b_1)^\# G(b_2)u_\lambda, S u_\lambda)_\lambda 
        S u_\lambda.
\end{equation*}
Now the assertion follows from \eqref{eq:psi}.
\end{proof}

\begin{Lemma}\label{lem:dual}
Let $S$ correspond to $s\in \Irr G_\lambda$ as in
Definition~\ref{def:pureim}. We consider it as an element in $G(\B(\U
a_\lambda))$. Then $S^\#$ corresponds to the dual representation of
$s$, which will be denoted by $s^\#$ hereafter.
\end{Lemma}

\begin{proof}
By \lemref{lem:psi} we have $(S^\# u_\lambda, S_1 u_\lambda) =
(u_\lambda, S S_1 u_\lambda)$ for any $s_1\in\Irr G_\lambda$.
By \cite[4.10]{extrem} we have $(u_\lambda, S S_1 u_\lambda) =
(\Phi_\lambda(S u_\lambda), \Phi_\lambda(S S_1 u_\lambda)) = (\utl,
s(z) s_1(z)\utl)^\sim$, where $(\ ,\ )^\sim$ is a bilinear form on
$\Vtl$ defined in [loc.cit.\ \S4].
By a property of $(\ ,\ )^\sim$, we have
\[
  (\utl, s(z)s_1(z)\utl)^\sim = (s(z^{-1}) \utl, s_1(z) \utl)^\sim
  = 
  \begin{cases}
     1 & \text{if $s(z^{-1}) = s_1(z)$},
     \\
     0 & \text{otherwise}.
  \end{cases}
\]
(See an equation in the middle of [loc.cit., Proof of 4.9] and
[loc.cit., Proof of 4.10]). Since $\{ S u_\lambda \}$ is an orthonormal
basis of $V(\lambda)_\lambda$, this means $S^\# u_\lambda$ corresponds
to the Schur Laurent polynomial $s(z^{-1})$, which corresponds to the
dual representation of $s$ in turn.
\end{proof}

\begin{Lemma}\label{lem:S}
Let $(s,b)\in \Irr G_\lambda\times\B_W(\lambda)$ and let $\tilde b$ be
the corresponding element in $B(\lambda)$ under the identification
\eqref{eq:ident}. We have
\begin{equation*}
   G(b) S \equiv G(\tilde b) \mod\Umgl, \quad
   S^\# G(b)^\# \equiv G(\tilde b)^\# \mod\Um[{}^>\lambda].
\end{equation*}
\end{Lemma}

\begin{proof}
By the construction of \eqref{eq:ident}, $\Phi_\lambda(G(\tilde
b)u_\lambda) = s(z)\Phi_\lambda(G_\lambda(b)) = s(z)G(b)\utl$. By
Definition~\ref{def:pureim}, this is equal to
$G(b)\Phi_\lambda(Su_\lambda) =
\Phi_\lambda(G(b)Su_\lambda)$. Therefore we have $G(b) S u_\lambda =
G(\tilde b) u_\lambda$.

Let us consider the first equation of the assertion. Since both sides
of the equation are in $\U a_\lambda$, the result follows from
\lemref{lem:a}(ii) if we prove that they define the same operators on
$V(\lambda)$. Since $V(\lambda)_\lambda$ is spanned by $\{ S_1
u_\lambda \mid s_1\in\operatorname{Irr}G_\lambda\}$, it is enough to
show $G(b) S S_1 u_\lambda = G(\tilde b)S_1 u_\lambda$ for every
$S_1$.  Consider the embedding $\Phi_\lambda \colon V(\lambda)
\rightarrow \Wtl$ as in \corref{cor:inj}. By \propref{schurcor} and
the fact that $s_{\boc_0}(z^{-1})$ is $\Um$-linear, we have
\begin{equation*}
   \Phi_\lambda(G(b) S S_1 u_\lambda) 
   = s_1(z^{-1}) \Phi_\lambda(G(b)S u_\lambda)
   = s_1(z^{-1}) \Phi_\lambda(G(\tilde b) u_\lambda)
   = \Phi_\lambda(G(\tilde b)S' u_\lambda).
\end{equation*}
We get the assertion by the injectivity of $\Phi_\lambda$.

The second equation follows from the first together with
\propref{prop:sharpinv}.
\end{proof}

\begin{Lemma}\label{lem:Schurmul}
Let $s_1, s_2\in \Irr G_\lambda$. Let $S_1$, $S_2\in G(\B(\U
a_\lambda))$ be the corresponding elements by
Definition~\ref{def:pureim}. We have
\begin{equation*}
   S_1 S_2 \equiv \sum_{s} c_{s_1s_2}^s S \mod \Um[{}^>\lambda],
\end{equation*}
where $s_1 s_2 = \sum_s c_{s_1s_2}^s s$.
\end{Lemma}

\begin{proof}
Since both sides are in $\Umgel$, it is enough to show that they
define the same operator on $V(\lambda)$. As in the proof of
\lemref{lem:S}, it follows from $S_1 S_2 u_\lambda = 
\sum_{s} c_{s_1s_2}^s S u_\lambda$. But this follows directly from
Definition~\ref{def:pureim}.
\end{proof}

\begin{Definition}
For $b, b'\in\B_W(\lambda)$, $s\in\operatorname{Irr}G_\lambda$ we set
\begin{equation*}
   G_\lambda(b,s,b') = G(b) S G(b')^\# \bmod \Um[{}^>\lambda]
   = \Psi_{b'}(G(b)S u_\lambda) \in \Um[\lambda].
\end{equation*}
\end{Definition}

The following is a direct consequence of \lemref{lem:mul1}.
\begin{Lemma}\label{lem:mul2}
We have
\begin{equation*}
\begin{split}
  & G_\lambda(b_1,s_1,b'_1) G_\lambda(b_2,s_2,b'_2)
\\
 =\; &
   q^n \sum_{s''\in\operatorname{Irr}G_\lambda}
   (G(b_2)S_2 u_\lambda, G(b'_1)S'' u_\lambda)_\lambda
   G(b_1) S_1 S'' G(b'_2)^\#
   \bmod \Um[{}^>\lambda],
\end{split}
\end{equation*}
where $n = (\wt b'_1, 2\lambda + \wt b'_1)/2$.
\end{Lemma}

\begin{Lemma}\label{lem:linind}
\(
   \left.\left\{ G_\lambda(b,s,b')\, \right|
     b, b' \in \B_W(\lambda), s\in\operatorname{Irr}G_\lambda \right\}
\)
is linearly independent.
\end{Lemma}

\begin{proof}
Suppose that
\(
   \sum_{b,b',s} a_{b,b',s} G(b) S G(b')^\#
\)
acts on $V(\lambda)$ by $0$. Then for any $b_1\in\B(\lambda)$
{\allowdisplaybreaks
\begin{equation*}
\begin{split}
  0  & = \sum_{b,b',s} a_{b,b',s} G(b) S G(b')^\# G_\lambda(b_1)
\\
     & = \sum_{b,b',s} a_{b,b',s} G(b) \left(G(b')S^\#\right)^\#
     G_\lambda(b_1)
\\
     & = \sum_{b,b',s,s'} a_{b,b',s} q^{n(b')}
     (G_\lambda(b_1), G(b')S^\# S' u_\lambda)_\lambda
     G(b) S'u_\lambda,
\end{split}
\end{equation*}
where $n(b') = (\wt b', 2\lambda + \wt b')/2$.} Here we have
used \lemref{lem:mul1} in the third equality. Since
\(
   \{ G(b) S' u_\lambda \mid b\in\B_W(\lambda), s'\in\Irr G_\lambda \}
\)
is linearly independent, we have
\begin{equation*}
   \sum_{b', s} a_{b,b',s} q^{n(b')}
   (G_\lambda(b_1), G(b')S^\# S' u_\lambda)_\lambda
     = 0
\end{equation*}
for any $b$, $b_1$, $s'$. This equality for $s'=1$, together with the
nondegeneracy of $(\ ,\ )_\lambda$ and linearly independence of $\{
G(b')S^\#u_\lambda \}$ imply $a_{b, b',s} = 0$.
\end{proof}

\begin{Definition}
For $\xi,\xi'\in P_+$, let $\Lmgel_{\xi,\xi'} =
(\La(\xi)\otimes\La(-\xi'))\cap \Um[{}^\ge\lambda]_{\xi,\xi'}$ and
$\Lmgl_{\xi,\xi'} = (\La(\xi)\otimes\La(-\xi'))\cap
\Um[{}^>\lambda]_{\xi,\xi'}$. Let $\Lm[\lambda]_{\xi,\xi'} =
\Lmgel_{\xi,\xi'}/\Lmgl_{\xi,\xi'}$. Let $\Lm[\lambda]$ be the
$\A_0$-submodule of $\Um[\lambda]$ consisting of elements whose images 
under $\Um[\lambda]\to \Um[\lambda]_{\xi,\xi'}$ are in
$\Lm[\lambda]_{\xi,\xi'}$ for any $\xi,\xi'\in P_+$.
Let $\Lmgel$ be the $\A_0$-submodule of $\Umgel$ consisting of
elements whose images under $\Umgel\to \Umgel_{\xi,\xi'}$ are in
$\Lmgel_{\xi,\xi'}$ for any $\xi,\xi'\in P_+$. By the definition of
$\Lm$ (see \cite[Part IV]{Lu-Book}), this is equal to
$\Lm\cap\Um[{}^\ge\lambda]$. We define $\Lmgl$ in the similar way.
\end{Definition}

\begin{comment}
In an earlier version:

\begin{Definition}
Let $\Lmgel = \Lm\cap\Um[{}^\ge\lambda]$ (resp. $\Lmgl =
\Lm\cap\Um[{}^>\lambda]$).  Define $\Lm[\lambda]$ to be the
$\Azero$-submodule of $\Um[\lambda]$ given by $\Lmgel/\Lmgl.$
\end{Definition}

A priori, it is not clear whether $\Lm[\lambda]$ is invariant under
$\te_i$, $\tf_i$. The trouble is that we do not have well-defined
$\te_i$, $\tf_i$ on $\Um$.

In turn, in this new definition it is not clear whether $\Lmgel\to
\Lmgel_{\xi,\xi'}$ is surjective or not.
\end{comment}

\begin{Lemma}\label{lem:crylatt}
$\Psi_{b'}(\La(\lambda)) \subset \Lm[\lambda]$
for $b'\in\B_W(\lambda)$.
\end{Lemma}

\begin{proof}
The homomorphism $\Psi_{b'}\colon V(\lambda)\to \Um[\lambda]$
intertwines the operators $\te_i$, $\tf_i$.
Since $\Lm[\lambda]$ is invariant under $\te_i$, $\tf_i$ (see
\lemref{lem:inv}(ii)), it is enough to show that $\Psi_{b'}(S
u_\lambda)\in\Lm[\lambda].$
This follows from \lemref{lem:S} as $\Psi_{b'}(S u_\lambda) = S
G(b')^\#\bmod\Um[{}^>\lambda]$.
\end{proof}

Most importantly, we have $G_\lambda(b,s,b')\in\Lm[\lambda]$.  Also
more generally, \lemref{lem:crylatt} holds for $b'\in\B(\lambda)$
thanks to \lemref{lem:Schurmul}.

\begin{Definition}\label{def:Bmlam}
Let $\Bm[\lambda]$ be the subset of $\Bm$ consisting of elements
which are connected to $\B(\lambda)^\#$ in the crystal graph.
It has an induced bi-crystal structure from $\Bm$. By
\thmref{thm:crystaldesc} $\Bm[\lambda] = \Bm[\mu]$ if $\mu\in\aW\lambda$,
$\Bm[\lambda]\cap\Bm[\mu] = \emptyset$ otherwise.
In particular, $\Bm[\lambda]$ ($\lambda\in \cP$) are mutually disjoint.
\end{Definition}

\begin{rem}
We will see that this definition coincides with a more obvious
definition of $\Bm[\lambda]$, i.e., the set of all elements $b\in\Bm$
such that $G(b)\in \Um[{}^\ge\lambda]$ and such that $G(b)$ acts on
$V(\lambda)$ non--trivially.  We will show, in fact, that
$\Um[\lambda]$ has a bi-crystal basis given by $\Bml$.
\end{rem}
As a
first step we have
\begin{Proposition}\label{prop:bi-cry}
There exists an isomorphism of bi-crystals
\begin{equation*}
\label{eq:bi-cry}
   \Bm[\lambda] \cong \B_W(\lambda)\times 
   \operatorname{Irr}G_\lambda\times \B_W(\lambda)
\end{equation*}
where $\te_i$, $\tf_i$, $\te_i^\#$, $\tf_i^\#$ of
the right hand side are defined by
\begin{equation*}
\begin{gathered}
   \te_i (b,s,b') = (\te_i(b,s), b'), \quad
   \tf_i (b,s,b') = (\tf_i(b,s), b'),
\\
   \te_i^\# (b,s,b') = (b,(\operatorname{id}\times\#)\circ\te_i
      \circ(\operatorname{id}\times\#)(b',s)),
\\
   \tf_i^\# (b,s,b') = (b,(\operatorname{id}\times\#)\circ\tf_i
      \circ(\operatorname{id}\times\#)(b',s)),
\end{gathered}
\end{equation*}
where we use the identification
\(
  \B(\lambda)\cong B_W(\lambda)\times\operatorname{Irr}G_\lambda
\)
of \eqref{eq:ident}. Here $s^\#$ denotes the dual representation of $s$.
\end{Proposition}

\begin{proof}
There is a set-theoretical bijection between the right hand side and
the subset
\begin{equation}\label{eq:base}
 \{ G_\lambda(b,s,b') \bmod q_s \Lm[\lambda] \mid
   b, b'\in\B_W(\lambda), s\in\operatorname{Irr}G_\lambda \}
   \subset \Lm[\lambda]/q_s\Lm[\lambda].
\end{equation}
The $\Um$-module structure on $\Um[\lambda]$ defines operators
$\te_i$, $\tf_i$ on $\Lm[\lambda]/q_s\Lm[\lambda]$. If we write
\begin{equation*}
   G_\lambda(b,s,b') = G(\tilde b)G(b')^\#\bmod\Um[{}^>\lambda]
   = \Psi_{b'}(G_\lambda(\tilde b)),
\end{equation*}
with $\tilde b = (b,s)\in\B(\lambda)$ defined by \eqref{eq:ident}, then
\begin{equation*}
   X_i G_\lambda(b,s,b') = \Psi_{b'}(X_i G_\lambda(\tilde b))
   \equiv \Psi_{b'}(G_\lambda(X_i \tilde b)) \mod q_s\Lm[\lambda],
\end{equation*}
for $X_i = \te_i$ or $\tf_i$. Here we have used \lemref{lem:crylatt}
for the second equality. This shows that the bijection respects the
crystal structure. The identification of $\#$-crystal structure
follows from the above discussion and
\begin{equation*}
   G_\lambda(b,s,b')^\# = G(b') S^\# G(b)^\#\bmod\Um[{}^>\lambda],
\end{equation*}
together with \lemref{lem:dual}. We identify $\B_W(\lambda)\times
\operatorname{Irr}G_\lambda\times \B_W(\lambda)$ with \eqref{eq:base}
hereafter.

Consider the natural $\Q$-linear map
\begin{equation*}
   \pi\colon\Lmgel/q_s\Lmgel
   \to 
   \Lm[\lambda]/q_s\Lm[\lambda],
\end{equation*}
induced by the projection $\Um[{}^\ge\lambda]\to \Um[\lambda]$.
Recall that the Kashiwara operators $\te_i$, $\tf_i$ on $\Lm/q_s\Lm$
are defined so that they are compatible with projections
$\Lm/q_s\Lm\to
\La(\xi)\otimes\La(-\xi')/q_s(\La(\xi)\otimes\La(-\xi'))$ induced from
$\Um\to V(\xi)\otimes V(-\xi')$ for all $\xi,\xi'\in P_+$ (see
\cite[Theorem~2.1.2]{Kas94} and \cite[Part~IV]{Lu-Book}). Since
$\Lmgel_{\xi,\xi'}$ is invariant under $\widetilde F_i^{(n)}$ by
\lemref{lem:inv}(ii), the subspace $\Lmgel/q_s\Lmgel\subset
\Lm/q_s\Lm$ is invariant under $\te_i$, $\tf_i$.
The map $\pi$ intertwines $\te_i$, $\tf_i$ as
$\Lmgel_{\xi,\xi'}/q_s\Lmgel_{\xi,\xi'}\to
\Lm[\lambda]_{\xi,\xi'}/q_s\Lm[\lambda]_{\xi,\xi'}$ does so.
It also intertwines $\#$, and hence $\te_i^\#$, $\tf_i^\#$.

Since $G(b)^\#\in\Um[{}^\ge\lambda]$ for $b\in\B(\lambda)$,
$\Bm[\lambda]$ is contained in $\Lmgel/q_s\Lmgel$ as it is invariant
under $\te_i$, $\tf_i$. We have $\pi(b^\#) = (u_\lambda,s,b')$, where
$b = (s^\#, b')$ under \eqref{eq:ident}. Therefore we have
\begin{equation*}
   \pi(\Bm[\lambda])\subset 
   \B_W(\lambda)\times\operatorname{Irr}G_\lambda\times\B_W(\lambda)
   \sqcup \{ 0 \}.
\end{equation*}
Consider $\Ker\,\pi\cap\Bm[\lambda]$. It is invariant under $\te_i$,
$\tf_i$, so every connected component contains an extremal vector, and
in particular an element in $\B(\lambda)^\#$. But $\B(\lambda)^\#$ is
mapped bijectively to $\{ u_\lambda\}\times \Irr G_\lambda\times
\B_W(\lambda)$ as we mentioned already. So $\Ker\,
\pi\cap\Bm[\lambda]$ is the empty set. Thus we have a map
\[
   \pi|_{\Bm[\lambda]}\colon\Bm[\lambda]\to
   \B_W(\lambda)\times\operatorname{Irr}G_\lambda\times\B_W(\lambda).
\]
It is enough to show that this map is bijective since it is clear
that bi-crystal operators are intertwined.

Note that any element of 
\(
   \B_W(\lambda)\times\operatorname{Irr}G_\lambda\times\B_W(\lambda)
\)
can be connected to a point in
\(
   \{u_\lambda\}\times\operatorname{Irr}G_\lambda\times \B_W(\lambda)
\)
in the crystal graph, since the assertion is so for $\B(\lambda)$
(\thmref{thm:crystaldesc}).
Since
\(
   \{u_\lambda\}\times\operatorname{Irr}G_\lambda\times \B_W(\lambda)
\)
is equal to $\pi(\B(\lambda)^\#)$, the map $\pi|_{\Bm[\lambda]}$ is
surjective.

Suppose that $b,b'\in\Bm[\lambda]$ satisfy $\pi(b) = \pi(b')$. We may
assume that $b\in\B(\lambda)^\#$. The condition $\pi(b) = \pi(b')$ and
$\Ker\, \pi\cap\Bm[\lambda] = \emptyset$ imply $b'$ is an extremal
vector with the weight $\wt\, b' = \wt\, b$. This means
$b'\in\B(\lambda)^\#$. Since $\pi$ is bijective on $\B(\lambda)^\#$,
we have $b = b'$.
\end{proof}

Hereafter, we identify $\Bm[\lambda]$ with
$\B_W(\lambda)\times\operatorname{Irr}G_\lambda\times\B_W(\lambda)$ by
\propref{prop:bi-cry}.

\begin{Lemma}\label{lem:union}
$\Bm = \bigsqcup_{\lambda\in P^{0,+}_\cl} \Bm[\lambda]$.
\end{Lemma}

\begin{proof}
Each connected component of $\Bm$ contains an extremal vector
\cite[9.3.4]{Kas94}. By definition, the set of extremal vectors is
equal to
\(
   \bigsqcup_{\lambda\in P^0_\cl} \B(\lambda)^*
    = \bigsqcup_{\lambda\in P^0_\cl} \B(\lambda)^\#.
\)
Furthermore, we have an isomorphism of crystals
\(
   S_w^*\colon \B(\lambda) \xrightarrow{\cong} \B(w\lambda)
\)
for $w\in \aW$. Therefore each component contains a vector in
$\B(\lambda)^\#$ with $\lambda\in \cP$.
Therefore $\Bm = \bigcup_{\lambda\in P^{0,+}_\cl} \Bm[\lambda]$.  By
the remark in Definition~\ref{def:Bmlam}, this is a disjoint union.
\end{proof}

For $\xi,\xi'\in P_+$, the kernel of the surjective homomorphism $\Um
a_{\xi-\xi'}\to V(\xi)\otimes V(-\xi')$ is the left ideal generated by
$f_i^{\langle \xi, h_i\rangle+1} a_{\xi-\xi'}$ and $e_i^{\langle \xi',
  h_i\rangle+1} a_{\xi-\xi'}$. Let $P(\xi,\xi')$ be the set of $\mu\in
\cP$ such that all $f_i^{\langle \xi, h_i\rangle+1} a_{\xi-\xi'}$ and
$e_i^{\langle \xi', h_i\rangle+1} a_{\xi-\xi'}$ act by $0$ on
$V(\mu')$ for any $\mu'\in\cP$ satisfying $\mu'\le\mu$. Therefore a
homomorphism $V(\xi)\otimes V(-\xi')\to \End \, V(\mu)$ is well-defined
if $\mu\in P(\xi,\xi')$. Moreover, any $\mu\in\cP$ is contained in
$P(\xi,\xi')$ for sufficiently large $\xi$, $\xi'$.

Let $\Bm[\mu]_{\xi,\xi'} = \Bm[\mu]\cap (\B(\xi)\otimes \B(-\xi))$.
Then we have $\B(\xi)\otimes \B(-\xi') =\bigsqcup_\mu
\Bm[\mu]_{\xi,\xi'}$ by \lemref{lem:union}.

Let $(b,s,b')\in\Bm[\mu]_{\xi,\xi'}$. We have
$G_\mu(b,s,b')(u_\xi\otimes u_{-\xi'})\in\Lm[\mu]_{\xi,\xi'}$ by
Lemma~\ref{lem:crylatt}. We choose lifts $\tilde
G_\mu(b,s,b')\in\Lm[{}^\ge\mu]_{\xi,\xi'}$.
(Note that we do not have $G(b) S G(b')^\#\in\Lm$ in general.)
Then
\begin{equation*}
    \bigsqcup_{\mu\in P^{0,+}_\cl}
    \left.\left\{ \tilde G_\mu(b,s,b')\, \right|
      (b,s,b')\in\Bm[\mu]_{\xi,\xi'} \right\}
\end{equation*}
is a $\A_0$-basis of $\La(\xi)\otimes \La(\xi')$, as it induces a
$\Q$-basis of $\La(\xi)\otimes \La(\xi')/q_s(\La(\xi)\otimes
\La(\xi'))$. (The transition matrix between this basis and
$\B(\xi)\otimes \B(-\xi')$ is upper triangular with $1$'s on the
diagonal with respect to the block decomposition induced by $\mu$ and
the order $\le$.)

\begin{Lemma}\label{lem:sum}
Let $\tilde G_\mu(b,s,b')\in\Lm[\mu]_{\xi,\xi'}$ as above. If
\begin{equation*}
   x = \sum_{\mu\in P^{0,+}_\cl}\sum_{(b,s,b')\in\Bm[\mu]_{\xi,\xi'}}
       a^\mu_{b,s,b'} \tilde G_\mu(b,s,b')\in\Um[{}^\ge\lambda]_{\xi,\xi'},
\end{equation*}
then $a_{b,s,b'}^\mu = 0$ if $\mu\ngeq\lambda$ and $\mu\in P(\xi,\xi')$.
\end{Lemma}

\begin{proof}
We take any minimal element $\mu_0$ with respect to $\le$ among
$\mu$'s with $a^\mu_{b,s,b'}\neq 0$ for some
$(b,s,b')\in\Bm[\mu]$.  Assume $\mu_0\ngeq\lambda$
and $\mu_0\in P(\xi,\xi')$. Then $V(\xi)\otimes V(-\xi')\to \End V(\mu_0)$
is well-defined, and the image of $x$ is $0$ by the assumption.
On the other hand, any $\mu$ with $a^\mu_{b,s,b'}\neq 0$ for some
$(b,s,b')\in\Bm[\mu]$ satisfies $\mu\nless\mu_0$ by the minimality
condition. Therefore among terms in $x$, only those summands with $\mu
= \mu_0$ act nontrivially on $V(\mu_0)$. Therefore, as operators on
$V(\mu_0)$, we have
\begin{equation*}
   x|_{V(\mu_0)} =
       \sum_{(b,s,b')\in\Bm[\mu_0]_{\xi,\xi'}} a^{\mu_0}_{b,s,b'}
       \tilde G_{\mu_0}(b,s,b')|_{V(\mu_0)} = 0.
\end{equation*}
However the operator $\tilde G_{\mu_0}(b,s,b')|_{V(\mu_0)}$ does
not depend on the choice of the lift $\tilde G_{\mu_0}(b,s,b')$ of
$G_{\mu_0}(b,s,b')$, and we have $a^{\mu_0}_{b,s,b'} = 0$ for
$(b,s,b')\in\Bm[\mu_0]_{\xi,\xi'}$, as $G_{\mu_0}(b,s,b')$ are
linearly independent (\lemref{lem:linind}). This is a contradiction.
Thus we have $\mu_0\ge\lambda$ or $\mu_0\notin P(\xi,\xi')$. Since
$\mu_0$ is any minimal element, this completes the proof.
\end{proof}

\begin{Proposition}\label{prop:global}
Let $\beta = (b,s,b')\in\Bm[\lambda]$. Considering $\beta$ as an
element of $\Bm$, let $G(\beta)$ be the corresponding element of the
global crystal basis of $\Um$. Then we have
\begin{equation*}
   G(\beta)\in\Um[{}^\ge\lambda], \qquad
   G(\beta)\bmod \Um[{}^>\lambda] = G_\lambda(b,s,b').
\end{equation*}  
\end{Proposition}

\begin{proof}
We will show the assertion by showing
\begin{aenume}
\item $G(\beta)(u_\xi\otimes u_{-\xi'})$ acts by $0$ on $V(\mu)$ if
$\mu\ngeq\lambda$ and $\mu\in P(\xi,\xi')$.
\item $G(\beta)(u_\xi\otimes u_{-\xi'})$ and 
$G_\lambda(b,s,b')(u_\xi\otimes u_{-\xi'})$ define the same operator
on $V(\lambda)$ if $\lambda\in P(\xi,\xi')$,  
\end{aenume}
for any $\xi,\xi'\in P_+$. We fix $\xi$, $\xi'$ hereafter.

Let $(b_1,s_1,b'_1)\in\Bm[\mu]_{\xi,\xi'}$. We choose and fix a lift
$\tilde G_\mu(b_1,s_1,b'_1)\in\Lm[{}^\ge\mu]_{\xi,\xi'}$ of
$G_\mu(b_1,s_1,b'_1)(u_\xi\otimes u_{-\xi'})$ as in \lemref{lem:sum}.
We write
\begin{equation}\label{eq:globalbase}
   G(\beta)(u_\xi\otimes u_{-\xi'})
   = \sum_{\mu\in \cP}\sum_{(b_1,s_1,b'_1)\in\Bm[\mu]_{\xi,\xi'}}
        a^\mu_{b_1,s_1,b'_1} \tilde G_\mu(b_1,s_1,b'_1)
\end{equation}
for some $a^\mu_{b_1,s_1,b'_1}\in\A_0$.

\begin{Claim}
{\rm (i)} If $\mu\ngeq\lambda$ and $\mu\in P(\xi,\xi')$, then
\(
   a^\mu_{b_1,s_1,b'_1} \in q_s\A_0
\)
for $(b_1,s_1,b'_1)\in\Bm[\mu]_{\xi,\xi'}$.

{\rm (ii)} If $\lambda\in P(\xi,\xi')$, then
\(
   a^\lambda_{b_1,s_1,b'_1} \in \delta_{(b,s,b'),(b_1,s_1,b'_1)} + q_s\A_0
\)
for $(b_1,s_1,b'_1)\in\Bm[\lambda]_{\xi,\xi'}$.
\end{Claim}

We can replace $G(\beta)(u_\xi\otimes u_{-\xi'})$ by an element which
is equal to it modulo $q_s(\La(\xi)\otimes\La(-\xi'))$ in showing (i)
and (ii). Let us write
\begin{equation*}
   (b,s) = X_1\cdots X_N (u_\lambda,s_0)
\end{equation*}
for some $X_i = \te_j$ or $\tf_j$. Then $G(\beta)(u_\xi\otimes
u_{-\xi'})$ is equal to 
\(
   X_1\cdots X_N S_0 G(b')^\#(u_\xi\otimes u_{-\xi'})
\)
modulo $q_s(\La(\xi)\otimes\La(-\xi'))$. We show the assertion for this
element.
By \lemref{lem:inv}, it is contained in $\Umgel_{\xi,\xi'}$.
By \lemref{lem:sum} we have (i).

To show (ii), suppose $\lambda\in P(\xi,\xi')$ and $V(\xi)\otimes
V(\xi')\to \End \, V(\lambda)$ is well-defined. If $\mu\notin
P(\xi,\xi')$, then $\lambda\ngeq\mu$ by the assumption. We combine
this with the above discussion to have $a^\mu_{b_1,s_1,b'_1} = 0$
unless $\lambda\ngtr\mu$. Therefore among terms in
\eqref{eq:globalbase}, only those summands with $\mu = \lambda$ act
nontrivially on $V(\lambda)$. Thus we have
\begin{multline*}
   X_1\cdots X_N S_0 G(b')^\#(u_\xi\otimes u_{-\xi'})
\\
  \equiv \sum_{(b_1,s_1,b'_1)\in\Bm[\lambda]_{\xi,\xi'}}
       a^\lambda_{b_1,s_1,b'_1} \tilde G_\lambda(b_1,s_1,b'_1)
    \mod \Um[{}^>\lambda]_{\xi,\xi'}.
\end{multline*}
On the other hand, we have
\begin{equation*}
  X_1\cdots X_N S_0 G(b')^\# \bmod \Um[{}^>\lambda]
   = \Psi_{b'}(X_1\cdots X_N S_0 u_\lambda).
\end{equation*}
Since $\Psi_{b'} (\Lm(\lambda)) \subset \Lml$ by \lemref{lem:crylatt},
this is equal to $G_\lambda(b,s,b')$ modulo
$q_s\Lm[\lambda]$. This completes the proof of the claim.

We take any minimal element $\mu_0$ with respect to $\le$ among
$\mu$'s with $a^\mu_{b_1,s_1,b'_1}\neq 0$ for some
$(b_1,s_1,b'_1)\in\Bm[\mu]$.  We obtain a contradiction under the
assumption $\mu_0\ngeq\lambda$ and $\mu_0\in P(\xi,\xi')$. As in the
proof of \lemref{lem:sum}, only summands with $\mu = \mu_0$ act
nontrivially on $V(\mu_0)$. Hence, for $b''\in\B(\mu_0)$ we have
\begin{equation*}
\begin{split}
   G(\beta) G_{\mu_0}(b'') &=
   \sum_{(b_1,s_1,b'_1)\in\Bm[\mu_0]} a^{\mu_0}_{b_1,s_1,b'_1} 
   \tilde G_{\mu_0}(b_1,s_1,b'_1) G_{\mu_0}(b'')
\\
     & = \sum_{(b_1,s_1,b'_1)\in\Bm[\mu_0]} a^{\mu_0}_{b_1,s_1,b'_1} 
     G(b_1) S_1 G(b'_1)^\# G_{\mu_0}(b'').
\end{split}
\end{equation*}
{}From $\overline{G(\beta) G_{\mu_0}(b'')} = G(\beta) G_{\mu_0}(b'')$,
we have
\begin{equation*}
   0 = \sum_{(b_1,s_1,b'_1)\in\Bm[\mu_0]}
     \left(a^{\mu_0}_{b_1,s_1,b'_1} -
     \overline{a^{\mu_0}_{b_1,s_1,b'_1}}\right) G(b_1) S_1 G(b'_1)^\#
     G_{\mu_0}(b'').
\end{equation*}
Arguing as in the proof of \lemref{lem:linind}, we get
\begin{equation*}
   a^{\mu_0}_{b_1,s_1,b'_1} = \overline{a^{\mu_0}_{b_1,s_1,b'_1}}.
\end{equation*}
However, by the assumption $\mu_0\ngeq\lambda$ and $\mu_0\in
P(\xi,\xi')$, we have $a^{\mu_0}_{b_1,s_1,b'_1} \in q_s\A_0$ using
Claim~(i). Thus this equality is impossible unless
$a^{\mu_0}_{b_1,s_1,b'_1} = 0$. This is a contradiction, and we get
(a).

Similarly, if $\lambda\in P(\xi,\xi')$, we have
\(
  a^{\lambda}_{b_1,s_1,b'_1} = \overline{a^{\lambda}_{b_1,s_1,b'_1}}.
\)
Then Claim~(ii) implies 
\(
  a^{\lambda}_{b_1,s_1, b'_1} = \delta_{(b,s,b'),(b_1,s_1,b'_1)}.
\)
We get (b).
\end{proof}

Now we have
\begin{Theorem}\label{thm:Umcry}
\textup{(i)} $\Um[\lambda]$ has a crystal base \( \left(\Lm[\lambda],
\Bm[\lambda]\right) \), where $\Bm[\lambda]$ is identified with a
subset of $\Lm[\lambda]/q_s\Lm[\lambda]$ as in the proof of
\propref{prop:bi-cry}.

\textup{(ii)} 
\(
  \left(\Lm[\lambda], \overline{\Lm[\lambda]}, \Umz[\lambda]\right)
\)
is a balanced triple of $\Um[\lambda]$ and
\(
  \{ G(b)\bmod \Um[{}^>\lambda] \mid b\in \Bm[\lambda] \}
\)
is its global crystal basis.
\end{Theorem}

\begin{proof}
By \propref{prop:global}, we can take $G(\beta)(u_\xi\otimes
u_{-\xi})$ for the lift $\tilde G_\lambda(b,s,b')$ in
\lemref{lem:sum}.  Therefore \lemref{lem:sum} implies that
\(
   \bigsqcup_{\mu\ge\lambda}\{ G(\beta)\mid \beta\in\Bm[\mu]\}
\)
is a basis of $\Um[{}^\ge\lambda]$. Hence $G(\Bm[\lambda])\bmod
\Um[{}^>\lambda]$ is a basis of $\Um[\lambda]$. Other axioms of the
global crystal basis obviously hold.
\end{proof}

\subsection{Cell structure of $\Um$}

We recall \cite{Lu-v} the definition of cells in $\Um$ with respect to
the global basis $\Gm = \{ G(\beta) \mid \beta \in \Bm \}.$
Let $\mathscr F$ be the collection of all subsets $K \subset \Bm$ with
the property: the $\Q(q_s)$--subspace of $\Um$ spanned by $G(K)$ is a
two--sided ideal of $\Um.$ If $\beta, \beta' \in \Bm$, we say that
$\beta \preceq \beta'$ if $\beta \in \cap_{K \in \mathcal F, \beta'
\in K} K$.  We say that $\beta \sim \beta'$ if $\beta \preceq \beta'$
and $\beta' \preceq \beta.$ The equivalence classes of $\sim$ are
called two--sided cells. Similarly, by considering left ideals or
right ideals above, we define $\preceq_{L}$ or $\preceq_R$.  The
equivalence classes are then called left cells or right cells
respectively.
Another equivalent definition of $\preceq$ is as follows: Let
$$G(\beta)G(\beta') = \sum_{\beta''} c_{\beta \beta'}^{\beta''}(q_s)
G(\beta'')$$
define the structure constants $c_{\beta\beta'}^{\beta''}(q_s) \in
\ca$ of $\Um$ with respect to the global basis.  For $\beta, \beta'
\in \Bm$ we say $\beta \preceq_L \beta'$ (resp. $\preceq_R$) if there
is a sequence $\be_1 = \be', \be_2, \dots, \be_N = \be$ in $\Bm$ and a
sequence $\gamma_1, \dots, \gamma_{N-1} \in \Bm$ such that
$c_{\gamma_s, \be_s}^{\be_{s+1}} \neq 0$ (resp.
$c_{\be_s,\gamma_s}^{\be_{s+1}} \neq 0$) for $s = 1, \dots, N-1.$ We
write $\be \preceq \be'$ if either of the above structure constants is
non--zero for all $\be_s, \gamma_s, \beta_{s+1}, s = 1, \dots, N-1.$

\begin{Proposition}\label{prop:twocell}
\textup{(i)} $\Bml$ is a two--sided cell of $\Bm$. 

\textup{(ii)} For any $b_2\in\B_W(\lambda)$, $\{ (b_1, s, b_2)\in \Bml \mid
s\in\Irr G_\lambda, b_1\in \B_W(\lambda)\}$ is a left cell.

\textup{(iii)} For any $b_1\in\B_W(\lambda)$, $\{ (b_1, s, b_2)\in \Bml \mid
s\in\Irr G_\lambda, b_2\in \B_W(\lambda)\}$ is a right cell.
\end{Proposition}

\begin{proof}
(i) Since $\Umgel$, $\Umgl$ are two--sided ideals of $\Um$, $\Bml$ is
a union of two--sided cells.
Let $b\in \Bm$. By \cite[6.4.3]{Kas94}, $f_i G(\beta)$ is equal to
$\left[\varepsilon_i(\beta) + 1\right] G(\tf_i \beta)$ plus a linear
combination of elements in $\Bm$ different from $\tf_i
\beta$. Therefore $G(\tf_i \beta)\preceq_L G(\beta)$ if $\tf_i
\beta\neq 0$. Similarly we have $G(\te_i \beta)\preceq_L G(\beta)$ if
$\te_i \beta\neq 0$. Thus if $\beta$ and $\beta'$ are in the same
connected component of $\Bm$, we have $G(\beta)\sim_L
G(\beta')$. Taking $\#$, we conclude that if $\beta$ and $\beta'$ are
in the same connected component of $\Bm$ as $\#$--crystal, we have
$G(\beta)\sim_R G(\beta')$. By \thmref{thm:crystaldesc} we have the
assertion if $a_\lambda\sim S a_\lambda$ for $s\in \Irr
G_\lambda$. This follows from a consideration of two--sided cells of
the based ring $(R(G_\lambda), \Irr G_\lambda)$, where $R(G_\lambda)$
is the representation ring of $G_\lambda$ and $\Irr G_\lambda$ is
considered as its basis. By using the Pieri formula, we can easily
check that it consists of a single two--sided cell, and hence our
assertion. The proofs of (ii), (iii) are contained in the above proof.
\end{proof}

Remark $\Bm = \bigsqcup_{\lambda\in\cP} \Bml$ by \lemref{lem:union}.
Therefore the above proposition gives a complete description of
two--sided, left, right cells of $\Bm$.

Next we define a function on $\Gm$ which will allow us construct a
limit algebra $\Um_0 = \bigoplus_{\lambda \in \cP} \Uml_0$ as $q
\rightarrow 0$.

\begin{Definition}  Let $\beta = (b, s, b') \in \Bml$.   Define 
$a(\beta) = -(\wt(b'), 2 \la + \wt(b'))/2.$
\end{Definition}

\begin{Remark}
$a(\beta)$ is equal to the half of the dimension of the quiver variety
$\mathfrak M(\mathbf v,\mathbf w)$ (or the dimension of the lagrangian
subvariety $\mathfrak L(\mathbf v,\mathbf w)$), where $\mathbf v =
-\wt(b')$, $\mathbf w = \lambda$. (See e.g.\ \cite[2.3.2(4)]{Nak00}.)
The inner product $(\ ,\ )$ induced on $P^0_\cl$ is equal to the
standard inner product for the untwisted $ADE$ case (see
\cite[Corollary~6.4]{Kac}), which is the only case for which quiver
varieties are defined.
\end{Remark}

\begin{Lemma}[\protect{cf.\ \cite[1.4, 4.11]{Lu-v}}] \label{lem:P1}
Let $\beta = (b, s, b') \in \Bml$. 
The following holds: 

\textup{(i)} $(q^{a(\beta)} \beta) \Lm[\lambda] \subset \Lm[\lambda]$
and for each $\beta$, $a(\beta) \ge 0$ is smallest integer with this
property.

\textup{(ii)} for a two--sided cell $\Bml$ and any $\lambda_1 \in \cP$, 
the restriction of $a$ to $\{ \beta\in \Bml \mid G(\beta)\in \U
a_{\lambda_1}\}$ is constant.
\end{Lemma}

\begin{proof} First note that using Lemmas~\ref{lem:Schurmul},
\ref{lem:mul2} and \propref{prop:global} we have the following
equalities in $\Uml = \Umgel/\Umgl$:
{\allowdisplaybreaks
\begin{align*}
  & G(b,s,b') G(b_1, s_1, b'_1) \bmod\Umgl 
\\
  =\; & q^{-a(\beta)}
   \sum_{s'\in\operatorname{Irr}G_\lambda} (G(b_1) S_1 u_\lambda,
   G(b') S' u_\lambda)_\lambda G(b) S S' G(b'_1)^\# \bmod \Umgl
   \notag
\\ = \; & q^{-a(\beta)}
   \sum_{s', s''\in\operatorname{Irr}G_\lambda} c_{s,s'}^{s''}
   (G(b_1) S_1 u_\lambda, G(b') S' u_\lambda)_\lambda G(b) S''
   G(b'_1)^\# \bmod \Umgl \notag
\end{align*}
where $s s' = \sum_{s''} c_{ss'}^{s''} s''$ and $c_{ss'}^{s''} \in
\Z_{\ge 0}.$} By \thmref{thm:main2} (and \lemref{lem:S}) this is
contained in $q^{-a(\beta)} \Lml$.
(i) now follows by multiplying both sides by $q^{a(\beta)}$.  Note
that when $b'_1 = b_2$ the product is in $\Lml \setminus q_s\Lml$ by
\thmref{thm:main2} and so $a(\beta)$ is the smallest integer with this
property.  We check that $a(\beta)$ is positive. We have $G(b,s,b')
G(b',1,b) = q^{-a(\beta)} G(b,s,b) \bmod (\Lmgl + q_s\Lmgel)$ by
\thmref{thm:main2}. It follows $c(q_s) = c_{\beta,
(b',1,b)}^{(b,s,b)}(q_s) = q^{-a(\beta)} \mod q_s \Z[q_s].$ Taking
$\setbox5=\hbox{A}\overline{\rule{0mm}{\ht5}\hspace*{\wd5}}\,$ of both
sides and using the
$\setbox5=\hbox{A}\overline{\rule{0mm}{\ht5}\hspace*{\wd5}}\,$
invariance of the global basis we also have $ c(q_s) = q^{a(\beta)}
\mod q_s^{-1} \Z[q_s^{-1}].$ This implies $a(\beta) \ge 0.$ (ii) is
clear from the definition: If $\beta = (b,s,b')\in\Bml$ satisfies
$G(\beta) \in \Um a_{\lambda_1}$, we have $\wt(b') = \lambda_1 -
\lambda$.
\end{proof}

Using Lemma \ref{lem:P1} for any two--sided cell, we define the ring
$\Uml_0$ which is a ``limit as $q \rightarrow 0$.''  For $\beta \in
\Bml$ define $\hb = q^{a(\beta)} G(\beta).$ Then $\{\hb \mid \be \in \Bml
\}$ is a $\Q(q_s)$--basis of $\Uml.$ For $\beta, \beta' \in \Bml$ such
that $G(\beta') \in \U a_{\lambda_1}$, we have $\hb \hb' = \sum_{\hb''}
q^{a(\be)} c_{\be,\be'}^{\be''} \hb''$ in $\Uml.$ (This follows from
$a(\beta') = a(\be'')$, which follows from Lemma~\ref{lem:P1}(ii)
above, since the only non--zero coefficients $c_{\be,\be'}^{\be''}$
are those for which $G(\be'') \in \U a_{\lambda_1}$).  Since
$q^{a(\beta)} c_{\be,\be'}^{\be''} \in \Z[q_s]$ by \lemref{lem:P1}(i),
the $\Z[q_s]$--submodule $\Uml^-$ of $\Uml$ generated by $\{\hb \mid
\be \in \Gml\}$ is a $\Z[q_s]$--subalgebra of $\Uml$.  We define
$\Uml_0 = \Uml^- / q_s \Uml^-$.  Define $t_\beta$ to be the image of
$\hb$ in $\Uml_0.$ Then $\Uml_0$ is a ring with a $\Z$--basis $\{
t_\beta \mid \beta \in \Gml \}.$ If we denote by $\{ \gamma_{\be,
\be'}^{\be''} \mid \be, \be', \be'' \in \Gml\}$ the structure
constants of $\Uml_0$ then $q^{a(\be)} c_{\be,\be'}^{\be''} \equiv
\gamma_{\be, \be'}^{\be''} \mod q_s\Z[q_s].$ We define $\Um_0 =
\bigoplus \Uml_0$ to be the direct sum of these rings.

\begin{Lemma}\label{lem:sharp}
$G(b,s,b')^\# = G(b', s^\#,b)$ for $(b,s,b')\in\Bml$.
\end{Lemma}

\begin{proof}
Note $\left(G(b)SG(b')^\#\right)^\# = G(b') S^\# G(b)^\#$. By
\lemref{lem:dual} $S^\#$ corresponds to the dual representation 
$s^\#$ modulo $\Uml$. We have $G_\lambda(b,s,b')^\# =
G_\lambda(b', s^\#, b)$ by the definition, then the assertion follows by
\propref{prop:global}.
\end{proof}

\begin{Definition}
Let $\lambda \in \cP.$ Let \( \mathcal D_{\Bm[\lambda]} = \{ (b,1,b)
\mid b\in \B_W(\lambda) \} \). Denote by $\Um(0)$ the subalgebra of
$\Um$ generated by $a_\lambda x a_\lambda$ for $x \in \U$.
\end{Definition}

\begin{Remarks}
(i) By the definition and \lemref{lem:sharp}, it is clear that
$G(\beta)^\# = G(\beta)$ for $\beta\in\mathcal D_{\Bm[\lambda]}$. It
is also clear that $G(\beta)$ is in $\Um(0)$.

\begin{comment}
These are \cite[Conj.~5.3(a)(c)]{Lu-v}.
\end{comment}

(ii) It is not true that $G(\beta)^\# = G(\beta)$ implies $\beta\in\mathcal
D_{\Bm[\lambda]}$. Consider $(b,s,b)$ with $s^\# = s$ and $s\neq 1$,
e.g. $s(z) = (z_1^2 + z_1z_2 + z_2^2)/z_1z_2$.

(iii) The identification \eqref{eq:ident} depends on the choice of the
section $\B_W(\lambda)\to \B_0(\lambda)$. But the ambiguity of
monomials $p(z) = \prod_i (z_{i,1}\cdots z_{i,\lambda_i})^{r_i}$
cancels in $G(\beta) = G(b,1,b)$. Therefore the bijection
$\B_W(\lambda)\to \mathcal D_{\Bm[\lambda]}$ is independent of the
section.
\end{Remarks}

The following gives a characterization of $\mathcal D_{\Bm[\lambda]}$: 
\begin{Proposition}[\protect{cf.\ \cite[Conj.~5.3(b)]{Lu-v}}]
Let $\lambda_1\in P^0_\cl$. If $G(\beta)\in \Um(0)a_{\lambda_1}\cap
G(\Bm[\lambda])$, then $q^{-a(\beta)}(a_{\lambda_1},G(\beta))\equiv 1
\mod q_s\Z[q_s]$ if $\beta\in\mathcal D_{\Bm[\lambda]}$
and $\equiv 0$ otherwise.
\end{Proposition}

\begin{proof}
By definition, $a$ takes the constant $(\lambda-\lambda_1,
\lambda+\lambda_1)/2$ on $\Um(0)a_{\lambda_1}\cap G(\Bm[\lambda])$. We
denote this constant by $a_0$. If $\mu > \lambda$, then
\begin{equation*}
   (\mu-\lambda_1,\mu+\lambda_1) - (\lambda-\lambda_1,\lambda+\lambda_1)
   = (\mu-\lambda, \mu+\lambda) > 0. 
\end{equation*}
Thus our assertion is equivalent to saying that if $G(\beta)\in
\Um(0)a_{\lambda_1}\cap \bigsqcup_{\mu\ge\lambda}G(\Bm[\mu])$, then
$q^{-a_0}(a_{\lambda_1},G(\beta))\equiv 1 \mod q_s\Z[q_s]$ if
$\beta\in\mathcal D_{\Bm[\lambda]}$ and $\equiv 0$ otherwise. We will
check this.

In this proof we replace $\Um$ by $\bigoplus_{\lambda\in P} \U
a_\lambda$ the modified quantum enveloping algebra defined for
$P$. Then by \cite[26.2.3]{Lu-Book}, $(\ ,\ )$ is the limit of the
inner product on $V(\xi)\otimes V(-\xi')$ (denoted also by $(\ ,\ )
$), where $\xi-\xi' = \lambda_1$ and $\xi$, $\xi'$ tend to
$\infty$.

For each $G(\beta) = G(b, s, b')\in \Um(0)a_{\lambda_1}\cap
\bigsqcup_{\mu\ge\lambda}G(\Bm[\mu])$, we choose and fix an expression
$(b,s) = X_1 X_2 \cdots X_N (u_\lambda,s_0)$ where $s_0\in\Irr
G_\lambda$, $X_i = \te_i$ or $\tf_i$. Then
{\allowdisplaybreaks
\begin{equation*}
\begin{split}
   & q^{-a(\beta)}(u_{\xi}\otimes u_{-\xi'},
   X_1 X_2 \cdots X_N S_0 G(b')^\# (u_{\xi}\otimes u_{-\xi'}))
\\
  =\; &
  \left((X_1 X_2 \cdots X_N S_0)^\# (u_{\xi}\otimes u_{-\xi'}),
   G(b')^\# (u_{\xi}\otimes u_{-\xi'})\right)
\\
  \equiv \; & \left(G(\tilde b) (u_{\xi}\otimes u_{-\xi'}),
   G(b')^\# (u_{\xi}\otimes u_{-\xi'})\right)
   \mod q_s\Z[q_s]
\\
   \equiv \; & \delta_{(b,s), (b', 1)}   \mod q_s\Z[q_s],
\end{split}
\end{equation*}
where $\tilde b = (b,s)$ by \eqref{eq:ident}.} Here we have
used \eqref{eq:psi} in the first equality, and the almost orthonomal
property of the global crystal basis (\cite[26.3.1]{Lu-Book}) in the
second and third equalities.
Since $\{ G(\beta) = G(b,s,b') \} =
\Um(0)a_{\lambda_1}\cap \bigsqcup_{\mu\ge\lambda}G(\Bm[\mu])$
is a $\A_0$-basis of $\Lm[{}^\ge\lambda]$, the set of
corresponding elements
\(
   \{ X_1 X_2 \cdots X_N S_0 G(b')^\#\linebreak[0]
   (u_{\xi}\otimes u_{-\xi'})
   \}
\)
spans $(\Um(0)a_{\lambda_1}\cap \Lm({}^\ge\lambda))(u_{\xi}\otimes
u_{-\xi'})$.
Therefore the above together with our previous remark
$a(\beta) \ge a_0$ implies
\[
   q^{-a_0}(u_{\xi}\otimes u_{-\xi'},
   (\Um(0)a_{\lambda_1}\cap \Lm({}^\ge\lambda))
   (u_{\xi}\otimes u_{-\xi'}))\in \Z[q_s].
\]
Since
\(
   G(b,s,b')(u_{\xi}\otimes u_{-\xi'})
   \equiv
   X_1 X_2 \cdots X_N S_0 G(b')^\# (u_{\xi}\otimes u_{-\xi'})
   \mod \Lm[{}^\ge\lambda](u_{\xi}\otimes u_{-\xi'}),
\)
the above shows
\[
   q^{-a_0}(u_{\xi}\otimes u_{-\xi'},
   G(b,s,b') (u_{\xi}\otimes u_{-\xi'}))
   \equiv \delta_{\lambda,\mu} \delta_{(b,s), (b', 1)} \mod q_s\Z[q_s].
\]
Taking limit $\xi,\xi'\to\infty$, we get the assertion.
\end{proof}

\begin{comment}
We need to modify the statement to include $\mu\ge\lambda$ for the
equality second last. $\Um(0)a_{\lambda_1}\cap \Lm[{}^\ge\lambda]$
contains element from $G(\Bm[\mu])$.
\end{comment}

\begin{Lemma}[\protect{cf.\ \cite[1.5]{Lu-v}}]  \label{lem:P2}
Let $d_1 = G(b_1, 1, b_1), d_2 = G(b_2, 1, b_2)$ be in $\mathcal
D_{\Bml}$. Let $\beta = G(b_3, s, b'_3).$ Then in $\Uml_0$,
$t_{d_1} t_\beta t_{d_2} = 
\delta_{b_1, b_3} \delta_{b'_3, b_2} t_\beta.
$
\end{Lemma}
\begin{proof} 
By \lemref{lem:mul2} we have 
\begin{align*} G(b_1, 1, b_1) & G(b_3,s,b'_3) \\
 & \equiv q^{-a(d_1)} \sum_{s''\in\operatorname{Irr}G_\lambda}
   (G(b_3) S u_\lambda, G(b_1) S'' u_\lambda)_\lambda G(b_1)
   S'' G(b'_3)^\# \bmod \Lmgl \\
 & \equiv q^{-a(d_1)} 
  \delta_{b_1, b_3} G(b_1) S G(b'_3)^\# \bmod (\Lmgl + q_s \Lmgel).
\end{align*}
The second equivalence follows from \thmref{thm:main2}.  A similar
calculation shows
\begin{equation*} 
G_\lambda(b_3,s,b'_3) G_\lambda(b_2, 1, b_2) \equiv
q^{-a(\beta)} \delta_{b'_3,b_2} G(b_3, s, b'_3) \bmod (\Lmgl + q_s \Lmgel).
\end{equation*}
Since $a(\beta) = a(d_2)$ if $b_3' = b_2$ the result follows
by multiplying both sides by $q^{2a(d_2)}$.
\end{proof}

\begin{rem} \lemref{lem:P2} says that the $\Z$--ring $\Uml_0$ 
has a generalized unit which is compatible with the basis $t_{\beta},
\beta \in \Gml.$ In particular, the ring $\Uml_0$ has an identity, $1
= \sum_{b_1 \in B_W(\lambda)} t_{G(b_1,1,b_1)}$.
\end{rem}

$\Uml$ is both a left $\Um$--module and a right $\Um$--module where
the respective module structures are given by:
\begin{equation} \label{sconsts}
G(\beta) G(\beta') = \sum_{\beta'' \in \Bml} c_{\beta
\beta'}^{\beta''}(q_s) G(\beta''), 
\end{equation}
when $\beta \in \Bm, \beta' \in \Bml$ in the first case and 
when $\beta \in \Bml, \beta' \in \Bm$ in the second case. 

We define $M_{\Uml}$ to be the vector space spanned over
$\Q(q_s,q'_s)$ by $\{ G(\beta) \mid \beta \in \Bml\}$.  $M_{\Uml}$ is
a $\Um$--bimodule where the left action is given by \eqref{sconsts}
and the right action is given by \eqref{sconsts} where $c_{\beta
\beta'}^{\beta''}(q_s)$ is replace by $c_{\beta
\beta'}^{\beta''}(q_s')$.  We now show that the left and right module
structures commute, i.e., $G(\beta_1) (G(\beta_2) G(\beta_3))=
(G(\beta_1) G(\beta_2)) G(\beta_3)$ where $\beta_2 \in \Bml$ and
$\beta_1, \beta_3 \in \Bm$.  In terms of the structure constants this
is equivalent to

\begin{Lemma}[cf.\ \protect{\cite[1.7, 4.15]{Lu-v}}] \label{lem:P3}
Let $\beta_4, \beta_2 \in \Bml$.  
Let $\beta_1, \beta_3 \in \Bm.$ Then
\begin{equation*} 
  \sum_{\beta \in \Bml} c_{\be_1, \be}^{\be_4}(q_s) c_{\be_2,
\be_3}^{\be}(q'_s) = \sum_{\beta \in \Bml} c_{\be_1, \be_2}^\be(q_s)
c_{\be, \be_3}^{\be_4}(q'_s).
\end{equation*}
\end{Lemma}

\begin{proof}  We fix $\be_i = (b_i, s_i, b'_i)$ for $i = 2, 4$.
For $\beta_j \in \Bm$ ($j=1,3$) and $b_i \in B_W(\lambda)$ ($i=2,4$),
we define $g_{\beta_j,b_i}^{(s,b)}\in\ca$ by
\begin{equation} \label{uconsts}
  G(\beta_j) G_\lambda(1,b_i)
  = \sum_{(s,b) \in \B(\lambda)} g_{\beta_j,b_i}^{(s,b)} G_\lambda(s,b).
\end{equation}

We check that 
\begin{aenume}
\item $\displaystyle G(\beta_j) G(b_i,s_i,b'_i) \equiv \sum_{(s,b) \in
\B(\lambda)} \sum_{s'\in\Irr G_\lambda} c_{s s_i}^{s'}
g_{\beta_j,b_i}^{(s,b)} G(b,s',b'_i) \mod \Umgl,$
\item $\displaystyle G(b_i,s_i,b'_i) G(\beta_j)\equiv \sum_{(s,b)
\in \B(\lambda)} \sum_{s'\in\Irr G_\lambda} c_{s
s_i^\#}^{s^{\prime\#}} g_{\beta_j^\#,b'_i}^{(s,b)} G(b_i,s^{\prime},b)
\mod \Umgl,$
\end{aenume}
where $c_{\bullet\bullet}^\bullet$ are structure constants of
$(R(G_\lambda), \Irr G_\lambda)$ as before.

We consider the $\Um$--homomorphism $\Psi_{\tilde b_i}\colon
V(\lambda) \rightarrow \Uml$ such that $\Psi_{\tilde b_i}(u_\lambda)
\linebreak[0] = S_i G(b'_i)^\#\bmod \Umgl$, where $\tilde b_i =
(s_i^\#, b'_i)$ by \eqref{eq:ident}. (The existence of $\Psi_{\tilde
b_i}$ is guaranteed by \lemref{lem:extrem}.)  The image of
\eqref{uconsts} under $\Psi_{\tilde b_i}$ gives
\begin{equation*}
   G(\beta_j) G(b_i,s_i,b'_i) \equiv \sum_{(s,b)}
   g_{\beta_j,b_i}^{(s,b)} G(b) S S_i G(b'_i)^\#.
\end{equation*}
Now (a) follows. (b) follows from (a) by taking $\#$.

By (a) and (b) we have
\begin{align*}
  & c_{\be_1,\be_2}^\be(q_s)
  = \sum_{s'} \delta_{b'_2,b'} c_{s's_2}^{s} g_{\be_1, b_2}^{(s', b)}(q_s),
\\
  &  c_{\be,\be_3}^{\be_4}(q_s')
  = \sum_{s'} \delta_{b,b_4} c_{s' s^\#}^{s_4^\#}
               g_{\be_3^\#,b'}^{(s',b'_4)}(q_s'),
\\
  & c_{\be_1,\be}^{\be_4}(q_s)
  = \sum_{s'} \delta_{b'_4,b'} c_{s's}^{s_4}
    g_{\be_1, b}^{(s', b_4)}(q_s),
\\
  & c_{\be_2,\be_3}^{\be}(q_s') = 
  \sum_{s'} \delta_{b,b_2} c_{s' s_2^\#}^{s^\#}
    g_{\be_3^\#, b'_2}^{(s', b')}(q_s'),
\end{align*} 
where $\beta = G(b,s,b')$.
This makes the identity of the lemma equivalent to
\begin{multline*}
   \sum_{(b,s,b') \in \Bm[\lambda]}\sum_{s',s''\in\Irr G_\lambda}
   \delta_{b_2', b'} c_{s's_2}^{s} g_{\be_1, b_2}^{(s', b)}(q_s)
   \delta_{b,b_4} c_{s'' s^\#}^{s_4^\#} g_{\be_3^\#, b'}^{(s'', b'_4)}(q_s')
\\
   = \sum_{(b,s,b') \in \Bm[\lambda]}\sum_{s',s''\in\Irr G_\lambda}
   \delta_{b'_4, b'} c_{s's}^{s_4} g_{\be_1, b}^{(s', b_4)}(q_s) 
   \delta_{b,b_2} c_{s'' s_2^\#}^{s^\#} g_{\be_3^\#,
   b'_2}^{(s'',b')}(q_s').
\end{multline*}
Note that
\begin{equation*}
    \sum_s c_{s's_2}^{s} c_{s'' s^\#}^{s_4^\#}
    = \sum_s c_{s's_2}^{s} c_{s^{\prime\prime\#} s}^{s_4}
\quad\text{and}\quad
    \sum_s c_{s's}^{s_4}c_{s'' s_2^\#}^{s^\#}
    = \sum_s c_{s's}^{s_4}c_{s^{\prime\prime\#} s_2}^{s}
\end{equation*}
are equal, since both are the multiplicity of $s_4$ in $s'\otimes
s_2\otimes s^{\prime\prime\#}$. Now the above identity is immediate.
\end{proof}

Fix a two--sided cell $\Gml$. We define a $\Q(q_s)$--linear map
$\Phi\colon \Um \rightarrow \Q(q_s) \otimes \Uml_0$ by
$$
 \Phi(G(\beta)) = \sum_{d \in \mathcal D_{\Bml}, \beta' \in \Bml}
 c_{\beta,d}^{\beta'}(q_s) t_{\beta'}, \quad (\beta \in \Bm)
$$
which is well-defined since $\mathcal D_{\Bml}$ is finite and for a
fixed $\beta, d$ there are only finitely many $c_{\beta,d}^{\beta'}
\neq 0$.

We have the following description of $\Phi$ due to
\cite[Prop.~1.9]{Lu-v}
\begin{Proposition} 
\textup{(i)} $\Phi$ is an algebra homomorphism.

\textup{(ii)} Let $P(\Bml)$ be the set of $\lambda \in P$ such that
$a_\lambda d = d$ for some $d \in \mathcal D_{\Bml}$.  
Then $\Phi(\sum_{\la \in P(\Bml)} a_\lambda) = 1$, and $\Phi(a_{\lambda'})
= 0$ for $\lambda \notin P(\Bml)$.  \qed
\end{Proposition} 

Next we describe an explicit realization of the ring structure of
$\Uml_0$.  As usual, $\lambda = \sum_i \lambda_i \varpi_i$. 
Let $T_\lambda$ be the set of triples $(d, d', s)$ where
$d,d' \in D_{\Bml}$ and $s \in \operatorname{Irr}
G_\lambda$. Let $J_\lambda$ be the free abelian group on $T_\lambda$
with a ring structure defined by
$$
  (d_1,d'_1, s) (d_2,d'_2, s') = \delta_{d'_1,d_2} \sum_{s'' \in
  \operatorname{Irr}G_\lambda} c_{ss'}^{s''} (d_1, d'_2, s''),$$
where
$c_{ss'}^{s''}$ is the multiplicity of $s''$ in the tensor product $s
\otimes s'$ as above. We have:
\begin{Theorem}\label{thm:limit}
\textup{(i)} There exists a ring isomorphism 
$\Uml_0 \overset \sim \longrightarrow J_\lambda$ which gives a
bijection between $\{t_\beta \mid \beta \in \Bml \}$ and $\{(d,d',
s) \mid d, d' \in \mathcal D_{\Bml}, s \in
\operatorname{Irr} G_\lambda\}$.

\textup{(ii)} For any $d_0 \in \mathcal D_{\Bml}$, the subset of $\Bml$
corresponding to $\{(d,d', s) \in T_\lambda \mid d' = d_0\}$
under the bijection in \textup{(i)} is a left cell. 

\textup{(iii)} For any $d_0 \in \mathcal D_{\Bml}$, the subset of $\Bml$
corresponding to $\{(d,d', s) \in T_\lambda \mid d = d_0\}$
under the bijection in \textup{(i)} is a right cell.
\end{Theorem}

\begin{proof}  By definition,  the elements 
$b \in \B_W(\lambda)$ are in $1$--$1$ correspondence
with the elements of $t_b = t_{(b,1,b)} \in D_{\Bml}$. The map which
sends $t_{(b_1,s,b'_1)} \mapsto (t_{b_1}, t_{b'_1}, s)$ is a bijection.
Using \lemref{lem:mul2} as in the proof of \lemref{lem:P1} we have
\begin{equation} \label{tprod}  t_{(b_1,s,b_2)} t_{(b'_1, s', b'_2)} = 
\delta_{b_2,b'_1} \sum_{s_1\in\operatorname{Irr}G_\lambda}
   c_{s,s'}^{s_1} t_{(b_1,s_1,b'_2)}.  
\end{equation} 
This implies (i).  (ii) and (iii) follow from \propref{prop:twocell}.
\end{proof}

The following proposition partly explains why our approach to the
limit algebra $J_\lambda$ via $V(\lambda)$ is natural.

\begin{Proposition}
Let $V(\lambda)_0$ be the free $\Z$-module with basis $\B(\lambda)
\cong \mathcal D_{\Bml}\times\Irr G_\lambda$ endowed with a
$J_\lambda$-module structure by
\begin{equation*}
   (d_1, d_2, s) \cdot (d',s') 
   = \delta_{d_2,d'} \sum_{s'' \in  \Irr G_\lambda} c_{ss'}^{s''}
   (d_1, s'').
\end{equation*}
Then $V(\lambda)_0\otimes_\Z \Q(q_s)$, pulled back by the composition
of $\Um\xrightarrow{\Phi}\Uml_0\overset \sim \longrightarrow
J_\lambda$, is isomorphic to $V(\lambda)$.
\end{Proposition}

\begin{proof}
For $b = (b_W,s)\in \B(\lambda) \cong \B_W(\lambda)\times\Irr
G_\lambda$, we define $\beta_{b}\in \Bml$ by $\beta_{b}
= (b_W,s,u_\lambda)$. For $\beta\in\Bm$, we have
\begin{equation*}
   G(\beta) G_\lambda(b) = \sum_{b'\in \B(\lambda)}
   c_{\beta\beta_b}^{\beta_{b'}} G_\lambda(b').
\end{equation*}
In fact, $G_\lambda(b) = G(\beta_b)u_\lambda$ implies
\(
   G(\beta) G_\lambda(b) = \sum_{\beta'\in\Bm}
   c_{\beta\beta_b}^{\beta'} G(\beta')u_\lambda
\)
and $G(\beta')u_\lambda$ is $G_\lambda(\beta')$ if
$\beta'\in\B(\lambda)$ and $0$ otherwise. On the other hand, the
$\Uml_0$-module structure on $V(\lambda)_0$ is given by
\begin{equation*}
   t_{\beta'} b = \sum_{b'\in \B(\lambda)}
   \gamma_{\beta'\beta_{b}}^{\beta_{b'}} b'
   \qquad\text{for $\beta'\in \Bml$}.
\end{equation*}
We define a $\Q(q_s)$--linear map $\Psi\colon V(\lambda)\to
V(\lambda)_0\otimes_\Z \Q(q_s)$ by
\begin{equation*}
   \Psi(G_\lambda(b)) = \sum_{d \in \mathcal D_{\Bml}, b'\in \B(\lambda)}
   c_{\beta_b d}^{\beta_{b'}} b', \quad (b \in \B(\lambda)).
\end{equation*}
Then we have
\(
   \Psi(G_\lambda(b)) = \Phi(G(\beta_b)) u_\lambda,
\)
since $t_{\beta'} u_\lambda = b'$ if $\beta' = \beta_{b'}$ and $0$
otherwise by \thmref{thm:limit}. We get
\begin{equation*}
   \Psi(G(\beta) G_\lambda(b)) = \Phi(G(\beta)G(\beta_b)) u_\lambda
   = \Phi(G(\beta)) \Phi(G(\beta_b)) u_\lambda
   = \Phi(G(\beta)) \Psi(G_\lambda(b)).
\end{equation*}
This shows that $\Psi$ is $\Um$--linear.

Let $\La(\lambda)_0$ be the $\A_0$-submodule of
$V(\lambda)_0\otimes_\Z \Q(q_s)$ generated by $q^{-a(\beta_b)} b$
($b\in \B(\lambda)$). Note that $a(\beta_b)$ is independent of $b$ by
\lemref{lem:P1}(ii). Since
\(
   q^{a(\beta_{b'})} c_{\beta_b d}^{\beta_{b'}} =
   q^{a(\beta_b)} c_{\beta_b d}^{\beta_{b'}}\in \Z[q_s]
\)
by \lemref{lem:P1}(i), we have
$\Psi(\La(\lambda))\subset\La(\lambda)_0$. The induced homomorphism
\(
   \Psi_0\colon \La(\lambda)/q_s \La(\lambda)\to 
    \La(\lambda)_0/q_s\La(\lambda)_0
\)
is given by
\[
   \Psi_0(b) = \sum_{d \in \mathcal D_{\Bml}, b'\in \B(\lambda)}
   \gamma_{\beta_b d}^{\beta_{b'}} b',
\] 
where $b'\in\B(\lambda)$ is identified with
$q^{-a(\beta_{b'})} b'\bmod q_s\La(\lambda)_0$. By \thmref{thm:limit}
we have $\gamma_{\beta_b d}^{\beta_{b'}} =
\delta_{bb'}\delta_{\beta_b, \beta_b d}$, hence the right hand side is
equal to $b$. This shows that $\Psi$ is an isomorphism as in the proof 
of \corref{cor:inj}.
\end{proof}

\begin{Proposition} For $\Um(A^{(1)}_n)$ we have $\# \mathcal D_{\Bml} =
\prod_{i=1}^{n} \binom{n+1}{i}^{\lambda_i}$.
\end{Proposition}
\begin{proof}
  We check the dimension of $W(\lambda) = \bigotimes_i
W(\varpi_i)^{\otimes \lambda_i}.$ The result follows if we show the
dimension of $W(\varpi_i)$ is $\binom{n+1}{i}$.  By considering the
Drinfeld polynomials of $W(\varpi_i)$ \cite[Remark 3.3]{extrem}, each
is an evaluation module where the underlying finite dimensional
representation of $\Um(A_{n})$ is $V(\omega_i)$ where the $\omega_i$
are the fundamental weights of type $A_{n}$
\cite[Prop.~12.2.13]{CP-book}.  These are of dimension
$\binom{n+1}{i}$ for $i = 1, \dots, n$.
\end{proof}

\begin{Remark}
  The above argument shows that the number of $\mathcal D_{\Bml}$ can
  be given if we know $\dim W(\varpi_i)$ for all $i\in I_0$. They are
  known for untwisted affine Lie algebras for classical groups and
  some exceptional groups \cite{CP}. The second author writes a
  computer program for them (for untwisted cases) by using
  Frenkel-Mukhin's algorithm \cite{FM}. The program gives us answers
  except one fundamental representation for $E_8$, corresponding to
  the triple node.
\end{Remark}

\end{document}